\documentclass[11pt]{amsart}
\usepackage[T1]{fontenc}
\usepackage{lmodern,amsmath,amssymb,mathtools}
\usepackage[margin=1.12in]{geometry}
\usepackage[colorlinks=true,linkcolor=blue,citecolor=blue,urlcolor=blue]{hyperref}
\usepackage{microtype}
\newtheorem{theorem}{Theorem}[section]
\newtheorem{lemma}[theorem]{Lemma}
\newtheorem{proposition}[theorem]{Proposition}
\newtheorem{corollary}[theorem]{Corollary}
\theoremstyle{remark}
\newtheorem{remark}[theorem]{Remark}
\newtheorem*{remark*}{Remark}
\numberwithin{equation}{section}
\newcommand{\R}{\mathbb R}

\DeclareMathOperator{\ind}{ind}
\DeclareMathOperator{\diag}{diag}
\title[Relative Morse index and the number of multi-peak solutions]{Existence, relative Morse index and the number of positive multi-peak solutions to subcritical Hamiltonian systems}
\author{Qing Guo}
\address{College of Science, Minzu University of China, Beijing 100081, China}
\email{guoqing0117@163.com}
\thanks{Corresponding author: Qing Guo. This work was supported by the National Natural Science Foundation of China (Grant No.~12271539).}
\date{}
\subjclass[2020]{35J47, 35B25, 35B40, 35P15}
\keywords{Hamiltonian elliptic system, multi-peak solutions, relative Morse index, nondegeneracy, local uniqueness, Schur complement}
\begin{document}
\begin{abstract}
We construct positive multi-peak solutions of a singularly perturbed
Hamiltonian elliptic system with a common uniformly positive potential
$a$ and powers $p,q>1$ below the Sobolev critical hyperbola. The peaks
concentrate at any prescribed set of $k$ distinct nondegenerate critical
points $P_i$ of $a$.

The Hessian of the action is strongly indefinite. A Schur complement
gives exactly $kN$ eigenvalues near zero and their leading terms,
which have order $\varepsilon^2$ and are determined by $D^2a(P_i)$.
The solutions are nondegenerate. A scalar Schur reduction and
localization of its inverse term give the relative Morse index
$k+\sum_{i=1}^k\ind D^2a(P_i)$.

Local Pohozaev identities give local uniqueness. A compactness argument
then shows that, for $m$ prescribed nondegenerate critical points, the
union of the fixed profile neighborhoods associated with their nonempty
subsets contains exactly $2^m-1$ positive solutions. Using uniqueness of positive
solutions of the limiting system up to translations, we obtain the same
count in classes defined by nonvanishing in fixed balls and a scaled
nonlinear mass bound that excludes an additional peak. We also give
second-order expansions of the peak locations and rescaled profiles.
\end{abstract}
\maketitle
\enlargethispage{8pt}

\section{Introduction}

We consider the system
\begin{equation}\label{eq:system}
 \begin{cases}
 -\varepsilon^2\Delta u+a(x)u=v^p&\text{in }\R^N,\\
 -\varepsilon^2\Delta v+a(x)v=u^q&\text{in }\R^N,\\
 u,v>0,\qquad u(x),v(x)\longrightarrow0&\text{as }|x|\longrightarrow\infty,
 \end{cases}
\end{equation}
where $\varepsilon>0$ is small. We assume throughout that
\begin{equation}\label{eq:exponents}
 N\ge2,\qquad p,q>1,\qquad
 \frac1{p+1}+\frac1{q+1}>\frac{N-2}{N}.
\end{equation}
For $N=2$ the last condition is automatic. For the potential, we assume
\begin{equation}\label{eq:potential}
 a\in C_b^4(\R^N),\qquad 0<a_*\le a(x)\le a^*<\infty,
\end{equation}
where $C_b^4$ means that the function and its derivatives through order four are continuous and bounded. No limit of $a$ at infinity is assumed. Let $P_1,\ldots,P_k$ be distinct points satisfying
\begin{equation}\label{eq:critical-points}
 \nabla a(P_i)=0,\qquad \det D^2a(P_i)\ne0.
\end{equation}
The number $k$ is fixed. The points may be minima, maxima or saddles.

The scalar semiclassical problem motivates the present analysis. Floer and Weinstein \cite{FW1986} constructed concentrated states near nondegenerate critical points of the potential, and Oh \cite{Oh1990} obtained positive multi-peak states near prescribed finite collections of such points. Both constructions use the nondegeneracy of the autonomous profile to solve for a correction. The peak positions are then determined.

For Hamiltonian systems, Ramos and Soares \cite{RS2006} proved concentration of least-energy solutions at global minima under a gap condition at infinity. Ramos and Tavares \cite{RT2008} constructed multiple peaks near potential wells. Ramos \cite{Ramos2009} treated a prescribed critical point that need not be a minimum, in dimensions $3\le N\le6$. Manna and Sahoo \cite{MS2023} studied the effect of variable coefficients under Neumann boundary conditions. In dimension two, Zhang, Yang, Zhang and Zhong \cite{ZYZZ2024} obtained multi-peak solutions at strict local minima for nonlinearities with exponential growth. The survey \cite{BDT2014} discusses the variational frameworks and related concentration results.

There is also substantial work on critical Hamiltonian systems. Frank, Kim and Pistoia \cite{FKP2021} proved nondegeneracy of the autonomous ground state for the massless critical Lane--Emden system. Kim and Pistoia \cite{KP2021} constructed blowing-up solutions in bounded domains. Guo, Hu and Peng studied nondegeneracy under linear perturbations \cite{GHPlinear2024}, and local uniqueness and nondegeneracy for a critical problem with radial coefficients in a ball \cite{GHP2024}. In the whole space, Guo, Liu and Peng \cite{GLP2023} constructed critical multi-bubbling solutions and proved nondegeneracy in a prescribed symmetry space. The concentration scales and geometries in these works are different from ours.

In this paper we construct positive multi-peak solutions of
\eqref{eq:system} and study their local uniqueness, small eigenvalues and
relative Morse index. The peaks may concentrate at minima, maxima or
saddles of $a$. The construction uses a positive ground state of the
limiting system and its nondegeneracy up to translations, recalled in
Lemma~\ref{lem:autonomous-nondegeneracy}. Uniqueness of all positive
solutions of the limiting system is used only in
Section~\ref{sec:positive-class}.

We solve the projected equation in $(H^2\cap W^{2,t})^2$, where
$t>\max\{N,2\}$, and then choose the centers by Brouwer degree.
One exponent may exceed the scalar Sobolev threshold. Also, for
$1<p<2$ or $1<q<2$, the corresponding power is not twice differentiable
at zero. The nonlinear estimates use only its first derivative as a map
from $H^2\cap W^{2,t}$ to $L^2\cap L^t$.

To prove local uniqueness, we compare two arbitrary solutions with the
prescribed profiles. A local Pohozaev identity first gives the center
estimate $P_i+o(\varepsilon)$. After rescaling at $P_i$, every limit of
the normalized difference belongs to the translation kernel of the
limiting system. Subtracting the two Pohozaev identities gives a linear
system with matrix $D^2a(P_i)$ for the translation coefficients. These
coefficients vanish, and an exterior estimate gives uniqueness.
This argument uses neither the small-spectrum theorem nor the relative
Morse index. Local uniqueness near prescribed scalar multi-peak configurations was
studied in \cite{CH2003}. Cao, Li and Luo \cite{CLL2015} used local
Pohozaev identities under local expansion conditions on the potential.
Related counting results can be found in \cite{Grossi2002,CLP2021}.

The spectral and index questions require a different argument. The
linearized operator of \eqref{eq:system} is in general not self-adjoint.
Exchanging its rows gives the Hessian of the action, which has infinitely
many positive and negative directions. Removing the approximate
translation modes still leaves an indefinite quadratic form. We prove
an inverse estimate on this complement. A Schur complement then reduces
the eigenvalue equation near zero to a symmetric $kN\times kN$ matrix.
Its dependence on the spectral parameter gives all $kN$ small
eigenvalues, including their multiplicities and their leading terms.
For related spectral reduction methods in scalar problems, see
\cite{BS2002}.

Finite indices for strongly indefinite elliptic systems were studied by
Angenent and van der Vorst \cite{AV2000} and Ramos
\cite{RamosIndex2009}. See also \cite{BDT2014,Amadori2024}.
Here the ordinary Morse index is infinite. We use the sum and difference
of the two components to separate a negative definite block. Completing
the square leaves a scalar Schur operator $\mathcal K_\varepsilon$,
bounded below, whose negative index equals the relative Morse index.
For unequal powers this operator contains a nonlocal inverse term.
We localize that term by an inverse-operator identity in which only the
cutoffs are differentiated. The variational characterization of the
ground state gives one negative direction per peak. Together with the
negative small eigenvalues, this yields
\[
 m_{\mathrm{rel}}(z_\varepsilon)
 =k+\sum_{i=1}^k\ind D^2a(P_i).
\]
Thus the formula separates the contribution of the limiting ground
states from that of the concentration points. For example, $k$ peaks at
nondegenerate minima have relative index $k$, whereas $k$ peaks at
nondegenerate maxima have relative index $k(N+1)$.

Finally, we count the solutions in fixed profile neighborhoods. Every
solution in such a neighborhood converges to the prescribed profiles.
The direct uniqueness theorem then gives exactly one solution for each
nonempty subset of the selected critical points. The neighborhoods are
pairwise disjoint, so their union contains exactly $2^m-1$ solutions.
In Section~\ref{sec:positive-class} we replace profile closeness by
nonvanishing in fixed balls and a bound on the scaled nonlinear mass.
The mass bound excludes additional peaks, including peaks escaping to
infinity, and uniqueness for the limiting system identifies the profiles.
The same count follows in these classes. The perturbed solutions
are not assumed to minimize the action.

\subsection{Main conclusions}

We first recall some facts about the limiting system
\begin{equation}\label{eq:autonomous}
 -\Delta U+U=V^p,\qquad -\Delta V+V=U^q \quad\text{in }\R^N.
\end{equation}
Let $W=(U,V)\in H^1(\R^N)^2$ be a positive ground state of
\eqref{eq:autonomous}, centered at the origin. The existence of a positive
ground state and its radial symmetry follow from
\cite[Theorems~1.8--1.9]{BMR2012} and \cite[Theorem~6.4]{BDT2014}.
The two components are smooth and strictly decreasing in the radial
variable. For some $C,\sigma>0$,
\[
 U(y)+V(y)+|\nabla U(y)|+|\nabla V(y)|\le Ce^{-\sigma|y|};
\]
see also \cite[Lemma~3.1]{CR2026}. Set $r_p=(p+1)/p$ and
$\mathcal R_1=(-\Delta+1)^{-1}$. The variational characterization of $W$
gives
\begin{equation}\label{eq:ground-state-quotient}
 \frac{\|\mathcal R_1(V^p)\|_{q+1}}{\|V^p\|_{r_p}}
 =\sup_{0\ne f\in L^{r_p}(\R^N)}
          \frac{\|\mathcal R_1 f\|_{q+1}}{\|f\|_{r_p}}.
\end{equation}
In \cite[formula~(6.5)]{BDT2014}, the exponents appear in the opposite order. In our notation, $U$ minimizes
$\|(-\Delta+1)w\|_{r_p}^{r_p}/\|w\|_{q+1}^{r_p}$.
Taking $w=\mathcal R_1 f$ then gives \eqref{eq:ground-state-quotient}.
We use this ground state throughout the paper.

For $b>0$, define the scaled profiles by
\begin{equation}\label{eq:profile-scaling}
 U_b(y)=b^{(p+1)/(pq-1)}U(\sqrt b\,y),\qquad
 V_b(y)=b^{(q+1)/(pq-1)}V(\sqrt b\,y).
\end{equation}
Here $b_i=a(P_i)$, $U_i=U_{b_i}$ and $V_i=V_{b_i}$. Fix $t>\max\{N,2\}$, and use the scaled norm
\begin{equation}\label{eq:scaled-norm}
 \|f\|_{W^{2,s}_\varepsilon}
 =\varepsilon^{-N/s}\bigl(\|f\|_s+\varepsilon\|\nabla f\|_s
                   +\varepsilon^2\|D^2f\|_s\bigr),\qquad
 \|z\|_{X_\varepsilon}=\|z\|_{H^2_\varepsilon}+\|z\|_{W^{2,t}_\varepsilon}.
\end{equation}
Here $H^2_\varepsilon=W^{2,2}_\varepsilon$, and pair norms are sums of
the component norms. After the change of variables $x=\varepsilon y$,
$X_\varepsilon$ is the norm of $(H^2\cap W^{2,t})^2$.

The following theorem gives existence and local uniqueness.
\begin{theorem}\label{thm:main-existence}
Assume \eqref{eq:exponents}--\eqref{eq:critical-points}. For all sufficiently small $\varepsilon>0$, system \eqref{eq:system} has a solution $z_\varepsilon=(u_\varepsilon,v_\varepsilon)$ of the form
\begin{equation}\label{eq:main-decomposition}
 z_\varepsilon(x)=\sum_{i=1}^k
 (U_{a(\xi_{i,\varepsilon})},V_{a(\xi_{i,\varepsilon})})
 \left(\frac{x-\xi_{i,\varepsilon}}{\varepsilon}\right)+w_\varepsilon(x),
 \qquad
 |\xi_{i,\varepsilon}-P_i|+\|w_\varepsilon\|_{X_\varepsilon}=O(\varepsilon^2).
\end{equation}
Moreover, let $\widetilde z_\varepsilon
=(\widetilde u_\varepsilon,\widetilde v_\varepsilon)$ be any family of
positive classical solutions of \eqref{eq:system}, decaying at infinity,
such that
\begin{equation}\label{eq:main-admissible}
 \left\|\widetilde z_\varepsilon-
 \sum_i(U_i,V_i)\left(\frac{\cdot-\eta_{i,\varepsilon}}{\varepsilon}\right)
 \right\|_\infty\longrightarrow0,
 \qquad \eta_{i,\varepsilon}\longrightarrow P_i.
\end{equation}
Then $\widetilde z_\varepsilon=z_\varepsilon$ for all sufficiently
small values of $\varepsilon$ for which $\widetilde z_\varepsilon$ is defined.
\end{theorem}

The convergence in \eqref{eq:main-admissible} is in the whole space.
Section~\ref{sec:construction} proves existence and
\eqref{eq:main-decomposition}. The uniqueness assertion is proved in
Section~\ref{sec:concentration-estimates} by comparing two solutions with
the local Pohozaev identity.

Define
\begin{equation}\label{eq:main-linear}
 \mathcal L_\varepsilon=
 \begin{pmatrix}
 -\varepsilon^2\Delta+a&-p v_\varepsilon^{p-1}\\
 -q u_\varepsilon^{q-1}&-\varepsilon^2\Delta+a
 \end{pmatrix},\qquad
 \mathcal H_\varepsilon=S\mathcal L_\varepsilon,\qquad S(h,\ell)=(\ell,h),
\end{equation}
Both operators have domain $H^2(\R^N)^2$ in $L^2(\R^N)^2$. The associated constants are defined by
\begin{equation}\label{eq:main-kappa}
 M_i=\int U_iV_i,\qquad
 d_i=\int\bigl((\partial_1U_i)^2+(\partial_1V_i)^2\bigr),\qquad
 \kappa_i=M_i/d_i>0.
\end{equation}
If no domain is indicated, the integral is taken over $\R^N$.

The next theorem concerns the self-adjoint Hessian
$\mathcal H_\varepsilon$.
\begin{theorem}\label{thm:main-spectrum}
Assume \eqref{eq:exponents}--\eqref{eq:critical-points}. For the constructed solution, and more generally for any positive solution family satisfying \eqref{eq:main-admissible}, there exists $\gamma>0$ such that, for small $\varepsilon$, $\mathcal H_\varepsilon$ has exactly $kN$ eigenvalues in $(-\gamma,\gamma)$, counted with multiplicity. If $\mu_1,\ldots,\mu_{kN}$ are the eigenvalues of
\[
 \diag\bigl(\kappa_iD^2a(P_i)\bigr)_{i=1}^k,
\]
then these eigenvalues can be numbered so that
\begin{equation}\label{sp:asymptotic}
 \lambda_{\varepsilon,j}=\varepsilon^2\mu_j+o(\varepsilon^2).
\end{equation}
Consequently $\ker\mathcal L_\varepsilon=\{0\}$, and the number of negative eigenvalues of $\mathcal H_\varepsilon$ in this interval equals
$\sum_i\ind(D^2a(P_i))$.
\end{theorem}

The full Hessian is still strongly indefinite, and this finite count only concerns the eigenvalues near zero. The next theorem gives the further information needed for a relative index.

Set $A_\varepsilon=-\varepsilon^2\Delta+a(x)$, acting on each component, and
$B_\varepsilon=\diag(qu_\varepsilon^{q-1},pv_\varepsilon^{p-1})$.
The free reference is $\mathcal H_{0,\varepsilon}=S A_\varepsilon$.
The relative Morse index $m_{\mathrm{rel}}(z_\varepsilon)$ counts, with multiplicity, the crossings from positive to negative through zero along the self-adjoint path on the fixed domain $H^2(\mathbb R^N)^2$
\[
 \mathcal H_{\varepsilon,t}=S A_\varepsilon-t B_\varepsilon,\qquad 0\le t\le1.
\]
Zero directions at $t=1$ are recorded separately by
$\dim\ker\mathcal H_\varepsilon$. Section~\ref{sec:relative-morse}
proves that the crossings are downward and that their total multiplicity
is finite. For the scalar Schur operator defined in
\eqref{mi:scalar-Schur}, we have
\[
 m_{\mathrm{rel}}(z_\varepsilon)=n_-(\mathcal K_\varepsilon),
 \qquad
 \dim\ker\mathcal H_\varepsilon=\dim\ker\mathcal K_\varepsilon.
\]
Here $n_-$ counts strictly negative eigenvalues with multiplicity.

\begin{theorem}\label{thm:main-morse}
Under \eqref{eq:exponents}--\eqref{eq:critical-points}, every positive solution family satisfying \eqref{eq:main-admissible} for the selected profiles satisfies
\begin{equation}\label{eq:main-morse}
 m_{\mathrm{rel}}(z_\varepsilon)
 =k+\sum_{i=1}^k\ind D^2a(P_i),\qquad
 \dim\ker\mathcal H_\varepsilon=0
\end{equation}
for all sufficiently small $\varepsilon$. In particular, the formula holds for the solutions of Theorem~\ref{thm:main-existence}.
\end{theorem}

For the counting result, consider $m$ distinct nondegenerate critical
points $P_1,\ldots,P_m$, with $W_j=(U_{a(P_j)},V_{a(P_j)})$ and
$\iota_j=\ind D^2a(P_j)$.  For a nonempty subset
$J\subset\{1,\ldots,m\}$, a neighborhood in $C_0(\R^N)^2$ is defined by
\begin{equation}\label{eq:count-tube}
 \mathcal U_{\varepsilon,J}(\delta,\rho)
 =\left\{z\in C_0(\R^N)^2:
 \inf_{\xi_j\in B_\rho(P_j),\ j\in J}
 \left\|z-\sum_{j\in J}W_j
       \left(\frac{\cdot-\xi_j}{\varepsilon}\right)\right\|_\infty
 <\delta\right\}.
\end{equation}
Here $C_0$ consists of continuous functions tending to zero at infinity,
and the norm of a pair is the sum of the component supremum norms.
The radius $\delta$ is fixed. Membership in this neighborhood does not assume an error
tending to zero with $\varepsilon$.

\begin{theorem}\label{thm:main-count}
Under \eqref{eq:exponents} and \eqref{eq:potential}, there exist
$\rho,\delta,\varepsilon_0>0$ such that, for every
$0<\varepsilon<\varepsilon_0$, each nonempty subset $J$ gives exactly
one positive classical solution $z_{J,\varepsilon}$ in
$\mathcal U_{\varepsilon,J}(\delta,\rho)$.
These neighborhoods are pairwise disjoint.
Let $\mathcal P_\varepsilon$ be the set of positive classical solutions
in their union.  Then $\#\mathcal P_\varepsilon=2^m-1$, and exactly
$\binom{m}{k}$ of these solutions have $k$ local maximum points of $u+v$.
Every solution is nondegenerate, and
\[
 m_{\mathrm{rel}}(z_{J,\varepsilon})
      =\sum_{j\in J}(1+\iota_j).
\]
\end{theorem}

Section~\ref{sec:counting} proves this theorem by compactness and the
direct uniqueness result. In the next theorem we use uniqueness of
positive solutions of \eqref{eq:autonomous} up to translations to replace
profile closeness by a mass bound and nonvanishing in fixed balls.

The scaled nonlinear mass of a positive classical decaying solution $z=(u,v)$ is defined by
\begin{equation}\label{eq:positive-mass}
 \begin{aligned}
 \mathcal Q_\varepsilon(z)
 &=\varepsilon^{-N}\int u^{q+1}
  =\varepsilon^{-N}\int v^{p+1},\\
 m(b)&=\int U_b^{q+1}=\int V_b^{p+1},\qquad
 m_* =\min_{b\in[a_*,a^*]}m(b)>0.
 \end{aligned}
\end{equation}
The integrals are finite, and the two component identities hold for every
such solution. These facts are proved in Section~\ref{sec:positive-class}. For solutions, the identity
$\mathcal Q_\varepsilon=\varepsilon^{-N}I_\varepsilon/\theta$ holds, where
$I_\varepsilon$ is the action in \eqref{eq:energy} and
$\theta=1-1/(p+1)-1/(q+1)>0$. The condition below is therefore also an energy
bound. The positive function $m(b)$ is continuous by
\eqref{eq:profile-scaling}.

The same $m$ points as in Theorem~\ref{thm:main-count} are used here. Fix
\[
 0<\kappa<m_*,\qquad
 0<\delta<\frac12\min_{1\le i\le m}
                 \bigl(U_{a(P_i)}(0)+V_{a(P_i)}(0)\bigr).
\]
For a radius $\rho>0$ and a nonempty subset $J$, let
\begin{equation}\label{eq:positive-class}
 \begin{split}
 \mathcal A_{\varepsilon,J}(\rho,\delta,\kappa)
 =\biggl\{z=(u,v):\;&z\text{ is a positive classical solution of
                         \eqref{eq:system}},\\
 &\max_{\overline{B_\rho(P_i)}}(u+v)\ge\delta
                                      \quad(i\in J),\\
 &\mathcal Q_\varepsilon(z)
          \le\sum_{i\in J}m(a(P_i))+m_*-\kappa\biggr\}.
 \end{split}
\end{equation}
The balls have a fixed physical radius. This definition does not include any profile approximation,
minimizing property or exact limiting mass.
It does not exclude additional peaks either; their exclusion is
part of the proof.

\begin{theorem}\label{thm:main-positive}
Assume \eqref{eq:exponents} and \eqref{eq:potential}, and suppose that
$P_1,\ldots,P_m$ are distinct nondegenerate critical points of $a$.
For every $\kappa,\delta$ in the ranges above, there exist
$\rho,\varepsilon_0>0$ such that, for all
$0<\varepsilon<\varepsilon_0$ and every nonempty
$J\subset\{1,\ldots,m\}$,
\[
 \mathcal A_{\varepsilon,J}(\rho,\delta,\kappa)
     =\{z_{J,\varepsilon}\},
\]
where $z_{J,\varepsilon}$ is the solution in
Theorem~\ref{thm:main-count}. Consequently these classes are pairwise
disjoint, and their union contains exactly $2^m-1$ solutions, with the
peak numbers and relative Morse indices stated there. The spectral
conclusions of Theorem~\ref{thm:main-spectrum} and the concentration
conclusions of Section~\ref{sec:quantitative-concentration} also hold.
Moreover,
\[
 \mathcal Q_\varepsilon(z_{J,\varepsilon})
       \longrightarrow\sum_{i\in J}m(a(P_i)).
\]
\end{theorem}

We also obtain a second-order location formula. With
\begin{equation}\label{eq:main-rho}
 \rho_i=\frac{\int |y|^2U_i(y)V_i(y)\,dy}{2NM_i},
\end{equation}
for the choice of centers made in Section~\ref{sec:quantitative-concentration}, we have
\begin{equation}\label{eq:intro-shift}
 \xi_{i,\varepsilon}=P_i-\varepsilon^2\rho_i
 [D^2a(P_i)]^{-1}\nabla\Delta a(P_i)+o(\varepsilon^2).
\end{equation}
Section~\ref{sec:quantitative-concentration} gives the even second-order
profile correction. The two component maxima in the $i$th core lie at
$\xi_{i,\varepsilon}+o(\varepsilon^3)$, and
$u_\varepsilon+v_\varepsilon$ has exactly $k$ local maximum points.
We also give the mass and energy concentration limits.

Section~\ref{sec:profiles} recalls the limiting profiles and the
functional setting. Section~\ref{sec:construction} gives the construction,
and Section~\ref{sec:concentration-estimates} proves local uniqueness.
The small spectrum and the relative Morse index are studied in
Sections~\ref{sec:potential-spectrum} and \ref{sec:relative-morse}.
Section~\ref{sec:quantitative-concentration} gives the location and profile
expansions. The counting results are proved in
Sections~\ref{sec:counting} and \ref{sec:positive-class}.

\section{Autonomous profiles and the functional setting}\label{sec:profiles}

\subsection{Basic properties and uniform estimates for the profiles}

For $b>0$, let $(U_b,V_b)$ be given by \eqref{eq:profile-scaling}. Write
\[
 \alpha=\frac{p+1}{pq-1},\qquad \beta=\frac{q+1}{pq-1}.
\]
Since $\alpha+1=p\beta$ and $\beta+1=q\alpha$, we have
\begin{equation}\label{eq:frozen-equations}
 (-\Delta+b)U_b=V_b^p,\qquad (-\Delta+b)V_b=U_b^q.
\end{equation}
Set $\mathcal R_b=(-\Delta+b)^{-1}$. The function $f_b=V_b^p$ attains
\begin{equation}\label{eq:profile-best-quotient}
 \sup_{0\ne f\in L^{r_p}(\mathbb R^N)}
 \frac{\|\mathcal R_b f\|_{q+1}}{\|f\|_{r_p}}.
\end{equation}
Indeed, for $f(y)=h(\sqrt b\,y)$,
\[
 \mathcal R_b f(y)=b^{-1}(\mathcal R_1h)(\sqrt b\,y),\qquad
 \frac{\|\mathcal R_b f\|_{q+1}}{\|f\|_{r_p}}
 =b^{-1+\frac N2(\frac1{r_p}-\frac1{q+1})}
   \frac{\|\mathcal R_1h\|_{q+1}}{\|h\|_{r_p}}.
\]
The change $h\mapsto f$ is a bijection on $L^{r_p}$, and $V_b^p$ is
a positive multiple of $V^p(\sqrt b\,\cdot)$. Thus
\eqref{eq:ground-state-quotient} gives the assertion.

\begin{lemma}[{\cite[Theorem~1.2(2)]{CR2026}}]\label{lem:autonomous-nondegeneracy}
Assume \eqref{eq:exponents}. For each $b>0$, let
\[
 L_b=\begin{pmatrix}
 -\Delta+b&-pV_b^{p-1}\\
 -qU_b^{q-1}&-\Delta+b
 \end{pmatrix},
 \qquad Z_{b,j}=\partial_j(U_b,V_b).
\]
Then the kernel of $L_b$ in $H^2(\R^N)^2$ is
\begin{equation}\label{eq:profile-kernel}
 \ker L_b=\operatorname{span}\{Z_{b,j}:1\le j\le N\}.
\end{equation}
\end{lemma}

\begin{proof}
For $b=1$, put $(\bar p,\bar q)=(p+1,q+1)$ in
\cite[Theorem~1.2(2)]{CR2026}. The condition
$N\ge2$, $\bar p,\bar q>2$ and $1/\bar p+1/\bar q>(N-2)/N$ then
becomes \eqref{eq:exponents}, and the linearized operator is $L_1$.
For general $b>0$, the change of variables in \eqref{eq:profile-scaling}
maps the kernel of $L_1$ onto that of $L_b$ and sends each translation
mode to a nonzero multiple of the corresponding translation mode.
\end{proof}

In the following estimates, $b$ varies in a fixed compact interval
$J\subset(0,\infty)$ containing the range of $a$.

The following estimates justify the polynomial moments, parameter dependence and translation derivatives used later. They follow from the decay of $(U,V)$ and the equations \eqref{eq:frozen-equations}.

\begin{lemma}\label{lem:profile-estimates}
There are $C,\sigma>0$, depending on $J$, such that
\begin{equation}\label{eq:profile-pointwise-decay}
 \sum_{j=0}^{2}\bigl(|D^jU_b(y)|+|D^jV_b(y)|\bigr)
 \le Ce^{-\sigma|y|}\qquad(b\in J).
\end{equation}
For every finite $r\ge2$ and integer $m\ge0$,
\begin{equation}\label{eq:profile-weighted-derivatives}
 \sup_{b\in J}\sum_{j=0}^{3}
 \bigl(\|\langle y\rangle^mD^jU_b\|_r+
       \|\langle y\rangle^mD^jV_b\|_r\bigr)<\infty,
 \qquad \langle y\rangle=(1+|y|^2)^{1/2}.
\end{equation}
For the fixed $t>\max\{N,2\}$ in the introduction, $b\mapsto(U_b,V_b)$ is $C^1$ with values in $(H^2\cap W^{2,t})^2$. Its derivative is
\begin{equation}\label{eq:profile-mass-derivative}
 \partial_bU_b=b^{-1}(\alpha U_b+\tfrac12 y\cdot\nabla U_b),\qquad
 \partial_bV_b=b^{-1}(\beta V_b+\tfrac12 y\cdot\nabla V_b).
\end{equation}
The matrix of inner products of the translations is
\begin{equation}\label{eq:profile-Gram}
 \langle Z_{b,j},Z_{b,l}\rangle=d_b\delta_{jl},\qquad
 d_b=\frac1N\int(|\nabla U_b|^2+|\nabla V_b|^2)>0.
\end{equation}
The numbers $d_b$ and $M_b:=\int U_bV_b$ have positive lower bounds on $J$.
\end{lemma}

\begin{proof}
By scaling, the profiles and their first derivatives decay exponentially, uniformly for $b\in J$. Differentiating \eqref{eq:frozen-equations} gives
\[
 (-\Delta+b)\partial_lU_b=pV_b^{p-1}\partial_lV_b,
 \qquad
 (-\Delta+b)\partial_lV_b=qU_b^{q-1}\partial_lU_b.
\]
Interior estimates therefore give
\[
 \|\partial_lU_b\|_{W^{2,r}(B_1(y_0))}
 +\|\partial_lV_b\|_{W^{2,r}(B_1(y_0))}
 \le C_r e^{-\sigma|y_0|},\qquad 2\le r<\infty.
\]
Taking $r>N$ proves \eqref{eq:profile-pointwise-decay}. Summation over a unit-ball covering with bounded overlap proves \eqref{eq:profile-weighted-derivatives}.

Differentiating \eqref{eq:profile-scaling} gives \eqref{eq:profile-mass-derivative}. The weighted estimates control its derivatives through order two in $L^2\cap L^t$. Local continuity under scaling and uniform control of the tails give continuity in these norms. This proves the asserted $C^1$ dependence. The same argument applies to translations.

Radial symmetry gives \eqref{eq:profile-Gram}. Positivity and continuity on the compact interval $J$ give the lower bounds for $d_b$ and $M_b$.
\end{proof}

\subsection{The adjoint kernel and solvability}

The next lemma identifies the range of the autonomous linearization. Its adjoint kernel determines the source terms in the projected equation and the solvability condition for the second-order profile correction.

\begin{lemma}\label{lem:frozen-range}
As an operator from $H^2(\R^N)^2$ to $L^2(\R^N)^2$, $L_b$ is Fredholm of index zero, and
\begin{equation}\label{eq:adjoint-kernel}
 \ker L_b^*=\operatorname{span}\{SZ_{b,j}:1\le j\le N\},\qquad
 \operatorname{Ran}L_b=\{h:\langle h,SZ_{b,j}\rangle=0\text{ for all }j\}.
\end{equation}
For every $h$ in this range there is a unique solution $\zeta\in H^2(\R^N)^2$ satisfying
\begin{equation}\label{eq:frozen-phase}
 L_b\zeta=h,\qquad \langle\zeta,Z_{b,j}\rangle=0\quad(1\le j\le N).
\end{equation}
Uniformly for $b\in J$, it satisfies $\|\zeta\|_{H^2}\le C\|h\|_2$. If also $h\in L^t$, then $\zeta\in W^{2,t}$ and
\begin{equation}\label{eq:frozen-X-inverse}
 \|\zeta\|_{H^2}+\|\zeta\|_{W^{2,t}}
 \le C(\|h\|_2+\|h\|_t).
\end{equation}
\end{lemma}

\begin{proof}
Put $A_b=-\Delta+b$. The bounded coefficients of $L_b-A_bI$ tend uniformly to zero at infinity. Rellich compactness on balls and smallness outside them show that this multiplication operator is compact from $H^2$ to $L^2$. Thus $L_b$ is Fredholm of index zero. Moreover,
\[
 L_b^*=SL_bS,\qquad D(L_b^*)=H^2(\R^N)^2,
\]
since the multiplication terms are bounded on $L^2$. The kernel identity \eqref{eq:profile-kernel} and the Fredholm alternative give \eqref{eq:adjoint-kernel}. The moments can be imposed because the translation Gram matrix is invertible. They also ensure uniqueness.

For uniformity, suppose that $b_n\to b\in J$, the moments of $\zeta_n$ vanish, $\|\zeta_n\|_{H^2}=1$, and $L_{b_n}\zeta_n\to0$ in $L^2$. Any weak limit lies in $\ker L_b$ and has zero moments; the latter pass to the limit by local compactness and the uniform $L^2$ tails of $Z_{b_n,j}$. Hence the limit is zero. Local compactness and uniform coefficient decay then give
\[
 \|(L_{b_n}-A_{b_n}I)\zeta_n\|_2\to0.
\]
The uniform bound for $A_{b_n}^{-1}:L^2\to H^2$ contradicts the normalization.

Finally, interpolate $h$ between $L^2$ and $L^t$ and apply
\[
 \|\zeta\|_{W^{2,r}}\le C_r(\|h\|_r+\|\zeta\|_r).
\]
Starting at $r=2$, Sobolev embedding with successive exponents satisfying $1/r'=\max\{1/t,1/r-1/N\}$ reaches $t$ in finitely many steps and proves \eqref{eq:frozen-X-inverse}.
\end{proof}

For later reference, we define the action on pairs for which it is finite by
\begin{equation}\label{eq:energy}
 I_\varepsilon(u,v)=\int\left(\varepsilon^2\nabla u\cdot\nabla v+a(x)uv
       -\frac{|u|^{q+1}}{q+1}-\frac{|v|^{p+1}}{p+1}\right).
\end{equation}
Here $I_b$ is the action with unit gradient coefficient and constant potential $b$. Multiplying the first autonomous equation by $V_b$ and the second by $U_b$, integrating, and using the decay gives
\begin{equation}\label{eq:profile-action}
 \begin{aligned}
 m(b)&=\int U_b^{q+1}=\int V_b^{p+1}
 =\int(\nabla U_b\cdot\nabla V_b+bU_bV_b),\\
 I_b(U_b,V_b)&=\left(1-\frac1{p+1}-\frac1{q+1}\right)m(b)>0.
 \end{aligned}
\end{equation}
The decay in Lemma~\ref{lem:profile-estimates} justifies these integrations.

\subsection{Spaces and the nonlinear remainder}

In the variable $y=x/\varepsilon$, set
\begin{equation}\label{eq:XY}
 \begin{gathered}
 X=(H^2\cap W^{2,t})(\R^N)^2,\qquad
 Y=(L^2\cap L^t)(\R^N)^2,\qquad t>\max\{N,2\}.
 \end{gathered}
\end{equation}
The intersections have sum norms. We use the continuous embedding $X\hookrightarrow(L^\infty\cap L^2\cap L^t)^2$. For $f_\varepsilon(y)=f(\varepsilon y)$, direct changes of variables give
\[
 \|D^j f_\varepsilon\|_s
   =\varepsilon^{j-N/s}\|D^j f\|_s\quad(j=0,1,2).
\]
Thus \eqref{eq:scaled-norm} is the usual Sobolev norm after scaling. In particular,
\begin{equation}\label{eq:scaled-embedding}
 \|z\|_\infty+\varepsilon\|\nabla z\|_\infty
 \le C\|z\|_{X_\varepsilon},
\end{equation}
with $C$ independent of $\varepsilon$. By the compact Sobolev embeddings for $t>N$, bounded sequences in $X$ have subsequences converging in $(C^1_{\rm loc})^2$ and in $(L^2_{\rm loc}\cap L^t_{\rm loc})^2$.

Every component of a function in $X$ tends to zero at infinity. If this were false, uniform continuity would give disjoint balls of a fixed radius on which its absolute value has a fixed positive lower bound, contradicting $L^2$ integrability. This fact will justify the maximum principle.

The nonlinear estimate below is the only global smoothness estimate on powers needed in the construction.

\begin{lemma}\label{lem:power-map}
For $s>1$, the map $P_s:f\mapsto(f_+)^s$ is $C^1$ from $H^2\cap W^{2,t}$ to $L^2\cap L^t$. Its derivative is $DP_s(f)h=s(f_+)^{s-1}h$. On a fixed bounded set in the domain,
\begin{equation}\label{eq:power-derivative-modulus}
 \|DP_s(f)-DP_s(g)\|_{\mathcal B(H^2\cap W^{2,t},L^2\cap L^t)}
 \le C\|f-g\|_{H^2\cap W^{2,t}}^{\min\{1,s-1\}}.
\end{equation}
Let $\mathcal N(u,v)=((v_+)^p,(u_+)^q)$, $\tau=\min\{1,p-1,q-1\}$, and
\[
 \mathcal R_W(w)=\mathcal N(W+w)-\mathcal N(W)-D\mathcal N(W)w.
\]
For uniformly bounded $W$ and $w,\widetilde w$ in a fixed small $X$ ball,
\begin{equation}\label{eq:nonlinear-estimate}
 \|\mathcal R_W(w)-\mathcal R_W(\widetilde w)\|_Y
 \le C(\|w\|_X^\tau+\|\widetilde w\|_X^\tau)
                         \|w-\widetilde w\|_X.
\end{equation}
The constant depends on a uniform $L^\infty$ bound for $W$, and not on the locations of its peaks.
\end{lemma}

\begin{proof}
For $r=2,t$, the bound $\|(f_+)^s\|_r\le\|f\|_\infty^{s-1}\|f\|_r$ shows that $P_s$ is well defined. Put $f_s(a)=(a_+)^s$ and $\tau_s=\min\{1,s-1\}$. On bounded intervals,
\[
 |f_s'(a)-f_s'(b)|\le C|a-b|^{\tau_s}.
\]
This follows from concavity when $1<s<2$ and the mean value theorem when $s\ge2$. Taking the coefficient in $L^\infty$ and the test function in $L^2\cap L^t$ proves \eqref{eq:power-derivative-modulus}. The identity
\[
 P_s(f+h)-P_s(f)-f_s'(f)h
 =\int_0^1[f_s'(f+\theta h)-f_s'(f)]h\,d\theta
\]
then proves Fr\'echet differentiability and continuity of the derivative. Applying the same estimate to
\[
 \mathcal R_W(w)-\mathcal R_W(\widetilde w)
 =\int_0^1[D\mathcal N(W+\widetilde w+\theta(w-\widetilde w))
                -D\mathcal N(W)](w-\widetilde w)\,d\theta
\]
gives \eqref{eq:nonlinear-estimate}, uniformly for bounded $W$.
\end{proof}

\section{Construction near prescribed critical points}\label{sec:construction}

We first solve the projected equation at fixed centers, using an inverse with finitely many auxiliary variables and moment conditions. We then choose the centers so that the auxiliary terms vanish. The local translation identity gives the position equation, and its continuous expansion yields existence by Brouwer degree. The argument also requires decay of the projected solution before its auxiliary terms vanish; we prove this decay below.

Choose $d>0$ so small that the closed balls $\overline{B_{4d}(P_i)}$ are
pairwise disjoint, and put
\[
 \mathcal D=\prod_{i=1}^k B_d(P_i),\qquad
 \xi=(\xi_1,\ldots,\xi_k)\in\overline{\mathcal D},\qquad
 b_{i,\xi}=a(\xi_i).
\]
The radius $d$ will be decreased once more in the degree argument. All parameters
$b_{i,\xi}$ belong to the fixed compact interval $[a_*,a^*]$. All constants
below are uniform for the allowed centers. The number of peaks is fixed.
The scaled variable is $y=x/\varepsilon$. We use the spaces $X,Y$ in
\eqref{eq:XY} and the real $L^2$ pairing for pairs. Define
\begin{align}
 \widehat A_\varepsilon&=-\Delta+a(\varepsilon y),\notag\\
 W_{\varepsilon,\xi}(y)&=\sum_{i=1}^k
       (U_{b_{i,\xi}},V_{b_{i,\xi}})(y-\xi_i/\varepsilon),
                                                   \label{con-ansatz}\\
 Z_{ij}(y)&=\partial_j(U_{b_{i,\xi}},V_{b_{i,\xi}})
                   (y-\xi_i/\varepsilon),\qquad 1\le i\le k,
                         \quad 1\le j\le N.\notag
\end{align}
Write $(\widehat U_i,\widehat V_i)$ for the $i$th profile in
\eqref{con-ansatz}, with parameter $b_{i,\xi}$ and center $\xi_i/\varepsilon$.
The differential equation map $\mathcal E_\varepsilon:X\to Y$ and its linearization are
\begin{equation}\label{con-F}
 \begin{split}
 \mathcal E_\varepsilon(u,v)&=
  \bigl(\widehat A_\varepsilon u-(v_+)^p,
        \widehat A_\varepsilon v-(u_+)^q\bigr),\\
 L_{\varepsilon,\xi}=\mathcal E'_\varepsilon(W_{\varepsilon,\xi})
 &=\begin{pmatrix}
 \widehat A_\varepsilon&-p(W_{\varepsilon,\xi}^{(2)})^{p-1}\\
 -q(W_{\varepsilon,\xi}^{(1)})^{q-1}&\widehat A_\varepsilon
 \end{pmatrix}.
 \end{split}
\end{equation}
Here the superscripts distinguish the two components of $W$.
Positive-part powers enter because the projected solution is not yet
known to be positive. The source directions are $SZ_{ij}$. The orthogonality conditions are imposed against $Z_{ij}$.

\subsection{Scalar estimates and interactions between different peaks}

The scalar inverse controls the part of a possible approximate kernel which
is not localized near any peak. Its estimate uses only the positive lower
bound for the potential.

\begin{lemma}\label{lem:con-scalar}
The map $\widehat A_\varepsilon:H^2\to L^2$ is an isomorphism with
$\|u\|_{H^2}\le C\|\widehat A_\varepsilon u\|_2$.
Its restriction from $H^2\cap W^{2,t}$ onto $L^2\cap L^t$ satisfies
\begin{equation}\label{con-scalar-inverse}
 \|u\|_{H^2}+\|u\|_{W^{2,t}}
 \le C\bigl(\|\widehat A_\varepsilon u\|_2+
                    \|\widehat A_\varepsilon u\|_t\bigr)
\end{equation}
with $C$ independent of $\varepsilon>0$.
\end{lemma}

\begin{proof}
Lax--Milgram gives a unique $H^1$ solution of
$\widehat A_\varepsilon u=f$ for $f\in L^2$. Testing with $u$
gives $\|u\|_2\le a_*^{-1}\|f\|_2$. Since
$-\Delta u=f-a(\varepsilon y)u$, the global elliptic estimate gives
$\|u\|_{H^2}\le C\|f\|_2$.

If $f\in L^t$ as well, use the admissible
test $u\min\{|u|,K\}^{t-2}$ and set
$Q_K=|u|^2\min\{|u|,K\}^{t-2}$. The gradient term is nonnegative, so
\[
 a_*\int Q_K
 \le\int |f|\,|u|\min\{|u|,K\}^{t-2}
 \le\|f\|_t\left(\int Q_K\right)^{(t-1)/t}.
\]
Here the last inequality follows from
$|u|\min\{|u|,K\}^{t-2}\le Q_K^{(t-1)/t}$.
Letting $K\to\infty$ gives $\|u\|_t\le a_*^{-1}\|f\|_t$.
The global estimate
\[
 \|u\|_{W^{2,t}}\le C_t(\|\Delta u\|_t+\|u\|_t)
\]
then proves \eqref{con-scalar-inverse}. All constants depend only on
$a_*,a^*,N,t$ and are therefore independent of $\varepsilon$.
\end{proof}

The next lemma estimates the error of the sum of profiles and its
translation derivatives. The estimates also cover powers below two. Set
\begin{equation}\label{con-r-definition}
 r_\varepsilon(\xi)=\varepsilon\sum_i|\nabla a(\xi_i)|
                     +\varepsilon^2+e^{-\kappa/\varepsilon},
\end{equation}
where the positive constant $\kappa$ can be decreased in the proof below.

\begin{lemma}\label{lem:con-basic-estimates}
For a fixed $\kappa>0$, uniformly for small $\varepsilon$ and allowed
centers, one has
\begin{align}
 \|\mathcal E_\varepsilon(W_{\varepsilon,\xi})\|_Y
       &\le C r_\varepsilon(\xi),                       \label{con-residual}\\
 \|L_{\varepsilon,\xi}Z_{ij}\|_Y
       &\le C\varepsilon+Ce^{-\kappa/\varepsilon},       \label{con-approx-kernel}\\
 \langle Z_{ij},Z_{h\ell}\rangle
       &=\delta_{ih}\delta_{j\ell}d_{b_{i,\xi}}
                              +O(e^{-\kappa/\varepsilon}),
                                                            \label{con-Gram}
\end{align}
where
\[
 d_b=\int\bigl((\partial_1U_b)^2+(\partial_1V_b)^2\bigr)>0.
\]
In particular the Gram matrices in \eqref{con-Gram} and their inverses are
uniformly bounded. Also $\|W_{\varepsilon,\xi}\|_X$ and
$\sum_{i,j}\|Z_{ij}\|_X$ are uniformly bounded.
\end{lemma}

\begin{proof}
The scaled separations $D_{ih}=|\xi_i-\xi_h|/\varepsilon$ satisfy $D_{ih}\ge c/\varepsilon$ for $i\ne h$. For $\alpha,\beta>0$, $m=\min\{\alpha,\beta\}$ and $1\le s<\infty$, the inequality
\[
 \alpha|y-\xi_i/\varepsilon|+\beta|y-\xi_h/\varepsilon|
 \ge \tfrac m2D_{ih}+\tfrac m2|y-\xi_i/\varepsilon|
\]
gives
\begin{equation}\label{con-exponential-product}
 \big\|e^{-\alpha|\cdot-\xi_i/\varepsilon|}
         e^{-\beta|\cdot-\xi_h/\varepsilon|}\big\|_s
       \le C_{\alpha,\beta,s}e^{-mD_{ih}/2}.
\end{equation}
For $s>1$ and bounded $A,B\ge0$,
\[
 |(A+B)^s-A^s-B^s|\le C_s(A^{s-1}B+AB^{s-1}).
\]
Iterating for the fixed number of profiles and applying \eqref{con-exponential-product} shows that the nonlinear interactions are $O(e^{-\kappa/\varepsilon})$ in $Y$. The remaining residual terms are controlled by
\[
 |a(\xi_i+\varepsilon z)-a(\xi_i)|
 \le\varepsilon|\nabla a(\xi_i)||z|+C\varepsilon^2|z|^2
\]
and the weighted profile estimates. This proves \eqref{con-residual}.

Differentiating the autonomous equation, the first component of $L_{\varepsilon,\xi}Z_{ij}$ is
\[
 (a(\varepsilon y)-b_{i,\xi})\partial_j\widehat U_i
 -p\left[\left(\sum_h\widehat V_h\right)^{p-1}
                              -\widehat V_i^{p-1}\right]
                         \partial_j\widehat V_i.
\]
The potential term is $O(\varepsilon)$ in $Y$. On bounded nonnegative intervals,
\[
 |(A+B)^\nu-A^\nu|\le C B^{\min\{1,\nu\}},\qquad\nu>0.
\]
Thus the second term is exponentially small by \eqref{con-exponential-product}; this uses no derivative of $s^{p-1}$ at zero. The other component is treated identically, proving \eqref{con-approx-kernel}.

For the Gram matrix, diagonal blocks are $d_{b_{i,\xi}}I_N$ by \eqref{eq:profile-Gram}, and off-diagonal blocks are exponentially small by \eqref{con-exponential-product} with $s=1$. The positive minimum of $d_b$ gives uniform invertibility. The remaining norm bounds follow from \eqref{eq:profile-weighted-derivatives} and translation invariance.
\end{proof}

\subsection{A uniform inverse with moment conditions}

We prove a uniform inverse estimate. The complementary quadratic form
need not be positive, so the proof uses local compactness and the
autonomous kernel instead of coercivity.

\begin{lemma}\label{lem:con-inverse}
There are $\varepsilon_0,C>0$ such that, for
$0<\varepsilon<\varepsilon_0$ and $\xi\in\overline{\mathcal D}$, the map
\begin{equation}\label{con-bordered}
 \mathcal T_{\varepsilon,\xi}(w,c)=
 \left(L_{\varepsilon,\xi}w-\sum_{i,j}c_{ij}SZ_{ij},
                       (\langle w,Z_{ij}\rangle)_{i,j}\right)
\end{equation}
is an isomorphism from $H^2(\R^N)^2\times\R^{kN}$ onto
$L^2(\R^N)^2\times\R^{kN}$. If
$\mathcal T_{\varepsilon,\xi}(w,c)=(h,\beta)$, then
\begin{equation}\label{con-L2-inverse}
 \|w\|_{H^2}+|c|\le C(\|h\|_2+|\beta|).
\end{equation}
Its restriction from $X\times\R^{kN}$ onto $Y\times\R^{kN}$ is also an
isomorphism, with
\begin{equation}\label{con-inverse-estimate}
 \|w\|_X+|c|\le C(\|h\|_Y+|\beta|).
\end{equation}
\end{lemma}

\begin{proof}
We first prove \eqref{con-L2-inverse}. If this estimate did not hold uniformly, there
would be $\varepsilon_n\to0$, allowed centers $\xi_n$, and pairs with
\[
 \|w_n\|_{H^2}+|c_n|=1,\qquad
 \mathcal T_{\varepsilon_n,\xi_n}(w_n,c_n)=(h_n,\beta_n),\qquad
 \|h_n\|_2+|\beta_n|\to0.
\]
Since $L^*S=SL$, testing the equation against $SZ_{ij}$ gives
\[
 \sum_{h,\ell}c_{h\ell}\langle Z_{h\ell},Z_{ij}\rangle
 =\langle w,SLZ_{ij}\rangle-\langle h,SZ_{ij}\rangle.
\]
The approximate-kernel and Gram estimates therefore imply
\begin{equation}\label{con-c-small}
 |c_n|\le C\bigl(\|h_n\|_2+
        (\varepsilon_n+e^{-\kappa/\varepsilon_n})\|w_n\|_2\bigr)\to0.
\end{equation}
For each $i$, take a subsequence such that $b_{i,\xi_n}\to b_i^\circ$.
The translated functions $w_n(\cdot+\xi_{i,n}/\varepsilon_n)$ converge
weakly in $H^2$ and strongly in $L^2$ on bounded sets to a pair
$w_i^\circ$. Passing to the equation gives $L_{b_i^\circ}w_i^\circ=0$.
The moments pass to the limit by the uniform $L^2$ bounds and the
exponential tails of the profile derivatives. Thus
\[
 \langle w_i^\circ,\partial_j(U_{b_i^\circ},V_{b_i^\circ})\rangle=0,
 \qquad 1\le j\le N.
\]
Lemma~\ref{lem:autonomous-nondegeneracy} and the translation Gram matrix
give $w_i^\circ=0$.

Put $B_n=L_{\varepsilon_n,\xi_n}-\widehat A_{\varepsilon_n}I$.
On the union of fixed scaled balls about the centers, local compactness
gives $B_nw_n\to0$ in $L^2$. Outside those balls, the coefficient norm
is at most $Ce^{-\sigma_0R}$, where
$\sigma_0=\sigma\min\{p-1,q-1\}>0$. Taking first $n\to\infty$ and then
$R\to\infty$ yields $\|B_nw_n\|_2\to0$. The scalar $L^2\to H^2$
estimate, applied to
\[
 \widehat A_{\varepsilon_n}w_n
   =h_n+\sum_{i,j}c_{ij,n}SZ_{ij,n}-B_nw_n,
\]
now contradicts the normalization.

For fixed $\varepsilon,\xi$, the decaying multiplication coefficients
are compact from $H^2$ to $L^2$. The moment rows and auxiliary columns
have finite rank. Hence $\mathcal T_{\varepsilon,\xi}$ is Fredholm of
index zero, as a compact perturbation of
$(w,c)\mapsto(\widehat A_\varepsilon w,c)$.
Estimate \eqref{con-L2-inverse} gives injectivity and therefore
surjectivity.

If $h\in Y$, its norms in every intermediate $L^r$, $2\le r\le t$,
are controlled by $\|h\|_Y$. The equation and the uniform boundedness of
the coefficients give
\[
 \|w\|_{W^{2,r}}
 \le C_r\bigl(\|h\|_r+|c|+\|w\|_r\bigr).
\]
Starting with $r_0=2$, choose successive exponents by
\[
 \frac1{r_{j+1}}=\max\left\{\frac1t,\frac1{r_j}-\frac1N\right\}.
\]
Sobolev embedding reaches $t$ in finitely many steps. The profile source
terms are uniformly bounded in all these norms. Together with
\eqref{con-L2-inverse}, this proves \eqref{con-inverse-estimate}
and the assertion on $X\times\R^{kN}$.
\end{proof}

\subsection{Solving the projected nonlinear equation}

The next lemma gives a unique small correction for each fixed configuration
of centers. Uniqueness of the centers, and hence of the full solution, is
proved by the local Pohozaev identity in Theorem~\ref{thm:direct-uniqueness}.

\begin{lemma}\label{lem:con-correction}
For all sufficiently small $\varepsilon$ and every allowed $\xi$, there
is a unique pair $(w_{\varepsilon,\xi},c_{\varepsilon,\xi})$ in a fixed
small neighborhood of zero in $X\times\R^{kN}$ such that
\begin{equation}\label{con-projected}
 \mathcal E_\varepsilon(W_{\varepsilon,\xi}+w_{\varepsilon,\xi})
                  =\sum_{i,j}c_{ij}SZ_{ij},\qquad
 \langle w_{\varepsilon,\xi},Z_{ij}\rangle=0.
\end{equation}
For fixed $\varepsilon$, this pair depends continuously on $\xi$ in
$X\times\R^{kN}$, and
\begin{equation}\label{con-correction-size}
 \|w_{\varepsilon,\xi}\|_X+|c_{\varepsilon,\xi}|
 \le C r_\varepsilon(\xi).
\end{equation}
\end{lemma}

\begin{proof}
Let $\tau=\min\{1,p-1,q-1\}$ and put
\[
 N_\xi(w)=-\bigl[\mathcal N(W_{\varepsilon,\xi}+w)
       -\mathcal N(W_{\varepsilon,\xi})
       -D\mathcal N(W_{\varepsilon,\xi})w\bigr].
\]
Lemma~\ref{lem:power-map} gives, uniformly in the centers,
\begin{equation}\label{con-nonlinear}
 \|N_\xi(w)-N_\xi(\widetilde w)\|_Y
 \le C_N(\|w\|_X^\tau+\|\widetilde w\|_X^\tau)
                              \|w-\widetilde w\|_X.
\end{equation}
Let $C_T$ bound the inverse of Lemma~\ref{lem:con-inverse}, and let
$C_R$ be the constant in \eqref{con-residual}. On $X\times\R^{kN}$,
with norm $\|w\|_X+|c|$, define
\[
 \mathcal P_{\varepsilon,\xi}(w,c)
 =\mathcal T_{\varepsilon,\xi}^{-1}
       \bigl(-\mathcal E_\varepsilon(W_{\varepsilon,\xi})-N_\xi(w),0\bigr).
\]
Choose a fixed $\delta>0$ with $2C_TC_N\delta^\tau\le1/2$. For small
$\varepsilon$, this map preserves the ball of radius
$2C_TC_Rr_\varepsilon(\xi)\le\delta$ and has Lipschitz constant at most
$1/2$. The contraction theorem gives a fixed point and
\eqref{con-correction-size}.

Any solution in the fixed $\delta$ ball satisfies
\[
 \|w\|_X+|c|
 \le C_TC_Rr_\varepsilon(\xi)+C_TC_N\|w\|_X^{1+\tau}
 \le C_TC_Rr_\varepsilon(\xi)+\tfrac12(\|w\|_X+|c|).
\]
It therefore lies in the contraction ball. This proves uniqueness in
the fixed neighborhood.

For fixed $\varepsilon$, the profiles and modes are continuous in the
norms used in \eqref{con-bordered}. The power estimate and the inverse
identity give operator-norm continuity of
$\mathcal T_{\varepsilon,\xi}^{-1}$. Subtracting the fixed-point
equations at nearby centers and absorbing the contraction term proves
continuity of $(w_{\varepsilon,\xi},c_{\varepsilon,\xi})$ in
$X\times\R^{kN}$.
\end{proof}

\subsection{Decay of the projected pair}

The projected solution may change sign and still contains the auxiliary
source terms. We therefore prove its decay directly. The estimate is used
to control the cutoff terms in the position equation.

\begin{lemma}\label{lem:con-tail}
Let $z_{\varepsilon,\xi}=(u_{\varepsilon,\xi},v_{\varepsilon,\xi})
=W_{\varepsilon,\xi}+w_{\varepsilon,\xi}$. There are $C,\eta>0$,
independent of the allowed centers and small $\varepsilon$, such that
\begin{equation}\label{con-tail}
 |u_{\varepsilon,\xi}(y)|+|v_{\varepsilon,\xi}(y)|
 +|\nabla u_{\varepsilon,\xi}(y)|+|\nabla v_{\varepsilon,\xi}(y)|
       \le C\sum_i e^{-\eta|y-\xi_i/\varepsilon|}.
\end{equation}
\end{lemma}

\begin{proof}
Write $z_{\varepsilon,\xi}=(u,v)$ and $T=|u|+|v|$. Choose $h>0$ with
$h^{p-1},h^{q-1}\le a_*/2$. Profile decay and
$\|w_{\varepsilon,\xi}\|_\infty=O(\varepsilon)$ imply $|u|,|v|\le h$
on
\[
 \Omega_{\varepsilon,\xi}
 =\R^N\setminus\bigcup_i\overline{B_R(\xi_i/\varepsilon)}
\]
for a fixed large $R$ and all sufficiently small $\varepsilon$.
The source in \eqref{con-projected} is bounded by
\[
 C|c_{\varepsilon,\xi}|\sum_i e^{-\sigma|y-\xi_i/\varepsilon|}.
\]
Applying the absolute-value inequality to both equations gives
\begin{equation}\label{con-Kato-exterior}
 (-\Delta+a_*/2)T
 \le C|c_{\varepsilon,\xi}|\sum_i e^{-\sigma|y-\xi_i/\varepsilon|}
 \quad\text{in }\Omega_{\varepsilon,\xi}.
\end{equation}
This inequality follows by replacing $|s|$ by
$j_\delta(s)=\sqrt{s^2+\delta^2}-\delta$ and passing to the limit:
$j_\delta''\ge0$ gives the required sign of the Laplacian term.

Choose $0<\eta<\min\{\sigma,\sqrt{a_*/2}\}$ and put
$\Phi=K\sum_i e^{-\eta|y-\xi_i/\varepsilon|}$. For $r>0$,
\[
 (-\Delta+a_*/2)e^{-\eta r}
 =\left(a_*/2-\eta^2+\frac{(N-1)\eta}{r}\right)e^{-\eta r}.
\]
The uniform bounds for $T$ and $c_{\varepsilon,\xi}$ allow a fixed $K$
so that $\Phi\ge T$ on the boundary and its equation dominates
\eqref{con-Kato-exterior}. Testing the difference with
$(T-\Phi)_+\in H_0^1(\Omega_{\varepsilon,\xi})$ gives
\[
 \int_{\Omega_{\varepsilon,\xi}}
 \left(|\nabla(T-\Phi)_+|^2+(a_*/2)(T-\Phi)_+^2\right)\le0.
\]
The test is justified by cutoffs at infinity, since $T,\Phi\in H^1$.
Thus $T\le\Phi$ outside the balls. Increasing $K$ gives the same bound
inside them.

The exponential weight changes by a bounded factor on each ball
$B_2(y_0)$. The equations, the bound just proved and the source decay
therefore give
\[
 \|u\|_{W^{2,t}(B_1(y_0))}+\|v\|_{W^{2,t}(B_1(y_0))}
 \le C\sum_i e^{-\eta|y_0-\xi_i/\varepsilon|}.
\]
Since $t>N$, Sobolev embedding gives the gradient estimate in
\eqref{con-tail}.
\end{proof}

\subsection{The localized translation identity and its expansion}

Choose $\chi_i\in C_c^\infty(B_{3d}(P_i))$ with
$0\le\chi_i\le1$ and $\chi_i=1$ on $B_{2d}(P_i)$. Their supports are
disjoint. Every point in the support of a derivative of any $\chi_i$ has
a fixed positive distance from every allowed center. In the next formulas,
$u,v$ mean the two components of $z_{\varepsilon,\xi}$ in scaled
coordinates. Subscripts on $a$ and $\chi_i$ denote physical derivatives.
Define
\begin{align}
 e_{\varepsilon,\xi}(y)
 &=\nabla u\cdot\nabla v+a(\varepsilon y)uv
           -\frac{(u_+)^{q+1}}{q+1}-\frac{(v_+)^{p+1}}{p+1},\notag\\
 \mathcal F_{ij}(\xi)
 &=\int\chi_i(\varepsilon y)\partial_j a(\varepsilon y)uv\,dy\notag\\
 &\quad-\int\left[\sum_h(\partial_h\chi_i)(\varepsilon y)
        (\partial_hu\partial_jv+\partial_hv\partial_ju)
       -(\partial_j\chi_i)(\varepsilon y)e_{\varepsilon,\xi}\right]dy,
                                                         \label{con-force}\\
 H_{ij,h\ell}(\xi)
 &=\int\chi_i(\varepsilon y)\partial_jz_{\varepsilon,\xi}\cdot Z_{h\ell}\,dy.
                                                         \label{con-H}
\end{align}
All these integrals are finite. For example, $u,v\in H^1\cap L^\infty$,
and $|u|^{q+1}\le\|u\|_\infty^{q-1}|u|^2$ and the analogous inequality
for $v$ imply integrability of the nonlinear potential terms.

The next identity reduces the vanishing of the auxiliary terms to a
finite-dimensional equation for the centers. Its error estimate uses
weighted moments of the profiles and unweighted norms of the correction.

\begin{lemma}\label{lem:con-force}
The matrix $H(\xi)$ is uniformly invertible for small $\varepsilon$, and
\begin{equation}\label{con-force-exact}
 H(\xi)c_{\varepsilon,\xi}=-\varepsilon \mathcal F(\xi).
\end{equation}
The map $\mathcal F$ is continuous on $\overline{\mathcal D}$. With
$M_b=\int U_bV_b>0$, one has
\begin{equation}\label{con-force-error}
 \left|\mathcal F_i(\xi)-M_{b_{i,\xi}}\nabla a(\xi_i)\right|
 \le C\left[\varepsilon^2+
 (|\nabla a(\xi_i)|+\varepsilon)r_\varepsilon(\xi)
 +r_\varepsilon(\xi)^2+e^{-\kappa_1/\varepsilon}\right]
\end{equation}
for some $\kappa_1>0$. In particular,
$\mathcal F\to \mathcal F_0$ uniformly on $\overline{\mathcal D}$, where
$\mathcal F_0(\xi)=(M_{a(\xi_i)}\nabla a(\xi_i))_{i=1}^k$.
\end{lemma}

\begin{proof}
The localized translation identity can be written using
\[
 T_{hj}=\delta_{hj}e_{\varepsilon,\xi}
             -(\partial_hu\partial_jv+\partial_hv\partial_ju).
\]
The product rule and the projected equations give, in distributions,
\[
 \sum_h\partial_hT_{hj}
 =\varepsilon\partial_ja(\varepsilon y)uv
       +\sum_{h,\ell}c_{h\ell}Z_{h\ell}\cdot\partial_jz.
\]
Indeed, the terms left after cancellation of the mixed second derivatives are
$\varepsilon\partial_ja(\varepsilon y)uv+
(\widehat A_\varepsilon u-(v_+)^p)\partial_jv+
(\widehat A_\varepsilon v-(u_+)^q)\partial_ju$.
Multiplying by $\chi_i(\varepsilon y)$ and integrating proves \eqref{con-force-exact}, with the signs in \eqref{con-force}. The calculation is justified by $H^2$ approximation and the Sobolev chain rule.

Substituting $z=W+w$ in \eqref{con-H}, the same-peak term is
$d_{b_{i,\xi}}\delta_{j\ell}+O(e^{-\kappa/\varepsilon})$. Other profile terms are exponentially small by \eqref{con-exponential-product} and separation of the cutoffs. The correction term is bounded by $C\|\nabla w\|_2\le Cr_\varepsilon(\xi)$. Thus
\begin{equation}\label{con-H-expansion}
 H(\xi)=\diag(d_{b_{i,\xi}}I_N)_{i=1}^k
                 +O(r_\varepsilon(\xi)+e^{-\kappa/\varepsilon}),
\end{equation}
which proves invertibility. The matrix $H$ need not be symmetric.

On each cutoff transition region, \eqref{con-tail} bounds the solution and its gradient by $Ce^{-c/\varepsilon}$. The scaled volume is $O(\varepsilon^{-N})$, so the stress terms in \eqref{con-force} are $O(e^{-\kappa_1/\varepsilon})$ after reducing $\kappa_1$.

For the remaining integral, decompose
\[
 uv=\sum_{h,\ell}\widehat U_h\widehat V_\ell
        +\sum_h\widehat U_hw_2+\sum_\ell\widehat V_\ell w_1+w_1w_2.
\]
All profile products except $\widehat U_i\widehat V_i$ contribute exponentially small errors. For this term put $z=y-\xi_i/\varepsilon$ and use
\[
 \nabla a(\xi_i+\varepsilon z)
 =\nabla a(\xi_i)+\varepsilon D^2a(\xi_i)z
                                      +O(\varepsilon^2|z|^2).
\]
The radial density $U_bV_b$ has zero first moment and a bounded second moment. Removing the cutoff therefore gives
\begin{equation}\label{con-pure-profile-force}
 \int\chi_i(\varepsilon y)\nabla a(\varepsilon y)
                    \widehat U_i\widehat V_i\,dy
 =M_{b_{i,\xi}}\nabla a(\xi_i)+O(\varepsilon^2)
                                      +O(e^{-\kappa_1/\varepsilon}).
\end{equation}
For the mixed terms at this center,
\[
 \left|\int\chi_i(\varepsilon y)\nabla a(\varepsilon y)\widehat U_iw_2\right|
 \le\bigl(|\nabla a(\xi_i)|\|U_{b_{i,\xi}}\|_2
       +C\varepsilon\||z|U_{b_{i,\xi}}\|_2\bigr)\|w_2\|_2
 \le C(|\nabla a(\xi_i)|+\varepsilon)r_\varepsilon(\xi).
\]
The term with $\widehat V_iw_1$ has the same bound. Other mixed terms have exponentially small profile factors on the cutoff support. Finally,
\[
 \left|\int\chi_i(\varepsilon y)\nabla a(\varepsilon y)w_1w_2\right|
 \le\|\nabla a\|_\infty\|w_1\|_2\|w_2\|_2
 \le Cr_\varepsilon(\xi)^2.
\]
This proves \eqref{con-force-error}; no weighted norm of $w$ is used.

Continuity follows from the continuity in $X$ of the profiles and correction. The gradient products use $H^1$ continuity; the nonlinear potentials use $L^2$ continuity and the uniform $L^\infty$ bound. Since $r_\varepsilon=O(\varepsilon)$ uniformly, \eqref{con-force-error} also gives uniform convergence to $\mathcal F_0$.
\end{proof}

\subsection{Solving the position equation}

The limiting position map has an isolated nondegenerate zero at the
prescribed configuration. The following argument produces a zero of the
perturbed map and then improves its location from $o(1)$ to
$O(\varepsilon^2)$.

\begin{proposition}\label{prop:con-existence}
For every sufficiently small $\varepsilon>0$ there are centers
$\xi_{i,\varepsilon}\in B_d(P_i)$ with
$c_{\varepsilon,\xi_\varepsilon}=0$. The associated projected pair is
strictly positive, solves the scaled system, and satisfies
\begin{equation}\label{con-final-size}
 |\xi_{i,\varepsilon}-P_i|\le C\varepsilon^2,\qquad
 \|w_{\varepsilon,\xi_\varepsilon}\|_X\le C\varepsilon^2.
\end{equation}
Rescaling to the physical variable gives the solution and the expansion
\eqref{eq:main-decomposition} with the stated estimates.
\end{proposition}

\begin{proof}
Decrease $d$ so that $P_i$ is the only critical point in $\overline{B_d(P_i)}$ and
\begin{equation}\label{con-gradient-lower}
 |\nabla a(\xi_i)|\ge c_0|\xi_i-P_i|
 \quad(\xi_i\in\overline{B_d(P_i)}).
\end{equation}
Both assertions follow from the invertibility of $D^2a(P_i)$ and Taylor expansion. At its unique zero $P=(P_1,\ldots,P_k)$, the limiting map has derivative
$D\mathcal F_0(P)=\diag(M_{b_i}D^2a(P_i))$. Thus
\[
 \deg(\mathcal F_0,\mathcal D,0)
 =\prod_i\operatorname{sgn}\det D^2a(P_i)\ne0.
\]
Uniform convergence ensures that the straight-line homotopy from $\mathcal F_0$ to $\mathcal F$ has no boundary zeros. Therefore $\mathcal F$ has a zero $\xi_\varepsilon\in\mathcal D$ for every small $\varepsilon$. By \eqref{con-force-exact}, its auxiliary coefficients vanish.

At any such zero set $g=\max_i|\nabla a(\xi_i)|$. Since $\min_b M_b>0$, \eqref{con-force-error} yields
\[
 g\le C\varepsilon^2+C(g+\varepsilon)r_\varepsilon
                         +Cr_\varepsilon^2+Ce^{-\kappa_1/\varepsilon}.
\]
Absorbing the exponential errors into $O(\varepsilon^2)$, we have
$r_\varepsilon\le C(\varepsilon g+\varepsilon^2)$. Since $g$ is uniformly bounded, substitution gives
\[
 g\le C\varepsilon^2+C\varepsilon g.
\]
For small $\varepsilon$, absorption gives $g=O(\varepsilon^2)$. Formula \eqref{con-gradient-lower} and \eqref{con-correction-size} now prove \eqref{con-final-size}.

With $c=0$, the first equation is $\widehat A_\varepsilon u=(v_+)^p$. Testing with $u_-:=\max\{-u,0\}$ gives
\[
 -\int(|\nabla u_-|^2+a(\varepsilon y)u_-^2)
 =\int(v_+)^pu_-\ge0,
\]
so $u\ge0$; similarly $v\ge0$. The positive profiles and the $O(\varepsilon^2)$ supremum error show that neither component is zero. Local elliptic regularity gives $C^{2,\alpha}_{\rm loc}$ regularity, and the strong maximum principle gives strict positivity. Differentiating the equations once also gives $H^3\cap W^{3,t}$ regularity, since the first derivatives of the powers and of $a$ are bounded on the solution range. Decay follows from \eqref{con-tail}.

Finally set $z_\varepsilon(x)=z_{\varepsilon,\xi_\varepsilon}(x/\varepsilon)$. The scaled norms agree with \eqref{eq:scaled-norm}, and the moment integrals acquire only the factor $\varepsilon^N$. This proves the existence and estimates in Theorem~\ref{thm:main-existence}.
\end{proof}

\section{Concentration estimates and local uniqueness}
\label{sec:concentration-estimates}

We first establish decay and center estimates for every admissible solution.
Subtracting the local Pohozaev identities then gives local uniqueness.

\subsection{Admissible families and estimates on the whole space}

We use the fixed profiles $U_i=U_{b_i}$ and $V_i=V_{b_i}$ at $b_i=a(P_i)$.
They have the decay and kernel specified in Section~\ref{sec:profiles}, and
\begin{equation}\label{eq:unique-frozen}
 (-\Delta+b_i)U_i=V_i^p,\qquad
 (-\Delta+b_i)V_i=U_i^q.
\end{equation}
Recall $M_i=\int U_iV_i>0$. Here we assume convergence to these profiles. In Section~\ref{sec:positive-class} we derive it from an energy bound and nonvanishing near the prescribed points, using the uniqueness of positive solutions of the limiting system.

A family of positive classical solutions $(u_\varepsilon,v_\varepsilon)$
is called admissible if both components tend to zero at infinity and there
are centers $\eta_{i,\varepsilon}\to P_i$ such that
\begin{equation}\label{eq:unique-class}
 \left\|u_\varepsilon-
       \sum_{i=1}^k U_i\left(\frac{\cdot-\eta_{i,\varepsilon}}{\varepsilon}\right)
 \right\|_{\infty}
 +\left\|v_\varepsilon-
       \sum_{i=1}^k V_i\left(\frac{\cdot-\eta_{i,\varepsilon}}{\varepsilon}\right)
 \right\|_{\infty}\longrightarrow0.
\end{equation}
This is a global assumption.  Convergence on fixed scaled neighborhoods
of the listed centers would not exclude extra peaks elsewhere.  We do not
include an $H^1$ assumption in this definition.  The integrability needed
in the proof will follow from the equations and \eqref{eq:unique-class}.
The parameter may range over an interval or over any sequence tending to
zero.  In the estimates below, the constants depend only on the fixed
profiles, the potential and the prescribed points.  The threshold for
small $\varepsilon$ may depend on the rate of convergence in
\eqref{eq:unique-class}.

The exterior comparison below requires no global integrability assumption.

\begin{lemma}\label{lem:unique-exterior-comparison}
Let $Q_{1,\varepsilon},\ldots,Q_{m,\varepsilon}$ be a finite list of
points, let $R_0>0$, and set
\[
 E_\varepsilon=\R^N\setminus
       \bigcup_{\nu=1}^m\overline{B_{R_0\varepsilon}(Q_{\nu,\varepsilon})}.
\]
The balls in this definition need not be disjoint.  Suppose that $T$ is
nonnegative, continuous and locally in $H^1$, is bounded by $A$ on
$\R^N$, tends to zero at infinity, and satisfies
\[
 (-\varepsilon^2\Delta+c_0)T\leq0
       \quad\hbox{weakly in }E_\varepsilon,
       \qquad c_0>0.
\]
For any fixed $0<\sigma<\sqrt{c_0}$ there is $C=C(R_0,\sigma)$ such that
\begin{equation}\label{eq:unique-exterior-bound}
 T(x)\leq CA\sum_{\nu=1}^m
          e^{-\sigma|x-Q_{\nu,\varepsilon}|/\varepsilon}
          \quad\hbox{for all }x\in\R^N.
\end{equation}
No global Sobolev norm of $T$ is needed.
\end{lemma}

\begin{proof}
The case $A=0$ is immediate. Otherwise set
\[
 B(x)=2Ae^{\sigma R_0}\sum_{\nu=1}^m
                  e^{-\sigma|x-Q_{\nu,\varepsilon}|/\varepsilon}.
\]
For $r=|x-Q_{\nu,\varepsilon}|>0$,
\[
 (-\varepsilon^2\Delta+c_0)e^{-\sigma r/\varepsilon}
 =\left(c_0-\sigma^2+\frac{(N-1)\sigma\varepsilon}{r}\right)
 e^{-\sigma r/\varepsilon}\ge0.
\]
Thus $B$ is a supersolution on $E_\varepsilon$ and $B\ge2A$ on the
union of the closed balls. For $\delta>0$, the positive part of
$T-B-\delta$ has compact support in $E_\varepsilon$, since $T$ tends
to zero at infinity and is strictly below $B+\delta$ on the boundary.
It is therefore an admissible test, even if the balls intersect, and gives
\[
 \varepsilon^2\int_{E_\varepsilon}|\nabla(T-B-\delta)_+|^2
 +c_0\int_{E_\varepsilon}(T-B-\delta)_+^2\le0.
\]
Letting $\delta\downarrow0$ proves the bound outside the balls; inside
them it follows from $T\le A\le B$. No global integrability of $T$
is used.
\end{proof}

We apply this comparison to the sum of the two positive components.
Interior estimates then give the corresponding gradient bound.

\begin{lemma}\label{lem:unique-tails}
For an admissible family and every $0<\sigma<\sqrt{a_*}$ there is a
constant $C_\sigma>0$, independent of small $\varepsilon$, such that
\begin{equation}\label{eq:unique-base-tail}
 u_\varepsilon+v_\varepsilon
 +\varepsilon(|\nabla u_\varepsilon|+|\nabla v_\varepsilon|)
 \leq C_\sigma\sum_i e^{-\sigma|x-\eta_{i,\varepsilon}|/\varepsilon}.
\end{equation}
\end{lemma}

\begin{proof}
Fix $0<\sigma<\sqrt{a_*}$ and choose $\delta>0$ with
$\sigma^2<a_*-\delta$. Global convergence gives a uniform supremum
bound. Choose $s>0$ with $s^{p-1},s^{q-1}\le\delta$.
By profile decay and \eqref{eq:unique-class}, both components are below
$s$ outside $\bigcup_iB_{R_0\varepsilon}(\eta_{i,\varepsilon})$ for
fixed large $R_0$. There
\[
 (-\varepsilon^2\Delta+a)(u_\varepsilon+v_\varepsilon)
 =v_\varepsilon^p+u_\varepsilon^q
 \le\delta(u_\varepsilon+v_\varepsilon),
\]
so $(-\varepsilon^2\Delta+a_*-\delta)(u_\varepsilon+v_\varepsilon)\le0$.
Lemma~\ref{lem:unique-exterior-comparison} gives the zeroth-order bound
with the chosen $\sigma$. On each ball $B_{2\varepsilon}(x_0)$ the
exponential weights vary by a bounded factor. Rescaled interior
$W^{2,t}$ estimates and Sobolev embedding give the gradient bound
with the same $\sigma$, as in the proof of Lemma~\ref{lem:con-tail}.
\end{proof}

For each fixed $1\leq s<\infty$, integration of the exponential bound gives
\begin{equation}\label{eq:unique-integrability}
 \|u_\varepsilon\|_s+\|v_\varepsilon\|_s
 +\varepsilon(\|\nabla u_\varepsilon\|_s+
                     \|\nabla v_\varepsilon\|_s)
 \leq C_s\varepsilon^{N/s}.
\end{equation}
The power of $\varepsilon$ follows by scaling each exponential weight.
Thus every admissible solution belongs to $H^1(\R^N)^2$. The equations
also give its Laplacians in $L^2$, hence it belongs to $H^2(\R^N)^2$.

Rescaling about $\eta_{i,\varepsilon}$ also gives
\begin{equation}\label{eq:unique-core-c1}
 \bigl(u_\varepsilon(\eta_{i,\varepsilon}+\varepsilon y),
       v_\varepsilon(\eta_{i,\varepsilon}+\varepsilon y)\bigr)
          \longrightarrow(U_i(y),V_i(y))
          \quad\hbox{in }C^1_{\rm loc}(\R^N)^2.
\end{equation}
Indeed, \eqref{eq:unique-class} fixes the uniform limit on each bounded
ball, while the rescaled equations give uniform local $W^{2,t}$ bounds.
Compact embedding upgrades the convergence to $C^1$.

\subsection{A local Pohozaev identity and the peak centers}

The local translation identity determines the center error. We record
the boundary flux used in this calculation.
Write $u_h=\partial_hu$, $v_h=\partial_hv$, and put
\[
 e(u,v)=\varepsilon^2\nabla u\cdot\nabla v+a uv
                  -\frac{u^{q+1}}{q+1}-\frac{v^{p+1}}{p+1},\qquad
 T_{hj}=\delta_{hj}e(u,v)-\varepsilon^2(u_hv_j+v_hu_j).
\]
The calculation in Lemma~\ref{lem:con-force}, now in physical variables and without auxiliary sources, gives
\[
 \sum_h\partial_hT_{hj}
 =a_j uv+(-\varepsilon^2\Delta u+a u-v^p)v_j
            +(-\varepsilon^2\Delta v+a v-u^q)u_j
 =a_j uv.
\]
The classical regularity assumed here justifies integration over any fixed ball.

For an outward unit normal $\nu$ define the resulting flux by
\begin{align}\label{eq:unique-flux}
 \mathcal P_j(u,v;\nu)
 ={}&\nu_j\left(\varepsilon^2\nabla u\cdot\nabla v+a uv
       -\frac{u^{q+1}}{q+1}-\frac{v^{p+1}}{p+1}\right)\notag\\
 &-\varepsilon^2\bigl(\partial_\nu u\,\partial_jv
                   +\partial_\nu v\,\partial_ju\bigr).
\end{align}
The divergence theorem gives, for every fixed ball $B$,
\begin{equation}\label{eq:unique-local-pohozaev}
 \int_B\partial_ja\,uv\,dx
       =\int_{\partial B}\mathcal P_j(u,v;\nu)\,dS.
\end{equation}

The next lemma combines local convergence with the tail estimate to
pass to polynomial moments on the expanding domains.

\begin{lemma}\label{lem:unique-moment-convergence}
Fix $r>0$ such that $\overline{B_{2r}(P_i)}$ are pairwise disjoint, and let
\[
 D_{i,\varepsilon}=
       \{y:\eta_{i,\varepsilon}+\varepsilon y\in B_r(P_i)\},\quad
 u_i^\varepsilon(y)=u_\varepsilon(\eta_{i,\varepsilon}+\varepsilon y),\quad
 v_i^\varepsilon(y)=v_\varepsilon(\eta_{i,\varepsilon}+\varepsilon y).
\]
For every fixed polynomial $Q$ on $\R^N$,
\begin{equation}\label{eq:unique-polynomial-moments}
 \int_{D_{i,\varepsilon}}Q(y)u_i^\varepsilon(y)v_i^\varepsilon(y)\,dy
       \longrightarrow\int_{\R^N}Q(y)U_i(y)V_i(y)\,dy.
\end{equation}
In particular,
\begin{equation}\label{eq:unique-moments}
 \int_{D_{i,\varepsilon}}u_i^\varepsilon v_i^\varepsilon\longrightarrow M_i,
 \qquad
 \int_{D_{i,\varepsilon}}y_lu_i^\varepsilon v_i^\varepsilon\longrightarrow0,
 \qquad
 \int_{D_{i,\varepsilon}}|y|^2u_i^\varepsilon v_i^\varepsilon\leq C.
\end{equation}
\end{lemma}

\begin{proof}
For small $\varepsilon$, the domain $D_{i,\varepsilon}$ contains each fixed ball eventually and is contained in $B_{2r/\varepsilon}(0)$. All other centers stay a fixed physical distance from $B_r(P_i)$, so \eqref{eq:unique-base-tail} gives
\[
 u_i^\varepsilon(y)v_i^\varepsilon(y)
 \le C(e^{-\sigma|y|}+e^{-c/\varepsilon}),\qquad y\in D_{i,\varepsilon}.
\]
If $|Q(y)|\le C_Q(1+|y|^m)$, then
\[
 \int_{D_{i,\varepsilon}\setminus B_R}|Q|u_i^\varepsilon v_i^\varepsilon
 \le C_Q\int_{|y|>R}(1+|y|^m)e^{-\sigma|y|}\,dy
      +C_Q\varepsilon^{-N-m}e^{-c/\varepsilon}.
\]
On $B_R$, use \eqref{eq:unique-core-c1}. Taking $\varepsilon\to0$ and then $R\to\infty$ proves \eqref{eq:unique-polynomial-moments}. Radial symmetry and exponential decay of $U_iV_i$ give \eqref{eq:unique-moments}.
\end{proof}

The moment limits and local identity now give the center estimate.
The Hessian may be indefinite; only its invertibility is used.

\begin{lemma}\label{lem:unique-center-rate}
Every admissible family, with any choice of centers satisfying
\eqref{eq:unique-class}, satisfies
\begin{equation}\label{eq:unique-centers}
 \eta_{i,\varepsilon}=P_i+o(\varepsilon),\qquad i=1,\ldots,k.
\end{equation}
\end{lemma}

\begin{proof}
On $\partial B_r(P_i)$, the solution and its scaled gradient are exponentially small by \eqref{eq:unique-base-tail}. Thus every term of the flux in \eqref{eq:unique-flux} is exponentially small, and \eqref{eq:unique-local-pohozaev} gives
\begin{equation}\label{eq:unique-single-small-force}
 \left|\int_{B_r(P_i)}\partial_ja\,u_\varepsilon v_\varepsilon\right|
 \le Ce^{-c/\varepsilon}.
\end{equation}
Set $x=\eta_{i,\varepsilon}+\varepsilon y$ and divide by $\varepsilon^N$. Taylor's formula gives
\[
 \partial_ja(\eta_{i,\varepsilon}+\varepsilon y)
 =\partial_ja(\eta_{i,\varepsilon})
  +\varepsilon\sum_l\partial_{jl}a(\eta_{i,\varepsilon})y_l
  +O(\varepsilon^2|y|^2).
\]
The moments in \eqref{eq:unique-moments} therefore yield
\[
 \partial_ja(\eta_{i,\varepsilon})(M_i+o(1))
 +\varepsilon\sum_l\partial_{jl}a(\eta_{i,\varepsilon})o(1)
 +O(\varepsilon^2)=O(\varepsilon^{-N}e^{-c/\varepsilon}).
\]
Since $M_i>0$, this implies $\nabla a(\eta_{i,\varepsilon})=o(\varepsilon)$. Invertibility of $D^2a(P_i)$ gives $|\nabla a(P_i+h)|\ge c|h|$ for small $h$, proving \eqref{eq:unique-centers}. No choice of centers beyond \eqref{eq:unique-class} was used.
\end{proof}

\subsection{Local uniqueness}

We first estimate the difference of two admissible solutions. These
estimates will allow us to subtract the local Pohozaev identities.

\begin{lemma}\label{lem:normalized-difference-tail}
Let $z_\varepsilon^{(h)}=(u_\varepsilon^{(h)},v_\varepsilon^{(h)})$,
$h=1,2$, be two admissible families for the same profiles and points.
Suppose that $z_\varepsilon^{(1)}\ne z_\varepsilon^{(2)}$, and set
\[
 d_\varepsilon
 =\|u_\varepsilon^{(1)}-u_\varepsilon^{(2)}\|_\infty
  +\|v_\varepsilon^{(1)}-v_\varepsilon^{(2)}\|_\infty,\qquad
 (\phi_\varepsilon,\psi_\varepsilon)
 =\frac{z_\varepsilon^{(1)}-z_\varepsilon^{(2)}}{d_\varepsilon}.
\]
For some $C,\sigma>0$, independent of small $\varepsilon$,
\begin{equation}\label{eq:normalized-difference-tail}
 |\phi_\varepsilon|+|\psi_\varepsilon|
 +\varepsilon(|\nabla\phi_\varepsilon|+|\nabla\psi_\varepsilon|)
 \le C\sum_i e^{-\sigma|x-P_i|/\varepsilon}.
\end{equation}
\end{lemma}

\begin{proof}
By Lemma~\ref{lem:unique-center-rate}, the centers of both families
are $P_i+o(\varepsilon)$. Thus \eqref{eq:unique-base-tail} holds for
both families with $P_i$ in place of their centers, after increasing
the constant. Subtracting the equations gives
\begin{equation}\label{eq:normalized-difference-system}
 \begin{cases}
 (-\varepsilon^2\Delta+a)\phi_\varepsilon
                         =c_{p,\varepsilon}\psi_\varepsilon,\\
 (-\varepsilon^2\Delta+a)\psi_\varepsilon
                         =c_{q,\varepsilon}\phi_\varepsilon,
 \end{cases}
\end{equation}
where
\[
 c_{p,\varepsilon}
 =p\int_0^1\bigl(sv_\varepsilon^{(1)}
                   +(1-s)v_\varepsilon^{(2)}\bigr)^{p-1}\,ds,
\]
and $c_{q,\varepsilon}$ is defined in the same way using $u$ and $q$.
These nonnegative coefficients are uniformly bounded. Since $p,q>1$,
the tail estimates allow us to choose a fixed $R_0$ such that both
coefficients are at most $a_*/2$ outside
$\bigcup_iB_{R_0\varepsilon}(P_i)$. Kato's inequality gives
\[
 (-\varepsilon^2\Delta+a_*/2)
       (|\phi_\varepsilon|+|\psi_\varepsilon|)\le0
\]
there in the weak sense. The sum tends to zero at infinity and is
bounded by one everywhere. Lemma~\ref{lem:unique-exterior-comparison}
therefore proves its exponential bound. Interior $W^{2,t}$ estimates
for \eqref{eq:normalized-difference-system} on balls of radius
$2\varepsilon$, followed by Sobolev embedding, give the gradient
bound. Only boundedness of the coefficients is needed in this step.
\end{proof}

\begin{theorem}\label{thm:direct-uniqueness}
Any two admissible families for the same profiles and the same
nondegenerate critical points agree for every sufficiently small common
value of $\varepsilon$. The assertion also holds when $\varepsilon$
ranges over any sequence tending to zero.
\end{theorem}

\begin{proof}
If the assertion failed, we could choose $\varepsilon_n\to0$ and
two distinct admissible solutions at each $\varepsilon_n$. Use the
normalization in Lemma~\ref{lem:normalized-difference-tail}; in the
following we write $\varepsilon$ for this sequence.

For each $i$, rescale about $P_i$. The center estimate and admissibility
give
\[
 z_\varepsilon^{(h)}(P_i+\varepsilon y)\longrightarrow(U_i,V_i)
 \quad\text{in }C^1_{\rm loc}(\R^N)^2,\qquad h=1,2.
\]
The normalized difference is uniformly bounded. Interior estimates
and a subsequence give
\[
 (\phi_\varepsilon,\psi_\varepsilon)(P_i+\varepsilon y)
 \longrightarrow(\phi_i,\psi_i)
 \quad\text{in }C^1_{\rm loc}(\R^N)^2.
\]
The coefficients in \eqref{eq:normalized-difference-system} converge
locally uniformly to $pV_i^{p-1}$ and $qU_i^{q-1}$. Hence
$L_{b_i}(\phi_i,\psi_i)=0$. Estimate
\eqref{eq:normalized-difference-tail} gives exponential decay of the
limit, so it belongs to $L^2(\R^N)^2$; its equation then gives
$H^2$ regularity. By Lemma~\ref{lem:autonomous-nondegeneracy}, there
are constants $t_{i,l}$ such that
\begin{equation}\label{eq:normalized-difference-limit}
 (\phi_i,\psi_i)
 =\sum_{l=1}^N t_{i,l}\partial_l(U_i,V_i).
\end{equation}

Choose $r>0$ with the closed balls $\overline{B_{2r}(P_i)}$ pairwise
disjoint. Subtract \eqref{eq:unique-local-pohozaev} for the two
solutions on $B_r(P_i)$ and divide by $d_\varepsilon$. This yields
\begin{equation}\label{eq:difference-local-pohozaev}
 \int_{B_r(P_i)}\partial_j a\,
       \bigl(v_\varepsilon^{(1)}\phi_\varepsilon
                       +u_\varepsilon^{(2)}\psi_\varepsilon\bigr)
 =\int_{\partial B_r(P_i)}
       \frac{\mathcal P_j(z_\varepsilon^{(1)};\nu)
                    -\mathcal P_j(z_\varepsilon^{(2)};\nu)}
            {d_\varepsilon}\,dS.
\end{equation}
We estimate the quotient on the right before taking the limit.
For example,
\[
 \frac{\nabla u_\varepsilon^{(1)}\cdot\nabla v_\varepsilon^{(1)}
             -\nabla u_\varepsilon^{(2)}\cdot\nabla v_\varepsilon^{(2)}}
      {d_\varepsilon}
 =\nabla\phi_\varepsilon\cdot\nabla v_\varepsilon^{(1)}
       +\nabla u_\varepsilon^{(2)}\cdot\nabla\psi_\varepsilon,
\]
and
\[
 \frac{(u_\varepsilon^{(1)})^{q+1}
                      -(u_\varepsilon^{(2)})^{q+1}}{d_\varepsilon}
 =(q+1)\phi_\varepsilon
       \int_0^1(su_\varepsilon^{(1)}
                     +(1-s)u_\varepsilon^{(2)})^q\,ds.
\]
The other terms have the same product or mean-value form.
Equations \eqref{eq:unique-base-tail} and
\eqref{eq:normalized-difference-tail} show that every term in the
boundary quotient is $O(e^{-c/\varepsilon})$, since the factors
$\varepsilon^2$ in the gradient terms cancel the two gradient scales.
Consequently the right-hand side of
\eqref{eq:difference-local-pohozaev} is $O(e^{-c/\varepsilon})$.

Set $x=P_i+\varepsilon y$ in the left-hand side and divide by
$\varepsilon^{N+1}$. Since $\nabla a(P_i)=0$ and $D^2a$ is bounded,
\[
 \frac{\partial_j a(P_i+\varepsilon y)}{\varepsilon}
 \longrightarrow\sum_{m=1}^N\partial_{jm}a(P_i)y_m,\qquad
 \left|\frac{\partial_j a(P_i+\varepsilon y)}{\varepsilon}\right|
 \le C|y|.
\]
The bounds above give an integrable exponential for the contribution
of the $i$th peak. Contributions from the other peaks are
exponentially small on $B_r(P_i)$ and remain so after integration
and division by $\varepsilon^{N+1}$. Truncating first to a fixed
$y$-ball and then letting its radius tend to infinity gives
\[
 0=\sum_{m=1}^N\partial_{jm}a(P_i)
       \int_{\R^N}y_m(V_i\phi_i+U_i\psi_i).
\]
Substituting \eqref{eq:normalized-difference-limit} and integrating
by parts, we obtain
\[
 0=-M_i\sum_{l=1}^N\partial_{jl}a(P_i)t_{i,l},
 \qquad
 \int y_m\partial_l(U_iV_i)=-\delta_{ml}M_i.
\]
Thus $D^2a(P_i)t_i=0$. Its invertibility gives $t_i=0$ for every
$i$, and all the local limits of the normalized difference vanish.

For each fixed $R$, the normalized difference now tends uniformly to
zero on $\bigcup_iB_{R\varepsilon}(P_i)$. On its complement,
\eqref{eq:normalized-difference-tail} bounds its absolute value by
$Ck e^{-\sigma R}$. Taking $\varepsilon\to0$ and then $R\to\infty$
contradicts
$\|\phi_\varepsilon\|_\infty+\|\psi_\varepsilon\|_\infty=1$.
\end{proof}

\begin{proof}[Completion of the proof of Theorem~\ref{thm:main-existence}]
Proposition~\ref{prop:con-existence} gives the solution and
\eqref{eq:main-decomposition}. The center convergence and continuity
of the profile scaling show that this solution satisfies
\eqref{eq:main-admissible}. Any other family satisfying that condition
is admissible in the sense of \eqref{eq:unique-class}.
Theorem~\ref{thm:direct-uniqueness} gives the asserted uniqueness for
all sufficiently small values of the common parameter.
\end{proof}

\begin{remark}
A proof by degree would compare the total degree of the position map
with its local degree at every zero. The small spectrum can be used
for this purpose after its relation to the derivative of the position
map has been established. The argument above proves uniqueness from
the local Pohozaev identity before the spectral analysis. Related degree
arguments appear in \cite{BS2002,DFGX2026,ILY2025}.
\end{remark}

\section{The small spectrum of the strongly indefinite Hessian}
\label{sec:potential-spectrum}

We now determine the small eigenvalues of the constructed solution.
Theorem~\ref{thm:direct-uniqueness} then gives the same conclusions for
every admissible family. We use cutoffs of the derivatives of the
solution as approximate translation modes. The equations give an exact quadratic-form identity on these
modes, and the inverse from Lemma~\ref{lem:con-inverse} controls their
orthogonal complement.

\subsection{The Hessian and its domain}

Let $z_\varepsilon=(u_\varepsilon,v_\varepsilon)$ be the solution of
Proposition~\ref{prop:con-existence}, with the centers chosen there.
Write $b_i=a(P_i)$ and $W_i=(U_i,V_i)$. The construction and
$|a(\xi_{i,\varepsilon})-b_i|=O(\varepsilon^4)$ give
\begin{equation}\label{sp:concentration}
 z_\varepsilon(x)=\sum_{i=1}^k
 W_i\left(\frac{x-\xi_{i,\varepsilon}}{\varepsilon}\right)
 +r_\varepsilon(x),\qquad
 \|r_\varepsilon\|_{X_\varepsilon}=O(\varepsilon^2),\qquad
 |\xi_{i,\varepsilon}-P_i|=O(\varepsilon^2).
\end{equation}
The profiles satisfy
\begin{equation}\label{sp:profile-kernel}
 \ker\begin{pmatrix}-qU_i^{q-1}&-\Delta+b_i\\
                    -\Delta+b_i&-pV_i^{p-1}\end{pmatrix}
 =\operatorname{span}\{\partial_jW_i:1\le j\le N\}
 \quad\text{in }H^2(\R^N)^2,
\end{equation}
by Lemma~\ref{lem:autonomous-nondegeneracy}. The equations are
\begin{equation}\label{sp:equation}
 (-\varepsilon^2\Delta+a)u_\varepsilon=v_\varepsilon^p,
 \qquad
 (-\varepsilon^2\Delta+a)v_\varepsilon=u_\varepsilon^q.
\end{equation}
Put
\begin{equation}\label{sp:hessian}
 A_\varepsilon=-\varepsilon^2\Delta+a(x),\qquad
 \mathcal H_\varepsilon=
 \begin{pmatrix}-qu_\varepsilon^{q-1}&A_\varepsilon\\
 A_\varepsilon&-pv_\varepsilon^{p-1}\end{pmatrix},
 \qquad D(\mathcal H_\varepsilon)=H^2(\R^N)^2.
\end{equation}
All projections and inner products in this section are in
$L^2(\R^N)^2$. We use real functions in the computations and complexify
the operators for spectral statements. With
\begin{equation}\label{sp:bounded-multiplier}
 B_\varepsilon=\diag(qu_\varepsilon^{q-1},pv_\varepsilon^{p-1}),
 \qquad \mathcal H_\varepsilon=SA_\varepsilon-B_\varepsilon,
\end{equation}
the scalar operator $A_\varepsilon$ acts on both components.

\begin{lemma}\label{sp:operator-domain}
For each $\varepsilon>0$, the operators $A_\varepsilon$ and
$\mathcal H_\varepsilon$ are self-adjoint on $H^2$ and $H^2{}^2$,
respectively.  The graph norm of $\mathcal H_\varepsilon$ is equivalent
to the $H^2$ norm for fixed $\varepsilon$.  Moreover,
\begin{equation}\label{sp:scalar-lower}
 \|A_\varepsilon f\|_2\geq a_*\|f\|_2,
 \qquad f\in H^2(\mathbb R^N),
\end{equation}
and $A_\varepsilon$ maps $H^2$ onto $L^2$.
\end{lemma}

\begin{proof}
The scalar assertions follow by rescaling Lemma~\ref{lem:con-scalar}.
The Fourier transform and boundedness of the real potential also give
self-adjointness of $A_\varepsilon$ on $H^2$.

By \eqref{sp:concentration},
\begin{equation}\label{sp:uniform-amplitudes}
 \|u_\varepsilon\|_\infty+\|v_\varepsilon\|_\infty\le C,
 \qquad \|r_\varepsilon\|_\infty=O(\varepsilon^2).
\end{equation}
Thus $B_\varepsilon$ is bounded and self-adjoint. The orthogonal change
$(f,g)\mapsto2^{-1/2}(f+g,f-g)$ transforms $SA_\varepsilon$ into
$\diag(A_\varepsilon,-A_\varepsilon)$, so its bounded perturbation
$\mathcal H_\varepsilon$ is self-adjoint on $H^2{}^2$.
Finally,
$-\varepsilon^2\Delta w=S\mathcal H_\varepsilon w-aw+SB_\varepsilon w$
and the Fourier estimate give
\begin{equation}\label{sp:scaled-graph}
 \varepsilon\|\nabla w\|_2+\varepsilon^2\|D^2w\|_2
 \le C\bigl(\|\mathcal H_\varepsilon w\|_2+\|w\|_2\bigr).
\end{equation}
The reverse estimate follows from \eqref{sp:hessian}. This proves
the graph-norm assertion.
\end{proof}

For each profile, recall
\begin{equation}\label{sp:constants}
 M_i=\int U_iV_i,\qquad
 d_i=\int\bigl((\partial_1U_i)^2+(\partial_1V_i)^2\bigr),\qquad
 \kappa_i=M_i/d_i>0.
\end{equation}

\subsection{The translation modes}

Choose fixed real cutoffs $\chi_i$ equal to one near $P_i$, with
pairwise disjoint compact supports, and set
\begin{equation}\label{sp:true-modes}
 s_{j,\varepsilon}=\varepsilon\partial_jz_\varepsilon,\qquad
 T_{ij,\varepsilon}=\chi_i s_{j,\varepsilon},\qquad
 e_{ij,\varepsilon}=\varepsilon^{-N/2}T_{ij,\varepsilon}.
\end{equation}
Let $G_\varepsilon$ be their Gram matrix, and let
$\widehat C_\varepsilon$ have entries
$\langle e_{il,\varepsilon},\mathcal H_\varepsilon e_{hj,\varepsilon}\rangle$.

\begin{lemma}\label{sp:mode-estimates}
The modes belong to $H^2(\R^N)^2$, and
\begin{equation}\label{sp:gram-residual}
 G_\varepsilon\longrightarrow G_0=\diag(d_iI_N)_{i=1}^k,
 \qquad \|\mathcal H_\varepsilon e_{ij,\varepsilon}\|_2
       =O(\varepsilon^2).
\end{equation}
Moreover,
\begin{equation}\label{sp:form-matrix}
 \varepsilon^{-2}\widehat C_\varepsilon
 \longrightarrow\diag(M_iD^2a(P_i))_{i=1}^k.
\end{equation}
\end{lemma}

\begin{proof}
The bounded solution belongs to $H^2$, and the first derivatives of
$v_\varepsilon^p$ and $u_\varepsilon^q$ belong to $L^2$. Thus the
equations give $z_\varepsilon\in H^3$. Differentiating once and
multiplying by $S$ yields
\begin{equation}\label{sp:differentiated-equation}
 \mathcal H_\varepsilon s_{j,\varepsilon}
 =-\varepsilon(\partial_ja)Sz_\varepsilon.
\end{equation}
Only the first derivatives of the powers enter this identity.

If $K$ is a fixed compact set separated from the $P_i$, the decay in
Lemma~\ref{lem:con-tail} and rescaled interior estimates give
\begin{equation}\label{sp:annulus-smallness}
 \varepsilon^{-N/2}\bigl(\|z_\varepsilon\|_{L^2(K)}
 +\varepsilon\|\nabla z_\varepsilon\|_{L^2(K)}
 +\varepsilon^2\|D^2z_\varepsilon\|_{L^2(K)}\bigr)
 =O(e^{-c_K/\varepsilon}).
\end{equation}
For the last term, apply the interior $W^{2,2}$ estimate on balls of
radius $\varepsilon$ and use the zeroth-order bound on slightly larger
balls. Polynomial factors in $\varepsilon^{-1}$ are absorbed by
decreasing $c_K$.

The profile expansion gives
\begin{equation}\label{sp:mode-strong-limit}
 \chi_i(\xi_{i,\varepsilon}+\varepsilon y)
 s_{j,\varepsilon}(\xi_{i,\varepsilon}+\varepsilon y)
 \longrightarrow\partial_jW_i(y)
 \quad\text{in }L^2(\R^N)^2.
\end{equation}
Indeed, the scaled derivative of the correction is $O(\varepsilon^2)$
in $L^2$, and the other peaks have exponentially small contributions
on the cutoff support. Radial symmetry then gives the Gram limit.

Expanding the commutator gives
\begin{equation}\label{sp:mode-commutator}
 \mathcal H_\varepsilon T_{ij,\varepsilon}
 =-\varepsilon\chi_i(\partial_ja)Sz_\varepsilon
 -\varepsilon^2S\bigl((\Delta\chi_i)s_{j,\varepsilon}
          +2\nabla\chi_i\cdot\nabla s_{j,\varepsilon}\bigr).
\end{equation}
The cutoff terms are exponentially small after normalization, by
\eqref{sp:annulus-smallness}. For the first term, use
$|\nabla a(\xi_{i,\varepsilon})|=O(\varepsilon^2)$ and
$|\nabla a(\xi_{i,\varepsilon}+\varepsilon y)|
\le C(\varepsilon^2+\varepsilon|y|)$ on the $i$th profile.
Its polynomial moments are finite. The correction in
\eqref{sp:concentration} contributes $O(\varepsilon^2)$ before the
leading factor $\varepsilon$. This proves the residual estimate.

For real smooth $\chi$ and $f,g\in H^2(\R^N)^2$, integration by parts
gives
\begin{equation}\label{sp:bilinear-cutoff}
 \langle\chi f,\mathcal H_\varepsilon(\chi g)\rangle
 =\tfrac12\bigl(\langle\chi^2f,\mathcal H_\varepsilon g\rangle
       +\langle\mathcal H_\varepsilon f,\chi^2g\rangle\bigr)
       +\varepsilon^2\int|\nabla\chi|^2f\cdot Sg.
\end{equation}
Apply this with $f=s_{l,\varepsilon}$, $g=s_{j,\varepsilon}$ and
$\chi=\chi_i$. Since
$s_{l,\varepsilon}\cdot Sz_\varepsilon
 =\varepsilon\partial_l(u_\varepsilon v_\varepsilon)$,
\eqref{sp:differentiated-equation} and one more integration by parts
show that
\begin{align}\label{sp:cutoff-identity}
 &\langle\chi s_{l,\varepsilon},
               \mathcal H_\varepsilon(\chi s_{j,\varepsilon})\rangle
 \notag\\
 &=\varepsilon^2\int\chi^2(\partial_{lj}a)u_\varepsilon v_\varepsilon
   +\varepsilon^2\int\chi\bigl((\partial_l\chi)(\partial_ja)
          +(\partial_j\chi)(\partial_la)\bigr)u_\varepsilon v_\varepsilon
 \notag\\
 &\quad+\varepsilon^2\int|\nabla\chi|^2
                     s_{l,\varepsilon}\cdot Ss_{j,\varepsilon}.
\end{align}
The last two terms are exponentially small after normalization.
The cutoff solution rescaled about $\xi_{i,\varepsilon}$ converges to
$W_i$ in $L^2$, so the product of its components converges to $U_iV_i$
in $L^1$. Boundedness and continuity of $D^2a$ therefore give
\[
 \varepsilon^{-N}\int\chi_i^2(\partial_{lj}a)
             u_\varepsilon v_\varepsilon
 \longrightarrow M_i\partial_{lj}a(P_i).
\]
Entries associated with different peaks vanish because the cutoffs
have disjoint supports. This proves \eqref{sp:form-matrix}.
\end{proof}

\subsection{The inverse on the orthogonal complement}

Let $E_\varepsilon$ be the span of the modes in \eqref{sp:true-modes}.
Define
\[
 T_\varepsilon h=\sum_{i,j}h_{ij}e_{ij,\varepsilon},\qquad
 J_\varepsilon=T_\varepsilon G_\varepsilon^{-1/2},\qquad
 P_\varepsilon=J_\varepsilon J_\varepsilon^*,\qquad
 Q_\varepsilon=I-P_\varepsilon.
\]
Thus $J_\varepsilon$ is an isometry from $\R^{kN}$ onto
$E_\varepsilon$. On $E_\varepsilon^\perp$, set
\begin{equation}\label{sp:compression-domain}
 D_\varepsilon w=Q_\varepsilon\mathcal H_\varepsilon w,
 \qquad D(D_\varepsilon)=H^2(\R^N)^2\cap E_\varepsilon^\perp.
\end{equation}

\begin{lemma}\label{sp:compressed-inverse}
The operator $D_\varepsilon$ is self-adjoint on the stated domain.
There is $c>0$, independent of small $\varepsilon$, such that
\begin{equation}\label{sp:compressed-bound}
 \|D_\varepsilon w\|_2\ge c\|w\|_2,
 \qquad w\in D(D_\varepsilon).
\end{equation}
Consequently,
\begin{equation}\label{sp:shifted-inverse}
 \|(D_\varepsilon-\lambda)^{-1}\|_{L^2\to L^2}\le2/c,
 \qquad |\lambda|\le c/2,
\end{equation}
and each inverse maps onto $D(D_\varepsilon)$.
\end{lemma}

\begin{proof}
Since $E_\varepsilon\subset H^2$ is finite-dimensional,
$\mathcal H_\varepsilon P_\varepsilon$ is bounded on $L^2$, and
$P_\varepsilon\mathcal H_\varepsilon$ has its bounded adjoint
extension. Removing the bounded self-adjoint off-diagonal operator
$Q_\varepsilon\mathcal H_\varepsilon P_\varepsilon
 +P_\varepsilon\mathcal H_\varepsilon Q_\varepsilon$
from $\mathcal H_\varepsilon$ gives the self-adjoint direct sum of its
finite-dimensional compression and $D_\varepsilon$. This proves
self-adjointness on \eqref{sp:compression-domain}.

We apply \eqref{con-L2-inverse} to compare with the profile
linearization. Let
$\mathcal U_\varepsilon f(y)=\varepsilon^{N/2}f(\varepsilon y)$,
a unitary map on $L^2$, and use the actual centers in the construction.
The profile expansion and the mode definition give
\[
 \|\mathcal U_\varepsilon e_{ij,\varepsilon}-Z_{ij}\|_2
       =O(\varepsilon^2).
\]
Hence, for $w\in E_\varepsilon^\perp$ and
$\widehat w=\mathcal U_\varepsilon w$,
\[
 |(\langle\widehat w,Z_{ij}\rangle)_{i,j}|
       \le C\varepsilon^2\|w\|_2.
\]
The power estimate and \eqref{con-final-size} also yield
\[
 \|\mathcal U_\varepsilon\mathcal H_\varepsilon
             \mathcal U_\varepsilon^{-1}-SL_{\varepsilon,\xi}\|_{L^2\to L^2}
 \le C\|w_{\varepsilon,\xi}\|_\infty^\tau=o(1),
 \qquad \tau=\min\{1,p-1,q-1\}>0.
\]
The difference here is a bounded multiplication operator. Finally,
self-adjointness and \eqref{sp:gram-residual} give
\[
 \|P_\varepsilon\mathcal H_\varepsilon w\|_2
 \le\|\mathcal H_\varepsilon J_\varepsilon\|\|w\|_2
 \le C\varepsilon^2\|w\|_2.
\]
Apply \eqref{con-L2-inverse} to $(\widehat w,0)$. The last three
estimates imply
\[
 \|w\|_2\le\|\widehat w\|_{H^2}
 \le C\|D_\varepsilon w\|_2+o(1)\|w\|_2.
\]
Absorption proves \eqref{sp:compressed-bound}.

This estimate gives a closed range and zero kernel. Self-adjointness
then implies surjectivity. The Neumann series for
$I-\lambda D_\varepsilon^{-1}$ proves
\eqref{sp:shifted-inverse}.
\end{proof}

\subsection{The Schur complement and the eigenvalue count}

\begin{proof}[Proof of Theorem~\ref{thm:main-spectrum}]
Put $m=kN$ and define
\[
 C_\varepsilon=J_\varepsilon^*\mathcal H_\varepsilon J_\varepsilon,
 \qquad R_\varepsilon=Q_\varepsilon\mathcal H_\varepsilon J_\varepsilon.
\]
Under the unitary map $(h,w)\mapsto J_\varepsilon h+w$, the Hessian
has the block representation
\[
 \begin{pmatrix}C_\varepsilon&R_\varepsilon^*\\
                 R_\varepsilon&D_\varepsilon\end{pmatrix},
 \qquad \text{with domain }\R^m\oplus D(D_\varepsilon).
\]
Lemma~\ref{sp:mode-estimates} gives
\begin{equation}\label{sp:blocks}
 C_\varepsilon=G_\varepsilon^{-1/2}\widehat C_\varepsilon
                      G_\varepsilon^{-1/2},\qquad
 \|R_\varepsilon\|=O(\varepsilon^2),\qquad
 \varepsilon^{-2}C_\varepsilon\longrightarrow
 C_0=\diag(\kappa_iD^2a(P_i))_{i=1}^k.
\end{equation}
Fix $0<\gamma<\min\{a_*/2,c/2\}$. For real $|\lambda|\le\gamma$, put
$\mathcal R_\varepsilon(\lambda)=(D_\varepsilon-\lambda)^{-1}$ and
\begin{equation}\label{sp:schur}
 F_\varepsilon(\lambda)=C_\varepsilon-\lambda I
       -R_\varepsilon^*\mathcal R_\varepsilon(\lambda)R_\varepsilon.
\end{equation}
Eliminating the complementary component in
\[
 (C_\varepsilon-\lambda I)h+R_\varepsilon^*w=f,
 \qquad R_\varepsilon h+(D_\varepsilon-\lambda)w=g
\]
gives
\[
 F_\varepsilon(\lambda)h
       =f-R_\varepsilon^*\mathcal R_\varepsilon(\lambda)g,
 \qquad
 w=\mathcal R_\varepsilon(\lambda)(g-R_\varepsilon h).
\]
If $F_\varepsilon(\lambda)$ is invertible, the inverse in these
coordinates is
\begin{equation}\label{sp:block-inverse}
 \begin{pmatrix}
 F^{-1}&-F^{-1}R^*\mathcal R\\
 -\mathcal RRF^{-1}&\mathcal R+\mathcal RRF^{-1}R^*\mathcal R
 \end{pmatrix},
\end{equation}
where subscripts and the argument are suppressed. All blocks are
bounded, and the second row takes values in $D(D_\varepsilon)$.
If the finite matrix is singular, the same elimination gives
\begin{equation}\label{sp:kernel-schur}
 \ker(\mathcal H_\varepsilon-\lambda)
 =\{J_\varepsilon h-\mathcal R_\varepsilon(\lambda)R_\varepsilon h:
                                h\in\ker F_\varepsilon(\lambda)\}.
\end{equation}
This is an isomorphism of kernels, since its two terms are orthogonal
and $J_\varepsilon$ is an isometry. Thus the spectral points in the
chosen interval are precisely the singular parameters of
$F_\varepsilon$, and their eigenspace dimensions agree with the dimension of its kernel.

The resolvent identity gives
$\mathcal R_\varepsilon'(\lambda)=\mathcal R_\varepsilon(\lambda)^2$.
Consequently,
\begin{equation}\label{sp:schur-monotonicity}
 F_\varepsilon'(\lambda)
 =-I-R_\varepsilon^*\mathcal R_\varepsilon(\lambda)^2R_\varepsilon
 \le-I.
\end{equation}
This uses self-adjointness of the inverse, since its square is
nonnegative. Let $\theta_{\varepsilon,1}(\lambda)\le\cdots\le
\theta_{\varepsilon,m}(\lambda)$ be the ordered eigenvalues of this
symmetric matrix. The finite-dimensional min--max formula
\begin{equation}\label{sp:finite-minmax}
 \theta_{\varepsilon,j}(\lambda)
 =\min_{\dim V=j}\max_{0\ne h\in V}
       \frac{\langle F_\varepsilon(\lambda)h,h\rangle}{|h|^2}
\end{equation}
shows that they are continuous and that
\[
 \theta_{\varepsilon,j}(\lambda_2)
 \le\theta_{\varepsilon,j}(\lambda_1)-(\lambda_2-\lambda_1)
 \quad(\lambda_2>\lambda_1).
\]
By \eqref{sp:blocks},
$\|C_\varepsilon\|+(2/c)\|R_\varepsilon\|^2=O(\varepsilon^2)$.
It follows that $F_\varepsilon(-\gamma)$ is positive definite and
$F_\varepsilon(\gamma)$ is negative definite for small $\varepsilon$.
Each $\theta_{\varepsilon,j}$ therefore has exactly one zero,
$\lambda_{\varepsilon,j}$, with
$\lambda_{\varepsilon,1}\le\cdots\le\lambda_{\varepsilon,m}$.
At any common zero $\lambda_*$,
\[
 \#\{j:\lambda_{\varepsilon,j}=\lambda_*\}
 =\dim\ker F_\varepsilon(\lambda_*)
 =\dim\ker(\mathcal H_\varepsilon-\lambda_*).
\]
The zeros form a finite set. They are isolated eigenvalues, and
self-adjointness identifies their multiplicities with these dimensions.
There are exactly $m=kN$ eigenvalues in the interval.

To compute their leading terms, let $c_{\varepsilon,j}$ be the ordered
eigenvalues of $C_\varepsilon$. Uniformly in $|\lambda|\le\gamma$,
\[
 \|F_\varepsilon(\lambda)-(C_\varepsilon-\lambda I)\|
 \le(2/c)\|R_\varepsilon\|^2=O(\varepsilon^4).
\]
The min--max comparison and the zero equation yield
\begin{equation}\label{sp:zero-pairing}
 |\lambda_{\varepsilon,j}-c_{\varepsilon,j}|=O(\varepsilon^4).
\end{equation}
Since $\varepsilon^{-2}C_\varepsilon\to C_0$,
$\varepsilon^{-2}c_{\varepsilon,j}\to\mu_j$. This proves
\eqref{sp:asymptotic}, including repeated limiting eigenvalues.
The matrices $D^2a(P_i)$ are invertible and $\kappa_i>0$, so the
asymptotics determine the signs and exclude zero from the spectrum.
Since $\mathcal H_\varepsilon=S\mathcal L_\varepsilon$ on their common
domain and $S$ is invertible, both kernels are zero.

By Theorem~\ref{thm:direct-uniqueness}, every family satisfying
\eqref{eq:main-admissible} coincides with the constructed solution for
small $\varepsilon$. The same spectral conclusions therefore apply.
\end{proof}

\section{The relative Morse index of the multi-peak solutions}
\label{sec:relative-morse}

Section~\ref{sec:potential-spectrum} counts the small eigenvalues
arising from translations. To compute the relative Morse index, we also need to count the negative directions of each autonomous profile.
We first separate an infinite-dimensional strictly negative part of
the Hessian by completing the square. The remaining scalar operator
is bounded below and has finitely many negative eigenvalues.
The variational characterization \eqref{eq:ground-state-quotient} gives
one negative direction for each limiting profile. Localization gives the
sum of these contributions, and the small eigenvalues give the remaining
negative directions. This yields
$m_{\mathrm{rel}}=k+\sum_i\operatorname{ind}D^2a(P_i)$.

\subsection{Removing the strictly negative component}

We pass to the sum and difference of the components and eliminate
the strictly negative block. We then count the eigenvalues passing
through zero as the nonlinear coefficients increase from zero to
their final values. Zero eigenvalues at the endpoint are recorded
separately.

Fix $\varepsilon>0$ and abbreviate
\begin{equation}\label{mi:coefficients}
 A=-\varepsilon^2\Delta+a(x),\qquad
 b=qu_\varepsilon^{q-1},\qquad c=pv_\varepsilon^{p-1},\qquad
 H_t=\begin{pmatrix}-tb&A\\A&-tc\end{pmatrix}
 \quad(0\leq t\leq1).
\end{equation}
Thus $H_0=SA$, where $S(\phi,\psi)=(\psi,\phi)$ and $A$ acts
componentwise, and $H_1=\mathcal H_\varepsilon$.
Each $H_t$ is self-adjoint on $H^2(\mathbb R^N)^2$ in
$L^2(\mathbb R^N)^2$. The path parameter $t\in[0,1]$ is distinct
from the Sobolev exponent used earlier. The coefficients $b,c$ do
not vary with $t$; they are bounded, strictly positive and tend to
zero at infinity. The same argument applies
to an autonomous profile, with $A=-\Delta+m$, $m>0$.

On the infinite-dimensional subspace of pairs $(f,-f)$, the
Hessian form satisfies
\[
 \langle H_1(f,-f),(f,-f)\rangle
 =-2\|A^{1/2}f\|_2^2-\int(b+c)f^2<0
 \quad(0\ne f\in H^1).
\]
Thus the ordinary index is infinite. The free operator $SA$ acts as
$A$ on pairs $(f,f)$ and as $-A$ on pairs $(f,-f)$. This suggests
using the sum and difference of the components. Put $d_\pm=(b\pm c)/2$
and use the orthogonal matrix
\[
 P=\frac1{\sqrt2}\begin{pmatrix}1&1\\1&-1\end{pmatrix},
 \qquad (h_+,h_-)=P(\phi,\psi).
\]
Since $PSP=\operatorname{diag}(1,-1)$, direct multiplication gives
\[
 PH_tP=
 \begin{pmatrix}A-td_+&-td_-\\-td_-&-A-td_+\end{pmatrix}.
\]
Eliminate the lower row as in \eqref{sp:schur}.
Here the lower block is negative, which changes the sign of the
inverse term in the remaining operator.
Set
\[
 M_t=A+td_+\geq a_*,\qquad
 \mathcal K_t=A-td_++t^2d_-M_t^{-1}d_-,
 \qquad D(\mathcal K_t)=H^2(\mathbb R^N).
\]
For the final operator we write
\begin{equation}\label{mi:scalar-Schur}
 \mathcal K=\mathcal K_1
 =A-d_++d_-(A+d_+)^{-1}d_-.
\end{equation}
All pairings below are quadratic-form pairings on $H^1$ when the
arguments are not in $H^2$.  Completing the square gives
\begin{equation}\label{mi:complete-square}
 \langle H_t(\phi,\psi),(\phi,\psi)\rangle
 =\langle\mathcal K_t h_+,h_+\rangle
  -\bigl\|M_t^{1/2}
       (h_-+tM_t^{-1}d_-h_+)\bigr\|_2^2.
\end{equation}
The triangular change and its inverse preserve $H^1\times H^1$, since
$M_t^{-1}d_-$ maps $L^2$ into $H^2$. For fixed $h_+$, the form is
maximized at $h_-=-tM_t^{-1}d_-h_+$, and its maximum is
$\langle\mathcal K_t h_+,h_+\rangle$. The positive sign of the inverse
term in \eqref{mi:scalar-Schur} follows from the negative block $-M_t$.

Each $\mathcal K_t$ is a bounded self-adjoint perturbation of $A$,
with form domain $H^1$.  The difference is relatively compact.
As in Lemma~\ref{lem:frozen-range}, a bounded multiplier tending
to zero is compact from $H^2$ to $L^2$. This applies to $d_+$
and, since $M_t^{-1}:L^2\to H^2$ is bounded, also proves that
$d_-M_t^{-1}d_-$ is compact on $L^2$.
The essential spectrum of $\mathcal K_t$ is therefore contained
in $[a_*,\infty)$.  Also
$\mathcal K_t\geq A-t\|d_+\|_\infty$ as forms, so it has only
finitely many negative eigenvalues.  Write $n_-(\mathcal K_t)$
for their number, including multiplicity.

The next identity determines the direction of the crossings.
It will also give a negative direction for each autonomous profile.

\begin{lemma}\label{mi:schur-monotonicity-t}
Let $A=-\varepsilon^2\Delta+a(x)$ with $a$ bounded and
$a\geq a_*>0$, and let
$b,c$ be bounded nonnegative functions.  Put $d_\pm=(b\pm c)/2$.
For $0\leq t\leq1$ define
\[
 M_t=A+td_+,\qquad
 \mathcal K_t=A-td_++t^2d_-M_t^{-1}d_-.
\]
For $h\in H^1(\mathbb R^N)$ set
\[
 z_t=-tM_t^{-1}d_-h,\qquad
 \phi_t=\frac{h+z_t}{\sqrt2},\qquad
 \psi_t=\frac{h-z_t}{\sqrt2}.
\]
Then
\begin{equation}\label{mi:schur-form-derivative-t}
 \frac{d}{dt}\langle\mathcal K_t h,h\rangle
 =-\int\bigl(b\phi_t^2+c\psi_t^2\bigr).
\end{equation}
In particular, the forms of $\mathcal K_t$ are nonincreasing in $t$.
If $b,c>0$ almost everywhere, the derivative is strictly negative
for every nonzero $h$.
\end{lemma}

\begin{proof}
Write the form of
$H_t=\left(\begin{smallmatrix}-tb&A\\A&-tc\end{smallmatrix}\right)$
in the coordinates $(h,z)=2^{-1/2}(\phi+\psi,\phi-\psi)$:
\[
 Q_t(h,z)=\langle(A-td_+)h,h\rangle
       -2t\langle d_-h,z\rangle-\langle M_tz,z\rangle.
\]
Since $M_t\geq a_*$, this strictly concave quadratic function of
$z\in H^1$ has the unique maximizer $z_t=-tM_t^{-1}d_-h$, and
\begin{equation}\label{mi:schur-fibre-maximum}
 \langle\mathcal K_t h,h\rangle=Q_t(h,z_t)
       =\max_{z\in H^1}Q_t(h,z).
\end{equation}
The resolvent identity implies that $M_t^{-1}$ is continuously
differentiable in operator norm on $L^2$, with
\[
 \frac{d}{dt}M_t^{-1}=-M_t^{-1}d_+M_t^{-1}.
\]
For fixed $\varepsilon>0$, elliptic estimates and the same
identity give differentiability from $L^2$ to $H^2$. Thus $z_t$
is differentiable in $H^1$. Since $\partial_zQ_t(h,z_t)=0$,
the chain rule leaves only the derivative at fixed $(h,z)$:
\[
 -\int d_+(h^2+z^2)-2\int d_-hz
 =-\int\left[b\left(\frac{h+z}{\sqrt2}\right)^2
             +c\left(\frac{h-z}{\sqrt2}\right)^2\right].
\]
This proves \eqref{mi:schur-form-derivative-t}.  If $h\ne0$, the
pair $(\phi_t,\psi_t)$ is nonzero.  Strict positivity of both
weights makes the integral strictly positive.
\end{proof}

A parameter $t_0\in(0,1)$ is a crossing if
$\ker H_{t_0}\ne\{0\}$, and its multiplicity is
$\dim\ker H_{t_0}$. It is downward if the eigenvalues near zero
are positive just before $t_0$ and negative just after it.
For the path \eqref{mi:coefficients}, define
\begin{equation}\label{mi:definition}
 m_{\mathrm{rel}}(\mathcal H_\varepsilon;SA)
   =\sum_{0<t<1}\dim\ker H_t,\qquad
 \nu(\mathcal H_\varepsilon)=\dim\ker H_1.
\end{equation}
We also write $m_{\mathrm{rel}}(z_\varepsilon)$ when the reference
$SA$ is understood. The next proposition proves that the sum is
finite and all crossings are downward. It also justifies excluding
the endpoint kernel: directions that reach zero at $t=1$ have not
become negative. In the final application this kernel is zero.

\begin{lemma}\label{mi:local-crossing}
      Let $T_s=T_{s_0}+B_s$ be self-adjoint operators on a common
      domain $\mathcal D$ in a real Hilbert space. Assume that
      $B_{s_0}=0$ and $B_s$ is a family of bounded self-adjoint
      operators which is $C^1$ in the operator norm.
      Suppose that zero is an isolated eigenvalue of $T_{s_0}$, and set
      \[
       E=\ker T_{s_0},\qquad r=\dim E<\infty.
      \]
      If
      \[
       \langle B'_{s_0}e,e\rangle<0\qquad(0\ne e\in E),
      \]
      then there exists $\delta>0$ such that, for all $s$ sufficiently
      close to $s_0$, the spectrum of $T_s$ in $(-\delta,\delta)$
      consists of exactly $r$ eigenvalues, counted with multiplicity.
      These eigenvalues are positive for $s<s_0$ and negative for $s>s_0$.
      
      If, in addition, $T_{s_0}$ is bounded below and has finite
      negative index, then for such $s$,
      \begin{equation}\label{mi:local-index-jump}
       n_-(T_s)=
       \begin{cases}
       n_-(T_{s_0}),&s\le s_0,\\
       n_-(T_{s_0})+r,&s>s_0.
       \end{cases}
      \end{equation}
      At an endpoint, the derivative and the conclusions are understood
      from within the parameter interval.
      \end{lemma}

\begin{proof}
Let $P_E$ be the orthogonal projection onto $E$ and $Q_E=I-P_E$.
These projections preserve $\mathcal D$. In the decomposition
$E\oplus E^\perp$, write
\[
 C_s=P_ET_s|_E,\qquad R_s=Q_ET_s|_E,\qquad
 D_s=Q_ET_s|_{\mathcal D\cap E^\perp}.
\]
At $s=s_0$, $C_{s_0}=R_{s_0}=0$ and $D_{s_0}$ is self-adjoint and
invertible. Choose $\eta>0$ so that its spectrum avoids
$[-4\eta,4\eta]$. For nearby $s$, the Neumann series gives
$\|(D_s-\lambda)^{-1}\|\le(2\eta)^{-1}$ for $|\lambda|\le\eta$.
The same block elimination as in \eqref{sp:schur}--\eqref{sp:kernel-schur}
leaves
\begin{equation}\label{mi:hessian-finite-matrix}
 F_s(\lambda)=C_s-\lambda I_E
              -R_s^*(D_s-\lambda)^{-1}R_s.
\end{equation}
The kernels have the same dimension, and
\begin{equation}\label{mi:hessian-finite-derivative}
 \partial_\lambda F_s(\lambda)
 =-I_E-R_s^*(D_s-\lambda)^{-2}R_s\le-I_E.
\end{equation}
For nearby $s$, $F_s(-\eta)$ is positive definite and $F_s(\eta)$
is negative definite. The ordered-eigenvalue argument used in
Theorem~\ref{thm:main-spectrum} therefore gives exactly $r$ eigenvalues of $T_s$
in $(-\eta,\eta)$, counted with multiplicity.

Put $h=s-s_0$ and $\Gamma=P_EB'_{s_0}|_E$. Differentiability gives
\[
 C_s=h\Gamma+o(|h|),\qquad R_s=O(|h|),\qquad
 F_s(0)=h\Gamma+o(|h|).
\]
Since $\Gamma$ is negative definite on the finite-dimensional space $E$,
$F_s(0)$ is positive definite for $h<0$ and negative definite for $h>0$.
The strict decrease in \eqref{mi:hessian-finite-derivative} proves the
asserted signs of its zeros.

For the last assertion, completing the square at $\lambda=0$ gives
\[
 \langle T_s(e+w),e+w\rangle
 =\langle F_s(0)e,e\rangle
  +\langle D_s(w+D_s^{-1}R_se),w+D_s^{-1}R_se\rangle.
\]
The bounded triangular change preserves the form domain and is
invertible. Hence
\begin{equation}\label{mi:finite-inertia-identity}
 n_-(T_s)=n_-(D_s)+n_-(F_s(0)).
\end{equation}
The negative index of $D_s$ is constant near $s_0$. Indeed, its negative
spectral subspace at $s_0$ is finite dimensional, and the positive and
negative spectral subspaces have a common gap from zero. A sufficiently
small bounded perturbation preserves their strict signs. The dimension
argument gives both bounds for the negative index. Thus
$n_-(D_s)=n_-(T_{s_0})$, and the signs of $F_s(0)$ prove
\eqref{mi:local-index-jump}. Zero directions at $s=s_0$ are excluded
from the negative index.
\end{proof}

\begin{proposition}\label{mi:negative-block}
All crossings of $H_t$ are downward and there are finitely many.
The scalar reduction satisfies
\begin{equation}\label{mi:scalar-index}
 m_{\mathrm{rel}}(z_\varepsilon)=n_-(\mathcal K),
 \qquad \nu(\mathcal H_\varepsilon)=\dim\ker\mathcal K.
\end{equation}
\end{proposition}

\begin{proof}
The lower row of $PH_tP(h_+,h_-)=0$ gives
$h_-=-tM_t^{-1}d_-h_+$. Substitution in the upper row gives
$\mathcal K_t h_+=0$. Conversely, this formula recovers a unique
$H^2$ pair from every $h_+\in\ker\mathcal K_t$. Therefore
\begin{equation}\label{mi:kernel-bijection-t}
 \ker H_t\cong\ker\mathcal K_t\qquad(0\le t\le1).
\end{equation}
Both paths have fixed domains and continuously differentiable bounded
perturbations. Their essential spectra avoid zero. For $H_t$, this
follows from the gap of $SA$ and the relative compactness of
$\operatorname{diag}(b,c)$; for $\mathcal K_t$, it follows from its
positive essential-spectrum bound. At each nonzero kernel,
\[
 \langle H_t'(\phi,\psi),(\phi,\psi)\rangle
 =-\int(b\phi^2+c\psi^2)<0,
\]
and Lemma~\ref{mi:schur-monotonicity-t} gives the same strict sign
for $\mathcal K_t'$. Lemma~\ref{mi:local-crossing} therefore applies
to both paths. Their crossings are downward, have the same
multiplicities and are isolated.

The set of crossing parameters is closed in $[0,1]$, since
invertibility persists under a small bounded perturbation. It is
empty near zero because $\mathcal K_0=A\ge a_*$. Isolation also
holds from the left at $t=1$. Compactness of $[0,1]$ now shows that
there are finitely many crossings. Starting with
$n_-(\mathcal K_0)=0$ and applying \eqref{mi:local-index-jump} at each
one yields
\[
 n_-(\mathcal K_1)
 =\sum_{0<t<1}\dim\ker\mathcal K_t
 =\sum_{0<t<1}\dim\ker H_t.
\]
If $t=1$ is a crossing, the modes approaching zero there are still
positive for $t<1$ and are zero at the endpoint. They contribute to
$\dim\ker\mathcal K_1$ and not to $n_-(\mathcal K_1)$.
This proves both identities in \eqref{mi:scalar-index}.
\end{proof}

\subsection{The variational index of the positive autonomous profile}

We now compute the relative Morse index of the limiting profile.
The variational characterization \eqref{eq:profile-best-quotient} gives
one negative direction, and \eqref{eq:profile-kernel} gives the kernel.

We now work with a fixed scaled ground state. To distinguish
the constant mass from the coefficients of the Hessian, write
\[
 A=-\Delta+m,\qquad
 AU=V^p,\quad AV=U^q,\qquad
 b=qU^{q-1},\quad c=pV^{p-1},\quad m>0.
\]
By \eqref{eq:profile-best-quotient}, with $r=(p+1)/p\in(1,2)$ and
$s=q+1>2$, the function $f_0=V^p$ maximizes
$\|A^{-1}f\|_s/\|f\|_r$.

The constrained second variation gives the following quadratic
inequality.

\begin{lemma}\label{mi:profile-variation-direct}
For every real $\psi\in L^2(\mathbb R^N)$ satisfying
\[
 \int cV\psi=0,
\]
let $\phi=A^{-1}(c\psi)$.  Then
\begin{equation}\label{mi:profile-direct-inequality}
 \int b\phi^2\leq\int c\psi^2.
\end{equation}
\end{lemma}

\begin{proof}
The two profile equations, tested against the opposite
components, give
\[
 I_0:=\int U^{q+1}=\int V^{p+1}>0.
\]
First suppose that $\psi\in C_c^\infty$ satisfies the stated
condition and put $h=c\psi$.  On the compact support of $h$,
$V^p$ has a positive minimum and $h$ is bounded.  Hence
$f_0+\theta h>0$ for all sufficiently small $|\theta|$.
Moreover,
\[
 \int f_0^{r-1}h=\int Vh=\int cV\psi=0,
 \qquad
 \int U^{s-1}A^{-1}h
     =\int A^{-1}(U^q)h=\int Vh=0.
\]
The second identity follows from the self-adjointness of
$A^{-1}$ on $L^2$: the exponential profile decay gives
$U^q\in L^2$, and $h\in L^2$ has compact support.

The maximizing property of $f_0$ implies that
\[
 \theta\longmapsto
 \frac1s\log\int|A^{-1}(f_0+\theta h)|^s
       -\frac1r\log\int|f_0+\theta h|^r
\]
has a maximum at zero.  For the denominator, compact support and the positive lower bound
for $f_0$ justify two derivatives under the integral. For the numerator,
$s>2$ and $U,A^{-1}h\in L^s$ give the integrable bound
$(|U|+|A^{-1}h|)^{s-2}|A^{-1}h|^2$ for the second derivative.
Set
\[
 J(\theta)=\int|U+\theta A^{-1}h|^s,\qquad
 L(\theta)=\int|f_0+\theta h|^r.
\]
Then $J(0)=L(0)=I_0$, $J'(0)=L'(0)=0$, and
\[
 J''(0)=s(s-1)\int U^{s-2}(A^{-1}h)^2,\qquad
 L''(0)=r(r-1)\int f_0^{r-2}h^2.
\]
Consequently the maximum at zero gives
\[
 0\geq\left(\frac1s\log J-\frac1r\log L\right)''(0)
 =\frac1{I_0}\left[
 q\int U^{q-1}(A^{-1}h)^2
       -\frac1p\int V^{1-p}h^2\right].
\]
Since $h=pV^{p-1}\psi=c\psi$, this is exactly
\eqref{mi:profile-direct-inequality}.

To pass to arbitrary $\psi\in L^2$, note that $cV=pV^p\in L^2$.
Choose $\chi\in C_c^\infty$ with $\int cV\chi\ne0$.  If
$\widetilde\psi_n\in C_c^\infty$ converges to $\psi$ in $L^2$,
replace it by
\[
 \psi_n=\widetilde\psi_n-
       \frac{\int cV\widetilde\psi_n}{\int cV\chi}\chi.
\]
Then $\int cV\psi_n=0$ and $\psi_n\to\psi$ in $L^2$.
The coefficients $b,c$ are bounded, and $A^{-1}$ is bounded
on $L^2$.  Consequently
$A^{-1}(c\psi_n)\to A^{-1}(c\psi)$ in $L^2$, and both sides
of \eqref{mi:profile-direct-inequality} converge.  This proves the assertion. Only the compactly supported variations
above require a second derivative.
\end{proof}

To apply this inequality to $\mathcal K$, we choose a pair with
sum coordinate $h$ and $A\phi=c\psi$. Its Hessian form is then
$\int c\psi^2-\int b\phi^2$. The Schur form is at least this
value by \eqref{mi:schur-fibre-maximum}.

\begin{lemma}\label{mi:profile-nonnegative-hyperplane}
The form of
\[
 \mathcal K=A-d_++d_-(A+d_+)^{-1}d_-
\]
is nonnegative on a closed subspace of codimension one in
$H^1(\mathbb R^N)$.
\end{lemma}

\begin{proof}
For $h\in H^1$ define
\[
 \phi_h=\sqrt2\,(A+c)^{-1}(ch),\qquad
 \psi_h=\sqrt2 h-\phi_h.
\]
The operator $A+c$ is uniformly positive.  Thus these functions
are well defined, $\phi_h\in H^2$, $\psi_h\in H^1$, and
\[
 A\phi_h=c\psi_h,\qquad
 \frac{\phi_h+\psi_h}{\sqrt2}=h.
\]
Both maps are bounded from $L^2$ to $L^2$, and from $H^1$
to $H^1$.  Define the continuous linear functional
\[
 \ell(h)=\int cV\psi_h.
\]
This functional is not zero.  Indeed, for $h_*=(pU+V)/\sqrt2$,
\[
 (A+c)(pU)=pV^p+pcU=c(V+pU)=\sqrt2\,ch_*.
\]
The uniqueness of the inverse of $A+c$ gives
$\phi_{h_*}=pU$ and $\psi_{h_*}=V$, so that
$\ell(h_*)=p\int V^{p+1}>0$.
Therefore $\ker\ell$ has codimension one in $H^1$.

For $h\in\ker\ell$, Lemma~\ref{mi:profile-variation-direct}
applies to $\psi_h$, since $\phi_h=A^{-1}(c\psi_h)$.
The original Hessian form at this pair is
\[
 2\langle A\phi_h,\psi_h\rangle
          -\int b\phi_h^2-\int c\psi_h^2
       =\int c\psi_h^2-\int b\phi_h^2\geq0.
\]
Set $z_h=(\phi_h-\psi_h)/\sqrt2$ and $M=A+d_+$. Since this
pair has sum coordinate $h$, \eqref{mi:complete-square} gives
\[
 \begin{split}
 \langle\mathcal Kh,h\rangle
 &=\int c\psi_h^2-\int b\phi_h^2
   +\|M^{1/2}(z_h+M^{-1}d_-h)\|_2^2\geq0.
 \end{split}
\]
Thus the form is nonnegative on $\ker\ell$.
\end{proof}

\begin{proposition}\label{mi:profile-index}
For each $m>0$, the scalar Schur operator associated with the
positive autonomous profile has exactly one negative
eigenvalue, which is simple, and
\[
 \ker\mathcal K
 =\operatorname{span}\{\partial_j(U+V):1\leq j\leq N\}.
\]
In particular, its negative index is $1$ and its kernel has dimension $N$.
\end{proposition}

\begin{proof}
The scalar Schur operator is self-adjoint on $H^2$, with form
domain $H^1$, and is bounded below.  Its essential spectrum is
$[m,\infty)$, because the coefficients tend to zero and its
additional inverse term is compact on $L^2$.  The preceding
lemma implies that its negative index is at most one: any
two-dimensional strictly negative form subspace would have a
nonzero intersection with $\ker\ell$.

To produce a negative direction, let
\[
 t_0=(pq)^{-1/2}<1,\qquad
 Z_0=(\sqrt p\,U,\sqrt q\,V),\qquad
 h_0=\frac{\sqrt p\,U+\sqrt q\,V}{\sqrt2}.
\]
A direct substitution into the two profile equations shows that
\[
 \begin{pmatrix}-t_0b&A\\A&-t_0c\end{pmatrix}Z_0=0.
\]
The kernel correspondence in the Schur reduction gives
$\mathcal K_{t_0}h_0=0$.  Since $h_0\ne0$ and $b,c>0$,
Lemma~\ref{mi:schur-monotonicity-t} gives the explicit identity
\[
 \langle\mathcal K_1h_0,h_0\rangle
 =-\int_{t_0}^1\int_{\mathbb R^N}
      (b\phi_t^2+c\psi_t^2)\,dx\,dt<0,
\]
where $(\phi_t,\psi_t)$ is the pair defined in that lemma for the
fixed function $h=h_0$.  Its sum is $\sqrt2h_0\ne0$, so the
integrand has strictly positive integral for every $t$.
The negative index is therefore exactly one.  Since the
negative spectrum is discrete and is counted with multiplicity,
its single negative eigenvalue is simple.

The kernel correspondence sends a vector $(\phi,\psi)$ in the
autonomous Hamiltonian kernel to $(\phi+\psi)/\sqrt2$, and is
bijective because the eliminated block is invertible.  The identity
\eqref{eq:profile-kernel} gives the displayed kernel and its dimension $N$.

\end{proof}

The variational characterization \eqref{eq:ground-state-quotient}
gives $n_-(\mathcal K_m)=1$. The ordinary Morse index of the Hamiltonian
Hessian is infinite. On the diagonal $p=q$, one has $U=V$
and $\mathcal K_m=A-pU^{p-1}$. For unequal powers, the inverse
term remains, and the hyperplane argument above is needed.

\subsection{Localization of the reduced quadratic form}
\label{mi:k-localization-section}

We now localize the form of the scalar operator near each peak.
The inverse term requires a commutator estimate. Once this error is
controlled, the autonomous negative eigenspaces give a lower bound
for the index. Orthogonality to these spaces gives the matching
upper bound.

We keep $\mathcal H_\varepsilon$ fixed and vary the spectral
parameter $\lambda$ in
$(\mathcal H_\varepsilon-\lambda)(\phi,\psi)=0$.
This differs from \eqref{mi:coefficients}, where $t$ changes the
operator. We first choose a small negative $\lambda$ to avoid the
autonomous translation kernels, and then increase it to zero.
For real $|\lambda|<a_*$, put
\begin{equation}\label{mi:k-family}
 \mathcal K_\varepsilon(\lambda)
 =A_\varepsilon-d_{+,\varepsilon}-\lambda
  +d_{-,\varepsilon}
   (A_\varepsilon+d_{+,\varepsilon}+\lambda)^{-1}
    d_{-,\varepsilon},
 \qquad D(\mathcal K_\varepsilon(\lambda))=H^2(\mathbb R^N),
\end{equation}
where
\[
 d_{+,\varepsilon}
   =\frac{qu_\varepsilon^{q-1}+pv_\varepsilon^{p-1}}2,
 \qquad
 d_{-,\varepsilon}
   =\frac{qu_\varepsilon^{q-1}-pv_\varepsilon^{p-1}}2.
\]
For the $i$th autonomous profile define
\[
 d_{+,i}=\frac{qU_i^{q-1}+pV_i^{p-1}}2,
 \qquad
 d_{-,i}=\frac{qU_i^{q-1}-pV_i^{p-1}}2,
 \qquad A_i=-\Delta+b_i,
\]
and let $\mathcal K_i(\lambda)$ be the expression in
\eqref{mi:k-family} with these autonomous coefficients.
These operators are self-adjoint and bounded below, with form
domain $H^1$. The compactness argument following
\eqref{mi:complete-square} applies to their local and inverse terms.
They are relatively compact perturbations of $A_\varepsilon-\lambda$
and $A_i-\lambda$, respectively. Their essential spectra lie in
$[a_*-\lambda,\infty)$, and their negative indices are finite.

The next estimate requires no derivatives of $d_{+,\varepsilon}$
or $d_{-,\varepsilon}$.

\begin{lemma}\label{mi:k-form-localization}
Let $M=-\Delta+V$ on $H^2(\mathbb R^N)$, where $V$ is real and
bounded and $V\geq\delta>0$.  Let $D$ be a bounded real multiplication
operator, and let $W$ be another bounded real function.  Denote by
\[
 k[f]=\int(|\nabla f|^2+W|f|^2)
       +\langle M^{-1}Df,Df\rangle,
 \qquad f\in H^1(\mathbb R^N),
\]
the corresponding quadratic form.  Suppose that real smooth functions
$\chi_0,\ldots,\chi_k$ satisfy
\[
 \sum_{j=0}^k\chi_j^2=1,\qquad
 0\leq\chi_j\leq1,\qquad
 \|\nabla\chi_j\|_\infty\leq C_0R^{-1},\qquad
 \|\Delta\chi_j\|_\infty\leq C_0R^{-2},\qquad R\geq1.
\]
Then
\begin{equation}\label{mi:k-IMS}
 \left|k[f]-\sum_{j=0}^k k[\chi_j f]\right|
       \leq C R^{-2}\|f\|_2^2.
\end{equation}
Here $C$ depends only on $\delta$, $\|D\|_\infty$, $C_0$ and the
number of cutoffs.  In particular it is independent of the positions
of the cutoffs and of the values of $W$.
\end{lemma}

\begin{proof}
Write $G=M^{-1}$.  Testing $Mu=f$ by $u$ gives
\[
 \|Gf\|_2\leq\delta^{-1}\|f\|_2,
 \qquad
 \|\nabla Gf\|_2\leq\delta^{-1/2}\|f\|_2.
\]
Multiplication by each cutoff preserves $H^2$: for $u\in H^2$,
both $\chi u$ and
\[
 \Delta(\chi u)=\chi\Delta u+2\nabla\chi\cdot\nabla u
                              +(\Delta\chi)u
\]
belong to $L^2$, so $\chi u\in H^2$. Put $C_j=[\chi_j,M]=2\nabla\chi_j\cdot\nabla+\Delta\chi_j$.
The resolvent identity gives
\[
 [\chi_j,G]=-GC_jG,\qquad
 \|C_jG\|\le
 2\delta^{-1/2}\|\nabla\chi_j\|_\infty
       +\delta^{-1}\|\Delta\chi_j\|_\infty\le C/R.
\]
Since $C_j$ is skew-symmetric, $GC_j$ extends to the bounded operator
$-(C_jG)^*$. Commuting once more and using
$[\chi_j,C_j]=-2|\nabla\chi_j|^2$ gives
\begin{equation}\label{mi:k-inverse-commutator}
 [\chi_j,[\chi_j,G]]
 =2GC_jGC_jG+2G|\nabla\chi_j|^2G,
 \qquad \|[\chi_j,[\chi_j,G]]\|\le CR^{-2}.
\end{equation}
These identities hold first on $H^2$ and then on $L^2$ by the displayed
operator bounds.
Since $\sum_j\chi_j^2=1$,
\[
 \sum_j\chi_jG\chi_j-G
 =-\frac12\sum_j[\chi_j,[\chi_j,G]].
\]
The multiplication operator $D$ commutes with every cutoff, so the
absolute difference between the inverse term at $f$ and the sum of
its values at $\chi_jf$ is at most $CR^{-2}\|D\|_\infty^2\|f\|_2^2$.
For the differential part, expansion of the
gradients and $\sum_j\chi_j\nabla\chi_j=0$ give the usual identity
\[
 \sum_j\int|\nabla(\chi_j f)|^2
 =\int|\nabla f|^2+
       \int\sum_j|\nabla\chi_j|^2|f|^2.
\]
The local zeroth-order terms sum exactly.  Combining these facts
proves \eqref{mi:k-IMS}.  The computation first applies to smooth
functions and then to $H^1$ by density and continuity of the forms.
\end{proof}

We now count all negative directions at a fixed parameter avoiding the
autonomous kernels.  This also gives a spectral gap about zero, which
will justify moving the parameter after localization.

\begin{lemma}\label{mi:k-split-count}
For the constructed solution, fix a real $\lambda$ with
$|\lambda|<a_*$.  If $0\notin\sigma(\mathcal K_i(\lambda))$ for
every $i$, then, for all sufficiently small $\varepsilon$,
\begin{equation}\label{mi:k-split-count-formula}
 n_-(\mathcal K_\varepsilon(\lambda))
       =\sum_{i=1}^k n_-(\mathcal K_i(\lambda)),
 \qquad 0\notin\sigma(\mathcal K_\varepsilon(\lambda)).
\end{equation}
In fact there is $c_\lambda>0$, independent of sufficiently small
$\varepsilon$, such that
$\sigma(\mathcal K_\varepsilon(\lambda))\cap
(-c_\lambda,c_\lambda)=\varnothing$.
\end{lemma}

\begin{proof}
Use the unitary scaling
$(\mathcal U_\varepsilon f)(y)=\varepsilon^{N/2}f(\varepsilon y)$ and put
$y_{i,\varepsilon}=\xi_{i,\varepsilon}/\varepsilon$.
We work in these scaled coordinates.
Write
\[
 \widehat a_\varepsilon(y)=a(\varepsilon y),\quad
 \widehat d_{\pm,\varepsilon}(y)
       =d_{\pm,\varepsilon}(\varepsilon y),\quad
 \widehat M_\varepsilon
       =-\Delta+\widehat a_\varepsilon
                         +\widehat d_{+,\varepsilon}+\lambda.
\]
The form of the scaled operator is
\begin{equation}\label{mi:k-scaled-form}
 \widehat k_{\varepsilon,\lambda}[f]
 =\int\bigl(|\nabla f|^2+
       (\widehat a_\varepsilon-\widehat d_{+,\varepsilon}
                                      -\lambda)|f|^2\bigr)
 +\langle\widehat M_\varepsilon^{-1}
               \widehat d_{-,\varepsilon}f,
               \widehat d_{-,\varepsilon}f\rangle.
\end{equation}
Set $\delta=a_*-|\lambda|>0$.  The inverse in this formula is
uniformly bounded because $\widehat M_\varepsilon\geq\delta$.
All its potential coefficients are uniformly bounded.

\emph{Local convergence and the exterior bound.}
The construction estimate \eqref{eq:main-decomposition}, exponential
profile decay, and continuity of the positive powers give, on each
fixed ball,
\begin{equation}\label{mi:k-local-weights}
 \widehat d_{\pm,\varepsilon}(y+y_{i,\varepsilon})
       \longrightarrow d_{\pm,i}(y),\qquad
 a(\xi_{i,\varepsilon}+\varepsilon y)\longrightarrow b_i
 \quad\hbox{uniformly}.
\end{equation}
Moreover, there is $\eta(R)\to0$ such that
\begin{equation}\label{mi:k-exterior-weight}
 \limsup_{\varepsilon\to0}
 \sup_{y\notin\bigcup_i B_R(y_{i,\varepsilon})}
          \widehat d_{+,\varepsilon}(y)\leq\eta(R).
\end{equation}
Indeed, outside these balls both profile sums are bounded by
$C e^{-\omega R}$, and the remainder is $o(1)$ uniformly.  Raising
this bound to the positive powers $p-1$ and $q-1$ proves the assertion.
Only uniform continuity on a bounded nonnegative interval is needed;
no derivative of these powers is used.

Let $k_{i,\lambda}$ be the form of $\mathcal K_i(\lambda)$, and let
$(\tau_{i,\varepsilon}f)(y)=f(y-y_{i,\varepsilon})$.  For each fixed
$R$ there is $\omega_R(\varepsilon)\to0$ such that
\begin{equation}\label{mi:k-core-form-convergence}
 \bigl|\widehat k_{\varepsilon,\lambda}
                  [\tau_{i,\varepsilon} f]-k_{i,\lambda}[f]\bigr|
 \leq\omega_R(\varepsilon)\|f\|_2^2
 \quad\hbox{if }\operatorname{supp} f\subset B_{2R}.
\end{equation}
The differential parts agree, and their local potentials converge
uniformly by \eqref{mi:k-local-weights}.  To check the inverse terms,
put
\[
 \widetilde M_{i,\varepsilon}
 =-\Delta+a(\xi_{i,\varepsilon}+\varepsilon y)
       +\widehat d_{+,\varepsilon}(y+y_{i,\varepsilon})+\lambda,
 \qquad
 \widetilde M_i=-\Delta+b_i+d_{+,i}+\lambda.
\]
The difference $\widetilde M_i-\widetilde M_{i,\varepsilon}$ is multiplication by a uniformly bounded function tending pointwise to zero. Thus this difference converges strongly to zero on $L^2$.
The operator $\widetilde M_i^{-1}1_{B_{2R}}$ is compact: its adjoint
$1_{B_{2R}}\widetilde M_i^{-1}$ is compact by local $H^2$ regularity and Rellich
compactness.  Consequently
\[
 \|(\widetilde M_i-\widetilde M_{i,\varepsilon})\widetilde M_i^{-1}1_{B_{2R}}\|
       \longrightarrow0.
\]
This follows by approximating $\widetilde M_i^{-1}1_{B_{2R}}$ in norm by
finite-rank operators; strong convergence is uniform on each of
their finite-dimensional images.
The resolvent identity now gives
\[
 \|(\widetilde M_{i,\varepsilon}^{-1}-\widetilde M_i^{-1})1_{B_{2R}}\|
 \leq\delta^{-1}
  \|(\widetilde M_i-\widetilde M_{i,\varepsilon})\widetilde M_i^{-1}1_{B_{2R}}\|
 \longrightarrow0.
\]
The factors $\widehat d_{-,\varepsilon}(y+y_{i,\varepsilon})$
also converge uniformly to $d_{-,i}$ on $B_{2R}$. Together with the
uniform resolvent bounds, the last display gives convergence in norm
of the inverse terms compressed to $B_{2R}$. This proves
\eqref{mi:k-core-form-convergence} without any limit assumption on
$a(x)$ at infinity.

\emph{A positive form on a subspace of the required codimension.}
Choose a smooth quadratic partition
$\chi_{0,\varepsilon},\ldots,\chi_{k,\varepsilon}$ such that
\[
 \sum_{j=0}^k\chi_{j,\varepsilon}^2=1,\qquad
 \chi_{i,\varepsilon}=1\ \hbox{on }B_R(y_{i,\varepsilon}),\qquad
 \operatorname{supp}\chi_{i,\varepsilon}
       \subset B_{2R}(y_{i,\varepsilon})\quad(i\geq1),
\]
with the derivative bounds of Lemma~\ref{mi:k-form-localization}.
The balls are disjoint for small $\varepsilon$. Choose each core
cutoff as the cosine of a smooth angle increasing from $0$ to
$\pi/2$ on its annulus, and the exterior cutoff as its sine there.
Taking the exterior cutoff equal to one away from the balls gives
the required partition and derivative bounds.

Let $E_i^-$ be the negative spectral subspace of
$\mathcal K_i(\lambda)$ and $m_i=\dim E_i^-$.
There is $\rho>0$ such that
\[
 k_{i,\lambda}[f]\geq\rho\|f\|_2^2
 \quad\hbox{for }f\in H^1\cap(E_i^-)^\perp,
 \qquad 1\leq i\leq k.
\]
This follows from the assumed zero gaps, finite negative spectra and
the positive essential-spectrum lower bounds.  Impose the linear
conditions
\begin{equation}\label{mi:k-local-orthogonality}
 \tau_{i,\varepsilon}^{-1}(\chi_{i,\varepsilon}f)
       \perp E_i^-\quad(1\leq i\leq k),
\end{equation}
and call the resulting subspace $\mathcal Z_\varepsilon$.  Its codimension in
$L^2$ is at most $m:=\sum_i m_i$.
By \eqref{mi:k-exterior-weight}, for sufficiently large fixed $R$
and then small $\varepsilon$, the exterior form satisfies
\[
 \widehat k_{\varepsilon,\lambda}
                [\chi_{0,\varepsilon}f]
 \geq\frac\delta2\|\chi_{0,\varepsilon}f\|_2^2,
\]
because the inverse term in \eqref{mi:k-scaled-form} is nonnegative.
For each core, \eqref{mi:k-core-form-convergence} and
\eqref{mi:k-local-orthogonality} give
\[
 \widehat k_{\varepsilon,\lambda}
                [\chi_{i,\varepsilon}f]
 \geq(\rho-\omega_R(\varepsilon))
          \|\chi_{i,\varepsilon}f\|_2^2.
\]
Apply \eqref{mi:k-IMS} and sum.  First choose $R$ so large that its
$C/R^2$ error is smaller than one quarter of
$\min\{\rho,\delta/2\}$; then choose $\varepsilon$ sufficiently
small.  We obtain a constant $c_+>0$ independent of these
$\varepsilon$ such that
\begin{equation}\label{mi:k-complement-positive}
 \widehat k_{\varepsilon,\lambda}[f]\geq c_+\|f\|_2^2
 \qquad(f\in H^1\cap \mathcal Z_\varepsilon).
\end{equation}
Decrease $c_+$ so that $0<c_+<\delta$, below the essential
spectrum. The codimension bound and
\eqref{mi:k-complement-positive} allow at most $m$ eigenvalues
below $c_+$, as the min--max argument below makes explicit.

\emph{The matching negative subspace.}
If $m=0$, \eqref{mi:k-complement-positive} already proves the result.
Otherwise choose a smooth cutoff $\vartheta_R$ equal to one on
$B_{R/2}$ and supported in $B_R$, and put
$E_i^{-,R}=\vartheta_RE_i^-$. Cutoff convergence in $H^1$ is uniform
on the unit sphere of each finite-dimensional space $E_i^-$.
The Gram matrices therefore converge to the identity and the form
matrices converge to negative definite matrices. Hence, for all
large $R$, these spaces have dimension $m_i$ and
\[
 k_{i,\lambda}[g]\le-\nu\|g\|_2^2
 \qquad(g\in E_i^{-,R}),
\]
where $\nu>0$ is independent of large $R$ and empty spaces are omitted.
Increase $R$ so that the localization error is also smaller than
$\nu/4$.

The translated spaces $\tau_{i,\varepsilon}E_i^{-,R}$ have disjoint
supports and hence form an $m$-dimensional direct sum.  For
$f=\sum_i\tau_{i,\varepsilon}g_i$ in this sum, the partition has
$\chi_{i,\varepsilon}f=\tau_{i,\varepsilon}g_i$ and
$\chi_{0,\varepsilon}f=0$.  Using \eqref{mi:k-IMS} and
\eqref{mi:k-core-form-convergence} gives
\[
 \widehat k_{\varepsilon,\lambda}[f]
 \leq\sum_i k_{i,\lambda}[g_i]
       +(C/R^2+\omega_R(\varepsilon))\|f\|_2^2
 \leq-\frac\nu2\|f\|_2^2
\]
for small $\varepsilon$. For this semibounded operator, define
\[
 \Lambda_j=\inf_{\substack{L\subset H^1\\ \dim L=j}}
        \ \sup_{0\ne f\in L}
        \frac{\widehat k_{\varepsilon,\lambda}[f]}{\|f\|_2^2}.
\]
Below the bottom of the essential spectrum these are precisely the
eigenvalues, repeated according to multiplicity.  The negative trial
space just constructed gives $\Lambda_m\leq-\nu/2$.
For any $(m+1)$-dimensional $L\subset H^1$, the intersection
$L\cap \mathcal Z_\varepsilon$ contains a nonzero vector, so
\eqref{mi:k-complement-positive} gives $\Lambda_{m+1}\geq c_+$.
Thus
\[
 \Lambda_1\leq\cdots\leq\Lambda_m\leq-\nu/2<0<c_+
                           \leq\Lambda_{m+1}.
\]
This proves the exact count and the gap about zero. We may take
$c_\lambda=\frac12\min\{\nu/2,c_+\}$ when $m>0$, and
$c_\lambda=c_+/2$ when $m=0$.  Returning to the original variables proves
\eqref{mi:k-split-count-formula}.
\end{proof}

\subsection{Recovering the position contribution}
\label{mi:k-crossing-section}

The preceding count gives the autonomous contribution at a fixed
negative parameter. We now increase that parameter to zero. Each
negative small eigenvalue found in Section~\ref{sec:potential-spectrum}
adds its multiplicity to the scalar negative index.

For $|\lambda|<a_*$, the constant orthogonal change of components
transforms $\mathcal H_\varepsilon-\lambda$ into
\[
 \begin{pmatrix}
 A_\varepsilon-d_{+,\varepsilon}-\lambda&-d_{-,\varepsilon}\\
 -d_{-,\varepsilon}&-A_\varepsilon-d_{+,\varepsilon}-\lambda
 \end{pmatrix}.
\]
Solving the second row gives the kernel correspondence
\begin{equation}\label{mi:k-lambda-kernel}
 \ker(\mathcal H_\varepsilon-\lambda)
       \cong\ker\mathcal K_\varepsilon(\lambda),
 \qquad
 h_-=-(A_\varepsilon+d_{+,\varepsilon}+\lambda)^{-1}
                                     d_{-,\varepsilon}h_+.
\end{equation}
The correspondence preserves the $H^2$ domains and the kernel
dimensions. The bounded part of \eqref{mi:k-family} is continuously
differentiable in operator norm, and
\begin{equation}\label{mi:k-lambda-derivative}
 \mathcal K_\varepsilon'(\lambda)
 =-I-d_{-,\varepsilon}
       (A_\varepsilon+d_{+,\varepsilon}+\lambda)^{-2}
          d_{-,\varepsilon}\leq-I.
\end{equation}
The essential spectrum stays above $a_*-\lambda>0$. Thus
Lemma~\ref{mi:local-crossing} applies whenever
$\ker\mathcal K_\varepsilon(\lambda_0)\ne\{0\}$. If its dimension is $r$,
then, for sufficiently small $s>0$,
\begin{equation}\label{mi:k-lambda-jump}
 \begin{split}
 n_-(\mathcal K_\varepsilon(\lambda_0-s))
   &=n_-(\mathcal K_\varepsilon(\lambda_0)),\\
 n_-(\mathcal K_\varepsilon(\lambda_0+s))
   &=n_-(\mathcal K_\varepsilon(\lambda_0))+r.
 \end{split}
\end{equation}

\begin{proof}[Proof of Theorem~\ref{thm:main-morse}]
For the constructed solution, put $m_i=n_-(\mathcal K_i(0))$.
The autonomous kernel correspondence gives
$\dim\ker\mathcal K_i(0)=N$.  The nonzero spectrum has a gap
about zero.  The autonomous version of
\eqref{mi:k-lambda-derivative}, and the downward crossing argument
just proved, therefore give a sufficiently small common
$\gamma>0$ such that
\[
 0\notin\sigma(\mathcal K_i(-\gamma)),\qquad
 n_-(\mathcal K_i(-\gamma))=m_i\quad(1\leq i\leq k).
\]
Indeed, \eqref{mi:k-lambda-jump} with $\lambda_0=0$ shows that
the $N$ translation modes are positive at a small negative
parameter, while the complementary negative index is unchanged.
There are finitely many profiles, so we may choose a common
$\gamma<a_*/2$ inside the interval of Theorem~\ref{thm:main-spectrum}.
Lemma~\ref{mi:k-split-count}, applied at the fixed parameter
$-\gamma$, gives
\[
 n_-(\mathcal K_\varepsilon(-\gamma))=\sum_i m_i,
 \qquad 0\notin\sigma(\mathcal K_\varepsilon(-\gamma)).
\]

All $kN$ small eigenvalues in Theorem~\ref{thm:main-spectrum} still lie in
$(-\gamma,\gamma)$ for small $\varepsilon$, since they are
$O(\varepsilon^2)$.  That theorem excludes other spectral
points in this interval.  The nondegenerate Hessians $D^2a(P_i)$
also imply $0\notin\sigma(\mathcal H_\varepsilon)$.
As $\lambda$ increases from $-\gamma$ to zero, the kernel
correspondence \eqref{mi:k-lambda-kernel} therefore identifies
all crossings of $\mathcal K_\varepsilon(\lambda)$ with this
finite collection of negative small eigenvalues of
$\mathcal H_\varepsilon$, preserving their multiplicities.
Neither endpoint is a crossing.  More explicitly, list the distinct
negative eigenvalues in this interval as
$\ell_1<\cdots<\ell_s<0$ and let
$r_j=\dim\ker(\mathcal H_\varepsilon-\ell_j)$.
On each interval between two consecutive $\ell_j$, the scalar
operator is invertible and its negative index is constant: the
positive and negative subspace estimates used above apply at any
point of that interval.  At $\ell_j$, formula
\eqref{mi:k-lambda-jump} gives an increase of $r_j$.
Adding these finitely many increments yields
\[
 n_-(\mathcal K_\varepsilon(0))
 -n_-(\mathcal K_\varepsilon(-\gamma))=\sum_{j=1}^s r_j.
\]
Both endpoint kernels are zero. Theorem~\ref{thm:main-spectrum} identifies
the last sum with $\sum_i\operatorname{ind}D^2a(P_i)$: the
constants $\kappa_i$ are positive and the limiting eigenvalues
are nonzero. Hence
\begin{align}\label{mi:k-additivity}
 n_-(\mathcal K_\varepsilon(0))
 &=\sum_i m_i+
  \#\{\text{eigenvalues of }\mathcal H_\varepsilon
                      \text{ in }(-\gamma,0)\}\notag\\
 &=\sum_i m_i+\sum_i\operatorname{ind}D^2a(P_i).
\end{align}
Proposition~\ref{mi:negative-block} identifies the left side
with the relative Morse index.  The kernel assertion follows
from Theorem~\ref{thm:main-spectrum}, since
$\mathcal H_\varepsilon=S\mathcal L_\varepsilon$ and $S$ is
invertible. Each prescribed profile is a scaled copy of the ground state $W$,
so Proposition~\ref{mi:profile-index} gives $m_i=1$.
Every positive decaying family satisfying \eqref{eq:main-admissible}
coincides with the constructed family for small $\varepsilon$, by
Theorem~\ref{thm:direct-uniqueness}. The same index formula and kernel
conclusion therefore apply to all such families.
\end{proof}

\section{The location and shape of the concentrated peaks}
\label{sec:quantitative-concentration}

We compute the displacement of each center and the first correction to
its profile. The construction gives an $O(\varepsilon^2)$ remainder.
The local Pohozaev identity determines the center displacement, and
the even profile correction yields the positions of the two component
maxima. We also obtain the mass and energy limits.

We retain the centers obtained in Section~\ref{sec:construction} and write
\[
 b_{i,\varepsilon}=a(\xi_{i,\varepsilon}),\qquad b_i=a(P_i),
 \qquad H_i=D^2a(P_i),
\]
\[
 W_{i,\varepsilon}(x)=
 (U_{b_{i,\varepsilon}},V_{b_{i,\varepsilon}})
       \left(\frac{x-\xi_{i,\varepsilon}}{\varepsilon}\right),
 \qquad z_\varepsilon=\sum_iW_{i,\varepsilon}+w_\varepsilon.
\]
The orthogonality conditions in \eqref{con-projected}, after the
change of variables $x=\varepsilon y$, become
\begin{equation}\label{eq:properties-orthogonality}
 \int_{\R^N}w_\varepsilon(x)\cdot
 \partial_j(U_{b_{i,\varepsilon}},V_{b_{i,\varepsilon}})
 \left(\frac{x-\xi_{i,\varepsilon}}{\varepsilon}\right)\,dx=0,
 \qquad 1\le i\le k,\quad 1\le j\le N.
\end{equation}
Equation \eqref{con-final-size}, expressed in physical variables, gives
\begin{equation}\label{eq:properties-construction-bound}
 |\xi_{i,\varepsilon}-P_i|+
       \|w_\varepsilon\|_{X_\varepsilon}\le C\varepsilon^2.
\end{equation}
Since $\nabla a(P_i)=0$ and the second derivatives of $a$ are bounded,
this also gives
\begin{equation}\label{eq:properties-mass-rates}
 |\nabla a(\xi_{i,\varepsilon})|\le C\varepsilon^2,
 \qquad |b_{i,\varepsilon}-b_i|\le C\varepsilon^4.
\end{equation}
The second estimate follows from Taylor expansion at $P_i$, whose linear
term vanishes.  Constants in this section may depend on the fixed
configuration of critical points and on the profile, but not on
$\varepsilon$.

The following lemma shows that the orthogonality conditions determine
the centers locally. The location formulas in this section use these centers.

\begin{lemma}\label{lem:properties-phase-unique}
There is a fixed $\delta>0$ such that, for small $\varepsilon$, the
centers $\xi_\varepsilon$ are the unique configuration
$\eta$ satisfying
$\max_i|\eta_i-\xi_{i,\varepsilon}|<\delta\varepsilon$ and the orthogonality conditions \eqref{eq:properties-orthogonality}, with the remainder defined by subtracting
the profile sum centered at $\eta$.  In particular, any other
centers satisfying these orthogonality conditions and
$\eta_i=P_i+O(\varepsilon^2)$ coincide with
$\xi_{i,\varepsilon}$ for sufficiently small $\varepsilon$.
\end{lemma}

\begin{proof}
For variable centers $\eta$ put
\[
 W_{\eta,i}(x)=W_{a(\eta_i)}
       \left(\frac{x-\eta_i}{\varepsilon}\right),\qquad
 Z_{\eta,ij}(x)=\partial_jW_{a(\eta_i)}
       \left(\frac{x-\eta_i}{\varepsilon}\right),\qquad
 w_\eta=z_\varepsilon-\sum_hW_{\eta,h},
\]
where $W_b=(U_b,V_b)$. Use the scaled pairing
$\langle f,g\rangle_\varepsilon=\varepsilon^{-N}\int f\cdot g$
and its norm $\|\cdot\|_{2,\varepsilon}$. Define
$\Gamma_{\varepsilon,ij}(\eta)=
\langle w_\eta,Z_{\eta,ij}\rangle_\varepsilon$.
The profile estimates make $\Gamma_\varepsilon$ a $C^1$ map. Since
\[
 \partial_{\eta_{h,\ell}}W_{\eta,h}
 =-\varepsilon^{-1}Z_{\eta,h\ell}
  +\partial_\ell a(\eta_h)(\partial_bW_b)_{b=a(\eta_h)}
       \left(\frac{\cdot-\eta_h}{\varepsilon}\right),
\]
differentiating the pairings gives
\begin{equation}\label{eq:properties-phase-jacobian}
 \varepsilon D_\eta\Gamma_\varepsilon
 ={G}_\eta+{E}_{\varepsilon,\eta},\qquad
 \|{E}_{\varepsilon,\eta}\|
 \le C\bigl(\varepsilon+\|w_\eta\|_{2,\varepsilon}\bigr),
\end{equation}
where ${G}_\eta$ is the Gram matrix of the $Z_{\eta,ij}$.
Indeed, mass derivatives carry a factor $\varepsilon$, and the
remaining errors pair $w_\eta$ with uniformly bounded scaled
$L^2$ profile derivatives.

On the convex set $\max_i|\eta_i-\xi_{i,\varepsilon}|
\le\delta\varepsilon$, integration of the profile derivative gives
$\|w_\eta\|_{2,\varepsilon}\le C(\varepsilon^2+\delta)$.
Exponential separation and \eqref{eq:profile-Gram} imply, uniformly there,
${G}_\eta\to\operatorname{diag}(d_iI_N)$, where
$d_i=\int|\partial_1W_{b_i}|^2>0$.
Choosing $\delta$ and then $\varepsilon$ small yields
\[
 h\cdot D_\eta\Gamma_\varepsilon h
 \ge\frac{\min_i d_i}{2\varepsilon}|h|^2.
\]
Integration along line segments proves injectivity, so the known zero
$\xi_\varepsilon$ is unique. Two configurations which are
$P_i+O(\varepsilon^2)$ differ by $o(\varepsilon)$ and hence lie in
this same neighborhood.
\end{proof}

\subsection{The localized position identity}

Lemma~\ref{lem:unique-tails} gives uniform exponential decay of the
constructed solution and its scaled gradient.

\begin{lemma}\label{lem:potential-pure-tail}
Fix $r>0$ so that the balls $B_{3r}(P_i)$ are pairwise disjoint.
Then
\begin{equation}\label{eq:potential-local-force}
 \varepsilon^{-N}\int_{B_r(P_i)}
       (\partial_j a)u_\varepsilon v_\varepsilon
       =O(e^{-c_r/\varepsilon}),\qquad 1\le j\le N.
\end{equation}
The same estimate holds with any fixed smooth cutoff equal to one on
$B_r(P_i)$ and supported in $B_{2r}(P_i)$.
\end{lemma}

\begin{proof}
The constructed family satisfies \eqref{eq:main-admissible}, so
\eqref{eq:unique-single-small-force} gives the ball estimate.
The cutoff integral differs from it only on an annulus at a fixed
positive distance from all centers. The integrand on that annulus is
exponentially small by \eqref{eq:unique-base-tail}. Absorbing
$\varepsilon^{-N}$ into the exponential proves both assertions.
\end{proof}

\subsection{The displacement of the centers}

The position identity contains the product of the two components.
Consequently the relevant moment is that of $U_iV_i$:
\begin{equation}\label{eq:potential-second-moment}
 M_i=\int U_iV_i,\qquad
 \rho_i=\frac{1}{2NM_i}\int |y|^2U_i(y)V_i(y)\,dy.
\end{equation}
These agree with \eqref{eq:main-kappa} and \eqref{eq:main-rho}.
They are finite and positive by positivity and exponential decay.

\begin{proposition}\label{prop:potential-center-shift}
The centers selected by the translation moment conditions satisfy
\begin{equation}\label{eq:potential-center-shift}
 \xi_{i,\varepsilon}
  =P_i-\varepsilon^2\rho_i H_i^{-1}\nabla\Delta a(P_i)
      +o(\varepsilon^2).
\end{equation}
\end{proposition}

\begin{proof}
Fix $i$ and abbreviate $\xi=\xi_{i,\varepsilon}$,
$b=b_{i,\varepsilon}$, $U=U_b$ and $V=V_b$.
Choose a cutoff $\chi$ as in Lemma~\ref{lem:potential-pure-tail}.
The mixed-term estimates in the proof of Lemma~\ref{lem:con-force}
and \eqref{eq:properties-construction-bound} give
\begin{equation}\label{eq:properties-force-split}
 0=\int\chi(\xi+\varepsilon y)
       \partial_j a(\xi+\varepsilon y)UV\,dy+O(\varepsilon^3).
\end{equation}
Indeed, writing $r_\varepsilon=\|w_\varepsilon\|_{X_\varepsilon}
=O(\varepsilon^2)$, the terms containing the correction are bounded by
$C[(|\nabla a(\xi)|+\varepsilon)r_\varepsilon+r_\varepsilon^2]
=O(\varepsilon^3)$, while the other profiles and cutoff terms are
exponentially small.

Taylor's formula and $a\in C_b^4$ give
\begin{align*}
 \partial_j a(\xi+\varepsilon y)
 ={}&\partial_j a(\xi)
 +\varepsilon\sum_\ell\partial_{j\ell}a(\xi)y_\ell\\
 &+\frac{\varepsilon^2}{2}
       \sum_{\ell,m}\partial_{j\ell m}a(\xi)y_\ell y_m
 +\mathcal R_{j,\varepsilon}(y),
\end{align*}
where $|\mathcal R_{j,\varepsilon}(y)|\le C\varepsilon^3|y|^3$.
Replacing the cutoff by one in these profile integrals has an
exponentially small error, since its nonconstant region has
$|y|\ge c/\varepsilon$.  The product $UV$ is radial.  Hence
\[
 \int y_\ell UV=0,\qquad
 \int y_\ell y_mUV=\frac{\delta_{\ell m}}{N}\int|y|^2UV.
\]
The third moment bounds the integrated Taylor remainder by
$C\varepsilon^3$.  With $M_b=\int U_bV_b$, we obtain
\begin{align}\label{eq:properties-profile-force}
 \int\chi(\xi+\varepsilon y)\partial_j a(\xi+\varepsilon y)UV\,dy
 ={}&M_b\partial_j a(\xi)\notag\\
 &+\frac{\varepsilon^2}{2N}
       \left(\int|y|^2U_bV_b\right)\partial_j\Delta a(\xi)
       +O(\varepsilon^3).
\end{align}
Combining \eqref{eq:properties-force-split} and
\eqref{eq:properties-profile-force}, and using $M_b\to M_i>0$, gives
\[
 \nabla a(\xi_{i,\varepsilon})
    =-\varepsilon^2\rho_i\nabla\Delta a(P_i)+o(\varepsilon^2).
\]
Finally,
\[
 \nabla a(\xi_{i,\varepsilon})
 =H_i(\xi_{i,\varepsilon}-P_i)
       +O(|\xi_{i,\varepsilon}-P_i|^2)
 =H_i(\xi_{i,\varepsilon}-P_i)+O(\varepsilon^4).
\]
The invertibility of $H_i$ proves
\eqref{eq:potential-center-shift}.  
\end{proof}

\subsection{The second-order profile correction}

To locate the component maxima more precisely, we compute the
$\varepsilon^2$ correction to each profile. We retain the original
orthogonality conditions to fix translations.

Recall
\[
 L_b=\begin{pmatrix}-\Delta+b&-pV_b^{p-1}\\
                    -qU_b^{q-1}&-\Delta+b\end{pmatrix},
 \qquad Z_{b,j}=(\partial_jU_b,\partial_jV_b).
\]
The corrector $(\Phi_i,\Psi_i)$ is the unique solution of
\begin{equation}\label{eq:potential-core-corrector}
 \begin{cases}
 L_{b_i}(\Phi_i,\Psi_i)
       =-\dfrac12(y^TH_i y)(U_i,V_i),\\[2mm]
 \displaystyle\int(\Phi_i,\Psi_i)\cdot Z_{b_i,j}=0,
       \qquad 1\le j\le N.
 \end{cases}
\end{equation}
Here and below the complement is defined by the $L^2$ pairing with
$Z_{b_i,j}$, whereas the solvability conditions are the pairings with
$SZ_{b_i,j}$.  These are different roles because $L_{b_i}$ itself is
not self-adjoint.

We verify the solvability and parity explicitly.  The forcing
$F_i=-\frac12(y^TH_i y)(U_i,V_i)$ belongs to
$(L^2\cap L^t)^2$ by exponential decay.  It is even under the simultaneous
reflection $y\mapsto-y$.  Every $SZ_{b_i,j}$ is odd under this reflection,
so
\[
 \langle F_i,SZ_{b_i,j}\rangle=0,\qquad 1\le j\le N.
\]
Lemma~\ref{lem:frozen-range} gives a unique solution in
$(H^2\cap W^{2,t})^2$ satisfying the stated orthogonality conditions.
Reflection commutes with $L_{b_i}$ and preserves those conditions.
The reflected solution has the same forcing, so uniqueness gives
\begin{equation}\label{eq:properties-corrector-even}
 (\Phi_i,\Psi_i)(-y)=(\Phi_i,\Psi_i)(y),\qquad
 \nabla\Phi_i(0)=\nabla\Psi_i(0)=0.
\end{equation}
Only simultaneous reflection is needed; when $H_i$ has off-diagonal
entries, separate reflection in each coordinate need not preserve the
forcing.  Local elliptic regularity makes the corrector smooth, since
the profile components are strictly positive on every compact set.

\begin{proposition}\label{prop:potential-core-expansion}
For every fixed $R>0$,
\begin{align}\label{eq:potential-core-expansion}
 \binom{u_\varepsilon}{v_\varepsilon}
       (\xi_{i,\varepsilon}+\varepsilon y)
 =\binom{U_{b_{i,\varepsilon}}}{V_{b_{i,\varepsilon}}}(y)
       +\varepsilon^2\binom{\Phi_i}{\Psi_i}(y)
       +o(\varepsilon^2)
 \quad\text{in }C^2(B_R)^2.
\end{align}
Each component has exactly one critical point in
$B_{R\varepsilon}(\xi_{i,\varepsilon})$ for small $\varepsilon$.
These points are nondegenerate strict local maxima, denoted by
$x_{i,\varepsilon}^u$ and $x_{i,\varepsilon}^v$, and satisfy
\begin{equation}\label{eq:potential-actual-maxima}
 x_{i,\varepsilon}^u=\xi_{i,\varepsilon}+o(\varepsilon^3),
 \qquad
 x_{i,\varepsilon}^v=\xi_{i,\varepsilon}+o(\varepsilon^3).
\end{equation}
Their heights are
\begin{align}\label{eq:potential-height-correction}
 u_\varepsilon(x_{i,\varepsilon}^u)
  &=U_i(0)+\varepsilon^2\Phi_i(0)+o(\varepsilon^2),\notag\\
 v_\varepsilon(x_{i,\varepsilon}^v)
  &=V_i(0)+\varepsilon^2\Psi_i(0)+o(\varepsilon^2).
\end{align}
Moreover, $u_\varepsilon+v_\varepsilon$ has exactly $k$ local maximum
points in $\R^N$, one in each core.
\end{proposition}

\begin{proof}
We first identify the limit of the divided correction.
Fix $i$ and let
\[
 E_{i,\varepsilon}(y)
   =\varepsilon^{-2}w_\varepsilon(\xi_{i,\varepsilon}+\varepsilon y),
 \qquad P_{i,\varepsilon}(y)
   =(U_{b_{i,\varepsilon}},V_{b_{i,\varepsilon}})(y).
\]
Translation and scaling of \eqref{eq:properties-construction-bound} give
\begin{equation}\label{eq:properties-divided-bound}
 \|E_{i,\varepsilon}\|_{H^2}
       +\|E_{i,\varepsilon}\|_{W^{2,t}}\le C.
\end{equation}
Here the correction is obtained by subtracting all profiles, so the
other peaks do not contribute to these global norms.
For any sequence tending to zero, weak compactness and local Sobolev
compactness give a subsequence and $E_i\in(H^2\cap W^{2,t})^2$ such that
\begin{equation}\label{eq:properties-divided-convergence}
 E_{i,\varepsilon}\rightharpoonup E_i
     \text{ in }H^2\text{ and }W^{2,t},\qquad
 E_{i,\varepsilon}\longrightarrow E_i
     \text{ in }C^{1,\alpha}_{\mathrm{loc}},
\end{equation}
for every fixed $0<\alpha<1-N/t$.

On a fixed ball the contribution of the other profiles is
\[
 T_{i,\varepsilon}(y)
   =\sum_{h\ne i}(U_{b_{h,\varepsilon}},V_{b_{h,\varepsilon}})
       \left(y+\frac{\xi_{i,\varepsilon}-\xi_{h,\varepsilon}}
                          {\varepsilon}\right).
\]
For fixed $R$, its derivatives through order two are
$O_R(e^{-c/\varepsilon})$ on $B_{2R}$. The same estimate holds in
$C^{0,\alpha}$: the values and first derivatives decay exponentially,
and the autonomous equations and interior Schauder estimates give the
bound for second derivatives. This uses only the first derivative of
the powers near zero. Put
\[
 \delta_{i,\varepsilon}
       =\varepsilon^2E_{i,\varepsilon}+T_{i,\varepsilon},
 \qquad
 A_{i,\varepsilon}(y)=a(\xi_{i,\varepsilon}+\varepsilon y)
                                      -b_{i,\varepsilon}.
\]
Taylor expansion and \eqref{eq:properties-mass-rates} imply, on every
fixed ball,
\begin{equation}\label{eq:properties-potential-taylor}
 \varepsilon^{-2}A_{i,\varepsilon}(y)
   =\varepsilon^{-1}\nabla a(\xi_{i,\varepsilon})\cdot y
     +\tfrac12 y^TD^2a(\xi_{i,\varepsilon})y+O_R(\varepsilon)
   \longrightarrow\tfrac12 y^TH_i y
\end{equation}
in $C^{0,\alpha}$.  The first term is $O_R(\varepsilon)$, and the
same bound for the Taylor remainder in this norm follows by applying
Taylor's formula also to one derivative with respect to $y$.

Subtract the equation of $P_{i,\varepsilon}$ from the equation of the
rescaled solution.  If $\mathcal N_{i,\varepsilon}(\delta)$ denotes the
pair of power remainders
\[
 \begin{pmatrix}
  (V_{b_{i,\varepsilon}}+\delta_2)^p
      -V_{b_{i,\varepsilon}}^p
      -pV_{b_{i,\varepsilon}}^{p-1}\delta_2\\
  (U_{b_{i,\varepsilon}}+\delta_1)^q
      -U_{b_{i,\varepsilon}}^q
      -qU_{b_{i,\varepsilon}}^{q-1}\delta_1
 \end{pmatrix},
\]
the exact local equation is
\begin{align}\label{eq:properties-divided-equation}
 L_{b_{i,\varepsilon}}E_{i,\varepsilon}
 ={}&-\varepsilon^{-2}A_{i,\varepsilon}P_{i,\varepsilon}
       -A_{i,\varepsilon}E_{i,\varepsilon}
       +\varepsilon^{-2}\mathcal N_{i,\varepsilon}
                                      (\delta_{i,\varepsilon})\notag\\
 &-\varepsilon^{-2}
       (L_{b_{i,\varepsilon}}+A_{i,\varepsilon}I)
                                      T_{i,\varepsilon}.
\end{align}
Here powers are evaluated at the actual positive components.  The formula
is used only on fixed balls, on which both profiles have a positive
lower bound.  The nonlinear remainder divided by $\varepsilon^2$
therefore tends to zero there.  Equations
\eqref{eq:properties-divided-convergence}--\eqref{eq:properties-divided-equation}
show in distributions that
\[
 L_{b_i}E_i=-\tfrac12(y^TH_i y)(U_i,V_i)
     \quad\text{on }\R^N.
\]

The global weak $L^2$ convergence in
\eqref{eq:properties-divided-convergence} and the strong $L^2$
convergence of $Z_{b_{i,\varepsilon},j}$ to $Z_{b_i,j}$ give
$\langle E_i,Z_{b_i,j}\rangle=0$ from
\eqref{eq:properties-orthogonality}.
Uniqueness in \eqref{eq:potential-core-corrector} yields
$E_i=(\Phi_i,\Psi_i)$. This identifies every subsequential limit,
so the local convergence holds for the whole family.

To improve the convergence, fix $R>0$ and
$0<\alpha<\alpha_0<1-N/t$. On $B_{2R}$ the profile components have a
uniform positive lower bound and
$\|\delta_{i,\varepsilon}\|_{C^{1,\alpha_0}}\le C_R\varepsilon^2$.
For either power $f$, Taylor's formula
\[
 f(s+d)-f(s)-f'(s)d
 =d^2\int_0^1(1-\theta)f''(s+\theta d)\,d\theta
\]
is therefore used on a fixed compact subinterval of $(0,\infty)$.
The bounds on $f^{(2)}$ and $f^{(3)}$ there give
\begin{equation}\label{eq:properties-local-nonlinear-remainder}
 \|\mathcal N_{i,\varepsilon}(\delta_{i,\varepsilon})
             \|_{C^{0,\alpha}(B_{2R})}
 \le C_R\|\delta_{i,\varepsilon}\|_{C^{0,\alpha}(B_{2R})}^2
 \le C_R\varepsilon^4.
\end{equation}
By \eqref{eq:properties-potential-taylor} and
\eqref{eq:properties-divided-equation}, this estimate implies
\[
 L_{b_i}Y_{i,\varepsilon}\longrightarrow0
        \quad\text{in }C^{0,\alpha}(B_{2R})^2,
 \qquad Y_{i,\varepsilon}=E_{i,\varepsilon}-(\Phi_i,\Psi_i).
\]
Indeed, the remaining correction product is $O_R(\varepsilon^2)$,
the other profiles contribute $O_R(\varepsilon^{-2}e^{-c/\varepsilon})$,
and replacing $b_{i,\varepsilon}$ by $b_i$ costs
$O_R(\varepsilon^4)$. Componentwise interior Schauder estimates give
\[
 \|Y_{i,\varepsilon}\|_{C^{2,\alpha}(B_R)}
 \le C_R\left(\|L_{b_i}Y_{i,\varepsilon}\|_{C^{0,\alpha}(B_{2R})}
           +\|Y_{i,\varepsilon}\|_{C^{0,\alpha}(B_{2R})}\right)
 \longrightarrow0.
\]
This proves \eqref{eq:potential-core-expansion}, since the other
profiles are exponentially small. No second derivative of the powers
near zero has been used.

We next record the strict inequalities for the radial derivatives,
including the second derivatives at the origin. Put
$h(r)=-U_b'(r)/r$ and $\ell(r)=-V_b'(r)/r$ for $r>0$.
Taylor expansion of the smooth radial profiles at zero shows that
they extend as $C^2$ radial functions on $\R^{N+2}$, with
$h(0)=-U_b''(0)$ and $\ell(0)=-V_b''(0)$.
Both functions are nonnegative by radial monotonicity and are not
identically zero. Differentiating the radial equations gives
\[
 -h''-\frac{N+1}{r}h'+bh=pV_b^{p-1}\ell,
 \qquad
 -\ell''-\frac{N+1}{r}\ell'+b\ell=qU_b^{q-1}h.
\]
Regarding $h$ and $\ell$ as radial functions on $\R^{N+2}$, the strong
maximum principle gives $h,\ell>0$, including at the origin. Thus
$U_b',V_b'<0$ for $r>0$, and
\begin{equation}\label{eq:properties-profile-hessian}
 D^2U_b(0)=U_b''(0)I<0,\qquad
 D^2V_b(0)=V_b''(0)I<0.
\end{equation}

We can now locate the critical points of the solution. Fix $R>0$.  By \eqref{eq:properties-profile-hessian} there is a
$0<\rho<R$ such that each limiting profile is strictly concave on
$\overline{B_\rho}$.  The $C^2$ convergence just proved gives uniform
strict concavity for each rescaled component on this ball.  Their
outward radial derivatives on $\partial B_\rho$ are negative for small
$\varepsilon$, so each component attains its maximum over the closed
ball in its interior.  Strict concavity gives exactly one critical point
there.  On the compact annulus $\rho\le|y|\le R$, the radial derivative
of each limiting profile is bounded away from zero in absolute value.
The local $C^1$ convergence excludes any other critical point in the
annulus.  This proves existence and uniqueness of the component maxima
in the stated core.

Let $s_{i,\varepsilon}^u=(x_{i,\varepsilon}^u-\xi_{i,\varepsilon})/
\varepsilon$ and write
$\widehat u_{i,\varepsilon}(y)
=u_\varepsilon(\xi_{i,\varepsilon}+\varepsilon y)$.
The gradient expansion about the selected center at zero and
\eqref{eq:properties-corrector-even} give
\[
 \nabla\widehat u_{i,\varepsilon}(0)=o(\varepsilon^2).
\]
The mean value formula between $0$ and $s_{i,\varepsilon}^u$ gives
\[
 0=\nabla\widehat u_{i,\varepsilon}(0)
   +\left[\int_0^1D^2\widehat u_{i,\varepsilon}
        (\theta s_{i,\varepsilon}^u)
                   \,d\theta\right]s_{i,\varepsilon}^u.
\]
The matrix in brackets is uniformly negative definite.  Its inverse is
uniformly bounded, so $s_{i,\varepsilon}^u=o(\varepsilon^2)$.
Returning to physical coordinates proves the first assertion in
\eqref{eq:potential-actual-maxima}; the second is identical.
Evaluation of \eqref{eq:potential-core-expansion} at these points gives
\eqref{eq:potential-height-correction}, since
$b_{i,\varepsilon}-b_i=O(\varepsilon^4)$ and the profiles depend
smoothly on the mass near $b_i$.

The argument for one core applies also to the sum of the two components:
its limiting Hessian at zero is negative definite and its limiting radial
derivative is negative away from zero.  Choose now a fixed large $R$
such that, outside $\bigcup_i B_{R\varepsilon}(\xi_{i,\varepsilon})$,
\[
 u_\varepsilon^q+v_\varepsilon^p
       \le\frac{a_*}{2}(u_\varepsilon+v_\varepsilon).
\]
Such an $R$ exists by \eqref{eq:unique-base-tail}, uniformly in small
$\varepsilon$.  The equations imply on this exterior region, including
the core spheres,
\[
 \Delta(u_\varepsilon+v_\varepsilon)
    \ge\frac{a_*}{2\varepsilon^2}(u_\varepsilon+v_\varepsilon)>0.
\]
A $C^2$ function at a local maximum has nonpositive Laplacian, so no
local maximum of the sum lies there.  The disjoint cores each contain
exactly one such maximum.  Their total number is therefore exactly $k$.
\end{proof}

\begin{remark}
The component maxima have the same expansion in
\eqref{eq:potential-actual-maxima}, but they need not coincide exactly.
The argument does not exclude additional critical points of an individual
component outside the peak regions. If the masses $b_i$ differ, the
peak heights may also differ.
\end{remark}

\subsection{Concentration of mass and energy}

We now prove the mass and energy limits. Since the energy density may
change sign, we also control its absolute value. Put
\[
 m(b_i)=\int U_i^{q+1}=\int V_i^{p+1},\qquad
 c_i=I_{b_i}(U_i,V_i).
\]
The equality of the nonlinear integrals follows by testing each
autonomous equation against the opposite component and integrating.

\begin{corollary}\label{cor:potential-concentration-measures}
As finite measures tested against bounded continuous functions,
\begin{align}\label{eq:potential-concentration-measures}
 \varepsilon^{-N}u_\varepsilon^{q+1}\,dx
   &\rightharpoonup\sum_i m(b_i)\delta_{P_i},\notag\\
 \varepsilon^{-N}v_\varepsilon^{p+1}\,dx
   &\rightharpoonup\sum_i m(b_i)\delta_{P_i},\notag\\
 \varepsilon^{-N}u_\varepsilon v_\varepsilon\,dx
   &\rightharpoonup\sum_iM_i\delta_{P_i},\notag\\
 \varepsilon^{-N}e_\varepsilon\,dx
   &\rightharpoonup\sum_i c_i\delta_{P_i}.
\end{align}
The last assertion is convergence of signed measures.  In particular,
$\varepsilon^{-N}I_\varepsilon(u_\varepsilon,v_\varepsilon)
\to\sum_i c_i$.
\end{corollary}

\begin{proof}
Write $\overline u_\varepsilon=\sum_iW_{i,\varepsilon,1}$,
$\overline v_\varepsilon=\sum_iW_{i,\varepsilon,2}$ and
$w_\varepsilon=(\phi_\varepsilon,\psi_\varepsilon)$.
In this proof set
$\|f\|_{1,\varepsilon}=\varepsilon^{-N}\int|f|$.
The profile sums are bounded in scaled $H^1$, whereas their corrections
are $O(\varepsilon^2)$ there. Cauchy--Schwarz gives
\begin{align}\label{eq:properties-product-errors}
 \|u_\varepsilon v_\varepsilon
       -\overline u_\varepsilon\overline v_\varepsilon\|_{1,\varepsilon}
       &\le C\varepsilon^2,\notag\\
 \|\varepsilon^2\nabla u_\varepsilon\cdot\nabla v_\varepsilon
       -\varepsilon^2\nabla\overline u_\varepsilon
                          \cdot\nabla\overline v_\varepsilon
          \|_{1,\varepsilon}&\le C\varepsilon^2.
\end{align}
For $s=p+1,q+1>2$, use
\[
 \bigl||A+B|^s-|A|^s\bigr|
 \le C_s(|A|^{s-1}|B|+|B|^s),\qquad
 \varepsilon^{-N}\int|B|^s
 \le\|B\|_\infty^{s-2}\varepsilon^{-N}\|B\|_2^2.
\]
The profile sum $A$ and correction $B$ are respectively bounded and
$O(\varepsilon^2)$ in $L^\infty$ and scaled $L^2$. Thus
\begin{equation}\label{eq:properties-power-errors}
 \|u_\varepsilon^{q+1}-\overline u_\varepsilon^{q+1}\|_{1,\varepsilon}
 +\|v_\varepsilon^{p+1}-\overline v_\varepsilon^{p+1}\|_{1,\varepsilon}
 \le C\varepsilon^2.
\end{equation}
Interactions of different profiles are exponentially small by
Lemma~\ref{lem:con-basic-estimates}. Taylor expansion also gives
\[
 \varepsilon^{-N}\int|a(x)-b_{i,\varepsilon}|
               W_{i,\varepsilon,1}W_{i,\varepsilon,2}
 \le C\bigl(\varepsilon|\nabla a(\xi_{i,\varepsilon})|
                  +\varepsilon^2\bigr)\le C\varepsilon^2,
\]
using the first two moments of $U_bV_b$.
Consequently, with
\[
 e_b=\nabla U_b\cdot\nabla V_b+bU_bV_b
          -\frac{U_b^{q+1}}{q+1}-\frac{V_b^{p+1}}{p+1},
\]
we have
\begin{equation}\label{eq:properties-energy-L1}
 \left\|e_\varepsilon-\sum_i e_{b_{i,\varepsilon}}
       \left(\frac{\cdot-\xi_{i,\varepsilon}}{\varepsilon}\right)
             \right\|_{1,\varepsilon}\le C\varepsilon^2.
\end{equation}
The other three densities satisfy the analogous estimate with their
individual profile densities.

For $D_b=U_b^{q+1},V_b^{p+1},U_bV_b$ or $e_b$, and bounded continuous
$f$, a change of variables gives
\[
 \varepsilon^{-N}\int f(x)D_{b_{i,\varepsilon}}
       \left(\frac{x-\xi_{i,\varepsilon}}{\varepsilon}\right)dx
 =\int f(\xi_{i,\varepsilon}+\varepsilon y)D_{b_{i,\varepsilon}}(y)\,dy
 \longrightarrow f(P_i)\int D_{b_i}.
\]
Dominated convergence applies because the densities have a common
integrable exponential bound in absolute value. The scaled $L^1$
errors above transfer this limit to the actual densities. They also
bound the total variation and the tails of the signed energy measures.
Taking $f=1$ proves the action limit.
\end{proof}

\section{The number of positive solutions}\label{sec:counting}

We prove that each fixed profile neighborhood contains exactly one
solution and then count the solutions corresponding to the nonempty
subsets of the prescribed critical points.

\subsection{Compactness in fixed profile neighborhoods}

Suppose that $P_1,\ldots,P_m$ are distinct nondegenerate critical points
of $a$, where $m$ is fixed.  Put
\[
 W_i=(U_{a(P_i)},V_{a(P_i)}),\qquad
 \iota_i=\operatorname{ind}D^2a(P_i)\quad(1\le i\le m).
\]
Choose $\rho>0$ so small that the closed balls
$\overline{B_{4\rho}(P_i)}$ are pairwise disjoint and $P_i$ is the only
critical point of $a$ in $\overline{B_{3\rho}(P_i)}$.  Such a choice is
possible because the Hessians are invertible.  For a nonempty subset
$J\subset\{1,\ldots,m\}$ and centers
$\xi=(\xi_i)_{i\in J}$ with $\xi_i\in B_\rho(P_i)$,
write
\begin{equation}\label{ct:profile-sum}
 W_{J,\varepsilon,\xi}(x)
 =\sum_{i\in J}W_i\left(\frac{x-\xi_i}{\varepsilon}\right).
\end{equation}
The profiles in this definition use the fixed masses $a(P_i)$.
We use the fixed profile neighborhoods
$\mathcal U_{\varepsilon,J}(\delta,\rho)$ defined in
\eqref{eq:count-tube}.
Let $\mathcal S_\varepsilon$ be the set of positive classical solutions
of \eqref{eq:system} which tend to zero at infinity.  No Sobolev
closeness or finite-energy assumption is added to this solution class.
The estimates below give these properties when they are needed.

Nondegeneracy implies the following local isolation of the autonomous
profile, modulo translations.

\begin{lemma}\label{ct:autonomous-isolation}
There exists $\delta_{\mathrm{aut}}>0$ with the following property for
every $i\in\{1,\ldots,m\}$.  If a positive classical decaying solution
$Z$ of the autonomous system at mass $a(P_i)$ satisfies
$\|Z-W_i\|_\infty<\delta_{\mathrm{aut}}$, then
\[
 Z(y)=W_i(y-\tau)
\]
for a unique $\tau\in\mathbb R^N$.  The number
$\delta_{\mathrm{aut}}$ can be chosen so that $|\tau|<1$.
\end{lemma}

\begin{proof}
Fix one profile $W=W_i$ and set $b_0=a(P_i)$. If isolation failed,
there would be positive decaying autonomous solutions $Z_n$ which
are not translates of $W$ but satisfy $\|Z_n-W\|_\infty\to0$.
The exterior comparison and gradient estimates in
Lemma~\ref{lem:unique-tails}, applied at $\varepsilon=1$ and mass
$b_0$, give $Z_n+|\nabla Z_n|\le Ce^{-\sigma|y|}$ componentwise.
Indeed, uniform closeness to $W$ makes both components small outside
one fixed ball, where
$(-\Delta+b_0/2)(Z_{n,1}+Z_{n,2})\le0$; the boundary values are
uniformly bounded. This argument uses decay, not an initial Sobolev
assumption. Dominated convergence gives $Z_n-W\to0$ in $L^2\cap L^t$.
Subtracting the equations and using the Lipschitz bound of the powers
on their common bounded range then gives convergence in $X$ by the
scalar inverse estimate.

Consider
\[
 \Phi_j(Z,\tau)=\int[Z-W(\cdot-\tau)]\cdot
                          \partial_jW(\cdot-\tau).
\]
At $(W,0)$, its derivative in $\tau$ is the positive Gram matrix of
the translations. The implicit function theorem gives $\tau_n\to0$
such that $w_n=Z_n(\cdot+\tau_n)-W$ has zero translation moments
and tends to zero in $X$. It satisfies
$L_{b_0}w_n=\mathcal R_W(w_n)$, so Lemmas~\ref{lem:frozen-range}
and \ref{lem:power-map} yield
\[
 \|w_n\|_X\le C\|\mathcal R_W(w_n)\|_Y
 \le C\|w_n\|_X^{1+\vartheta},\qquad
 \vartheta=\min\{1,p-1,q-1\}>0.
\]
The compatibility condition is automatic from the equation. For large
$n$, this forces $w_n=0$, a contradiction.

The translate is unique because $W$ has a unique maximum at its center.
For $|\tau|\ge1$,
\[
 \|W(\cdot-\tau)-W\|_\infty
 \ge U_{b_0}(0)-U_{b_0}(1)>0.
\]
Decrease the neighborhood radius below this number, and take the
minimum over the finitely many profiles. This also gives $|\tau|<1$.
\end{proof}

The next lemma improves the fixed error in \eqref{eq:count-tube}
to an error tending to zero along a subsequence.

\begin{lemma}\label{ct:fixed-tube-compactness}
There exists $\delta_0>0$, independent of $\varepsilon$ and of the
nonempty subset $J$, such that the following holds.
Suppose $\varepsilon_n\to0$ and
$z_n\in\mathcal S_{\varepsilon_n}\cap
\mathcal U_{\varepsilon_n,J}(\delta_0,\rho)$.
Then, after passing to a subsequence, there are centers
$\eta_{i,n}\to P_i$, $i\in J$, such that
\begin{equation}\label{ct:automatic-admissibility}
 \left\|z_n-\sum_{i\in J}
 W_i\left(\frac{\cdot-\eta_{i,n}}{\varepsilon_n}\right)
 \right\|_\infty\longrightarrow0.
\end{equation}
\end{lemma}

\begin{proof}
All choices below can be made for the complete finite collection of
profiles, hence uniformly for its subsets.  Choose $h>0$ with
$h^{p-1},h^{q-1}\le a_*/2$, and choose
\begin{equation}\label{ct:radius-choice}
 0<\delta_0<\min\left\{h/2,\ \delta_{\mathrm{aut}}/2,
                \tfrac14\min_i U_i(0)\right\}.
\end{equation}
Let $\xi_{i,n}$ be any centers in the neighborhood definition.  Profile decay gives a
fixed large $R_0$ such that the sum of the profile components is at
most $h/2$ outside
$\bigcup_{i\in J}B_{R_0\varepsilon_n}(\xi_{i,n})$.
The tube bound therefore makes both components of $z_n$ at most $h$
in that exterior region.  Its components are also uniformly bounded
everywhere.

We first obtain uniform decay and local limits. For $T_n=u_n+v_n$,
the equations give
\[
 (-\varepsilon_n^2\Delta+a_*/2)T_n\le0
 \quad\text{outside the stated cores}.
\]
The supersolutions used in Lemma~\ref{lem:unique-tails}, with boundary
values uniformly bounded on the core spheres, yield
\begin{equation}\label{ct:uniform-tail}
 u_n+v_n+\varepsilon_n(|\nabla u_n|+|\nabla v_n|)
 \le C\sum_{i\in J}
       e^{-\sigma|x-\xi_{i,n}|/\varepsilon_n}.
\end{equation}
The constants depend on the fixed profiles, $\rho$ and $\delta_0$,
not on $n$ or the chosen solutions.  This use of exterior comparison
requires only boundedness and decay; finite energy has not been
assumed.  The gradient estimate follows from interior estimates on
balls of radius $\varepsilon_n$.

Pass to a subsequence with $\xi_{i,n}\to Q_i\in
\overline{B_\rho(P_i)}$ for every $i\in J$.  After scaling about
$\xi_{i,n}$, the uniformly bounded equations give local $W^{2,t}$
bounds and then local $C^{1,\alpha}$ bounds.  The power maps have
bounded first derivatives on the common bounded range; Schauder
estimates therefore give subsequential $C^2$ convergence on compact
sets to a nonnegative pair $Z_i=(\widehat U_i,\widehat V_i)$ solving
\[
 (-\Delta+a(Q_i))\widehat U_i=\widehat V_i^p,
 \qquad
 (-\Delta+a(Q_i))\widehat V_i=\widehat U_i^q.
\]
The contribution from the same center in \eqref{ct:uniform-tail} passes to the bound
$\widehat U_i+\widehat V_i\le C e^{-\sigma|y|}$; contributions of
other centers disappear on every fixed ball.
The tube bound also gives
\begin{equation}\label{ct:autonomous-close}
 \|Z_i-W_i\|_\infty\le\delta_0.
\end{equation}
In particular $\widehat U_i(0)>0$.  The equations and the strong
maximum principle imply that both components are strictly positive.

To locate the limiting centers, use the local Pohozaev identity
\eqref{eq:unique-local-pohozaev} on the fixed physical ball
$B_{2\rho}(P_i)$.  Its boundary stays a fixed positive distance from
all centers in the neighborhood definition.  Estimate \eqref{ct:uniform-tail} makes the
flux exponentially small, so
\[
 \varepsilon_n^{-N}\int_{B_{2\rho}(P_i)}
     (\partial_j a)u_nv_n\longrightarrow0.
\]
In coordinates $x=\xi_{i,n}+\varepsilon_n y$, local convergence
and the exponential bound give the limit
\[
 (\partial_j a)(Q_i)\int_{\mathbb R^N}
                  \widehat U_i\widehat V_i=0.
\]
For this limit, the contribution near the $i$th center is bounded by
an integrable exponential in $y$. Contributions from the other centers
are $O(e^{-c/\varepsilon_n})$ on the physical ball, and remain
exponentially small after integration and rescaling. Truncate in $y$,
pass to the local limit, and then let the truncation radius tend to infinity.
The product integral is positive.  Hence $\nabla a(Q_i)=0$.
By the choice of $\rho$, $Q_i=P_i$.

The limit $Z_i$ now has mass $a(P_i)$.  By
\eqref{ct:autonomous-close}, \eqref{ct:radius-choice} and
Lemma~\ref{ct:autonomous-isolation},
\begin{equation}\label{ct:core-translates}
 Z_i(y)=W_i(y-\tau_i),\qquad |\tau_i|<1.
\end{equation}
Here the local isolation of $W_i$ suffices, since the limit $Z_i$ is
already close to $W_i$.

Set $\eta_{i,n}=\xi_{i,n}+\varepsilon_n\tau_i$. Then
$\eta_{i,n}\to P_i$ and
\[
 z_n(\eta_{i,n}+\varepsilon_n y)\longrightarrow W_i(y)
 \quad\text{locally uniformly}.
\]
Since the shifts $\tau_i$ are bounded, the same centers can be used in
\eqref{ct:uniform-tail}, after increasing its constant. On each fixed
scaled core, the error in \eqref{ct:automatic-admissibility} tends
uniformly to zero. On the complement of these cores, the solution and
the profile sum are bounded by $Ce^{-\sigma_1R}$ for some
$\sigma_1>0$, where $R$ is the scaled core radius. Letting first
$n\to\infty$ and then $R\to\infty$ proves the assertion.
\end{proof}

\subsection{The exact count}

\begin{proof}[Proof of Theorem~\ref{thm:main-count}]
Choose $\rho$ as above and take $\delta=\delta_0$ as in
Lemma~\ref{ct:fixed-tube-compactness}, decreasing $\delta_0$ if necessary.
For each nonempty subset $J$, the construction gives a solution
$z_{J,\varepsilon}$ whose centers tend to $P_i$, $i\in J$.
The varying-mass profiles converge uniformly to the fixed profiles
in \eqref{ct:profile-sum}, and the correction tends uniformly to zero.
Thus $z_{J,\varepsilon}$ belongs to its fixed neighborhood for all
sufficiently small $\varepsilon$.

We next prove that the neighborhood contains only this solution.
If this failed for arbitrarily small $\varepsilon$, there would be
$\varepsilon_n\to0$ and solutions
\[
 z_n\in\mathcal S_{\varepsilon_n}\cap
       \mathcal U_{\varepsilon_n,J}(\delta_0,\rho),
 \qquad z_n\ne z_{J,\varepsilon_n}.
\]
By Lemma~\ref{ct:fixed-tube-compactness}, after passing to a subsequence,
$(z_n)$ is admissible for the profiles indexed by $J$. The constructed sequence
$(z_{J,\varepsilon_n})$ is admissible for the same profiles.
Theorem~\ref{thm:direct-uniqueness} gives
$z_n=z_{J,\varepsilon_n}$ for large $n$, a contradiction.
There are finitely many nonempty subsets, so one positive parameter
threshold works for all of them.

The neighborhoods are pairwise disjoint for small $\varepsilon$.
Indeed, suppose that a function belongs to the neighborhoods for
$J\ne J'$. After interchanging the sets, choose $i\in J\setminus J'$
and evaluate the first component at a center
$\xi_i\in B_\rho(P_i)$ for the first representation. Its value is at
least $U_i(0)-\delta_0$. The centers in the second representation
are in different fixed balls, so its value is at most
$\delta_0+Ce^{-c/\varepsilon}$. This contradicts
$2\delta_0<\min_iU_i(0)$ for small $\varepsilon$.

The solution in each neighborhood is the constructed one.
Theorems~\ref{thm:main-spectrum} and \ref{thm:main-morse} give its
nondegeneracy and the stated relative Morse index formula.
Proposition~\ref{prop:potential-core-expansion} gives exactly $|J|$
local maximum points of $u_{J,\varepsilon}+v_{J,\varepsilon}$.
Since there are finitely many subsets, all these conclusions hold
with a common threshold. There are $2^m-1$ nonempty subsets and
$\binom{m}{k}$ subsets of cardinality $k$, which proves the counts.
\end{proof}

\section{Positive solutions under an energy bound}
\label{sec:positive-class}

We replace profile closeness by nonvanishing in fixed balls and the
upper bound for the scaled nonlinear mass in \eqref{eq:positive-class}.
By \cite[Theorem~1.2(1)]{CR2026}, every positive $H^1(\R^N)^2$
solution of \eqref{eq:autonomous} is a common translate of $(U,V)$.
The exponents are related as in the proof of
Lemma~\ref{lem:autonomous-nondegeneracy}. After scaling, every positive
$H^1(\R^N)^2$ solution $(\widetilde U,\widetilde V)$ of
\eqref{eq:frozen-equations} with mass $b$ satisfies
\begin{equation}\label{eq:positive-autonomous-classification}
 (\widetilde U,\widetilde V)(y)=(U_b,V_b)(y-\tau)
 \quad\text{for some }\tau\in\R^N.
\end{equation}
We use this fact to identify the limiting profiles. The mass bound
rules out additional peaks. Uniform decay then allows us to locate
the centers by the local Pohozaev identity.

\subsection{Positive mass and uniform bounds}

Recall the masses in \eqref{eq:positive-mass}, and write
$W_b=(U_b,V_b)$ for the positive autonomous profile.
Before using the action or the mass bound, we verify that their integrals
are finite.  The initial decay argument is for one fixed solution; its
constants need not yet be uniform in the family.

\begin{lemma}\label{pc:finite-energy}
Every positive classical decaying solution $z=(u,v)$ of \eqref{eq:system} belongs
to $H^1(\R^N)^2$, has finite nonlinear integrals, and satisfies
\begin{equation}\label{pc:testing-identity}
 \int_{\R^N}(\varepsilon^2\nabla u\cdot\nabla v+a(x)uv)
 =\int_{\R^N}u^{q+1}=\int_{\R^N}v^{p+1}.
\end{equation}
For the action in \eqref{eq:energy} and
$\theta=1-1/(p+1)-1/(q+1)>0$, it follows that
\begin{equation}\label{pc:action-identity}
 \varepsilon^{-N}I_\varepsilon(z)
       =\theta\mathcal Q_\varepsilon(z).
\end{equation}
\end{lemma}

\begin{proof}
Fix the solution and choose $h>0$ with
$h^{p-1},h^{q-1}\le a_*/2$.  Since both components tend to zero at
infinity, there is $R>0$ such that $u,v\le h$ on $\R^N\setminus B_R$.
Writing $w=u+v$, the equations give
\[
 (-\varepsilon^2\Delta+a(x))w=u^q+v^p\le(a_*/2)w,
 \qquad
 (-\varepsilon^2\Delta+a_*/2)w\le0
       \quad\hbox{on }\R^N\setminus\overline{B_R}.
\]
The comparison argument in Lemma~\ref{lem:unique-exterior-comparison}
gives exponential decay of $u+v$. Here the constants may depend on the
fixed solution and on $\varepsilon$. Interior estimates on balls of
radius $\varepsilon$ give the same decay for the first derivatives,
as in Lemma~\ref{lem:unique-tails}. Hence $u,v\in H^1$ and both
nonlinear powers are integrable.

Let $\chi_R$ be a smooth cutoff equal to one on $B_R$, supported in
$B_{2R}$, with $|\nabla\chi_R|\le C/R$.  Test the first equation with
$\chi_Rv$ and the second with $\chi_Ru$.  The respective cutoff errors
are
\[
 \varepsilon^2\int v\nabla u\cdot\nabla\chi_R,
 \qquad
 \varepsilon^2\int u\nabla v\cdot\nabla\chi_R.
\]
Both tend to zero by the integrability just established.  Letting
$R\to\infty$ proves \eqref{pc:testing-identity}.  Substitution into
\eqref{eq:energy} yields \eqref{pc:action-identity}.
\end{proof}

The next lemma derives a uniform amplitude bound from the mass bound
and strict subcriticality.

\begin{lemma}\label{pc:uniform-bound}
For every $C_0>0$ there is $C>0$ such that every positive classical
solution of \eqref{eq:system} with $0<\varepsilon\le1$ and
$\mathcal Q_\varepsilon(z)\le C_0$ satisfies
\begin{equation}\label{pc:uniform-bound-estimate}
 \|u\|_\infty+\|v\|_\infty
       +\varepsilon\|\nabla u\|_\infty
       +\varepsilon\|\nabla v\|_\infty\le C.
\end{equation}
The constant depends only on $C_0,N,p,q$ and $a^*$.
\end{lemma}

\begin{proof}
Recall $\alpha=(p+1)/(pq-1)$ and $\beta=(q+1)/(pq-1)$ from
Section~\ref{sec:profiles}.  Put
\[
 d=2\alpha(q+1)=2\beta(p+1)
       =\frac{2(p+1)(q+1)}{pq-1}>N,
 \qquad
 M=\max\{\|u\|_\infty^{1/(2\alpha)},
                  \|v\|_\infty^{1/(2\beta)}\}.
\]
The inequality $d>N$ is equivalent to the strict subcritical condition
in \eqref{eq:exponents}.  If $M\le1$, the amplitude bound is immediate.
Otherwise let $\tau=M^{-1}<1$.  The component realizing $M$ attains its
maximum, since it is continuous and tends to zero at infinity.  Choose
such a point $x_0$ and set
\[
 \widetilde u(y)=\tau^{2\alpha}u(x_0+\varepsilon\tau y),
 \qquad
 \widetilde v(y)=\tau^{2\beta}v(x_0+\varepsilon\tau y).
\]
Then $0<\widetilde u,\widetilde v\le1$ and
$\max\{\widetilde u(0),\widetilde v(0)\}=1$.  Since
$\alpha+1=p\beta$ and $\beta+1=q\alpha$, the rescaled equations are
\[
 \begin{cases}
 -\Delta\widetilde u+
       \tau^2 a(x_0+\varepsilon\tau y)\widetilde u
          =\widetilde v^p,\\
 -\Delta\widetilde v+
       \tau^2 a(x_0+\varepsilon\tau y)\widetilde v
          =\widetilde u^q.
 \end{cases}
\]
For any fixed $r>N$, interior $W^{2,r}$ estimates on $B_2$ and the
Sobolev embedding bound both gradients uniformly on $B_1$.  Hence there
are $r_0,c_0>0$, independent of $\varepsilon$, $\tau$ and $x_0$, such that
\[
 \int_{B_{r_0}}
       (\widetilde u^{q+1}+\widetilde v^{p+1})\ge c_0.
\]
Indeed, the component equal to one at zero is at least $1/2$ on a fixed
small ball.  On the other hand, changing variables and using
\eqref{pc:testing-identity} gives
\[
 \int_{B_{r_0}}
       (\widetilde u^{q+1}+\widetilde v^{p+1})
 =\tau^{d-N}\varepsilon^{-N}
       \int_{B_{\varepsilon\tau r_0}(x_0)}(u^{q+1}+v^{p+1})
 \le2C_0\tau^{d-N}.
\]
Thus $M^{d-N}\le2C_0/c_0$, proving the amplitude bound.  Finally,
rescale the original equations on $B_{2\varepsilon}(x)$, for arbitrary
$x$.  Their right-hand sides are now uniformly bounded, so the same
interior estimate bounds $\varepsilon\nabla u$ and
$\varepsilon\nabla v$ uniformly.
\end{proof}

We next justify the use of the autonomous classification for limits of
arbitrary positive solutions.  The nonvanishing points in this lemma
may escape to infinity.

\begin{lemma}\label{pc:autonomous-limits}
Let $\varepsilon_n\downarrow0$ and let $z_n=(u_n,v_n)$ be positive
classical solutions of \eqref{eq:system} with
$\mathcal Q_{\varepsilon_n}(z_n)\le C_0$.  Suppose
$u_n(x_n)+v_n(x_n)\ge\delta_1>0$.  After taking a subsequence with
$a(x_n)\to b\in[a_*,a^*]$, a further subsequence satisfies
\begin{equation}\label{pc:limiting-profile}
 z_n(x_n+\varepsilon_n\,\cdot)
       \longrightarrow W_b(\cdot-\zeta)
       \quad\hbox{in }C^2_{\rm loc}(\R^N)^2
\end{equation}
for some $\zeta\in\R^N$, where $W_b=(U_b,V_b)$.  This limit belongs to
$H^1(\R^N)^2$ and has nonlinear mass $m(b)\ge m_*$.
\end{lemma}

\begin{proof}
Lemma~\ref{pc:uniform-bound} bounds the rescaled solutions and their
first derivatives.  For every fixed $R$,
\[
 \sup_{|y|\le R}|a(x_n+\varepsilon_n y)-b|
 \le \varepsilon_n R\|\nabla a\|_\infty+|a(x_n)-b|\longrightarrow0.
\]
Interior elliptic estimates and the local Lipschitz continuity of the
powers on bounded nonnegative intervals give uniform local
$C^{2,\nu}$ bounds for some $\nu>0$.  A diagonal subsequence therefore
converges in $C^2_{\rm loc}$ to a nonnegative autonomous solution
$(\widehat U,\widehat V)$.  Its value at zero has sum at least
$\delta_1$, so it is nonzero.  The two equations and the strong maximum
principle imply that both components are strictly positive.

Changing variables on each fixed ball and then increasing its radius
shows
\[
 \int\widehat U^{q+1}\le C_0,
 \qquad \int\widehat V^{p+1}\le C_0.
\]
The limit has globally bounded first derivatives.  If
$\widehat U+\widehat V$ did not tend to zero at infinity, there would be
points $y_j\to\infty$ where this sum stays above a positive constant.
The gradient bound would give disjoint balls of a fixed radius on each
of which one component stays above another positive constant.  Summing
their nonlinear integrals would contradict the preceding finite bound.
Consequently both components tend to zero.

Apply the barrier and cutoff proof of Lemma~\ref{pc:finite-energy} to
the autonomous equation with $\varepsilon=1$ and constant potential
$b$.  It gives exponential decay, membership in $H^1(\R^N)^2$ and
equality of the two nonlinear integrals.  Thus the limit belongs to the
class covered by \eqref{eq:positive-autonomous-classification}, which
identifies it as a translate of $W_b$.  Its mass is therefore exactly
$m(b)$.  The lower bound $m(b)\ge m_*$ follows from the definition in
\eqref{eq:positive-mass}.
\end{proof}

\subsection{Compactness in fixed physical neighborhoods}

Fix the $m$ critical points and the numbers $\delta,\kappa$ in
Theorem~\ref{thm:main-positive}.  We choose the physical radius before
letting $\varepsilon$ tend to zero.  With
$d=2(p+1)(q+1)/(pq-1)>N$, profile scaling gives
\[
 m(b)=b^{(d-N)/2}m(1),\qquad m_*=m(a_*)>0.
\]
This function is continuous and strictly increasing on $[a_*,a^*]$.
Decrease $\rho>0$ so that
\begin{equation}\label{pc:choice-radius}
 \begin{gathered}
 \overline{B_{4\rho}(P_i)}\ \hbox{are pairwise disjoint},\qquad
 P_i\ \hbox{is the only critical point of }a
       \hbox{ in }\overline{B_{3\rho}(P_i)},\\
 |m(a(x))-m(a(P_i))|<\frac{\kappa}{2m}
       \quad\hbox{for }x\in\overline{B_\rho(P_i)},
       \quad 1\le i\le m.
 \end{gathered}
\end{equation}
Such a choice is possible by continuity and nondegeneracy.  In particular,
no center convergence to $P_i$ is part of the definition
\eqref{eq:positive-class}.

The next proposition proves convergence to the prescribed profiles.
We first exclude additional peaks, then locate the remaining centers.

\begin{proposition}\label{pc:class-compactness}
Let $J\subset\{1,\ldots,m\}$ be nonempty, let
$\varepsilon_n\downarrow0$, and let
$z_n\in\mathcal A_{\varepsilon_n,J}(\rho,\delta,\kappa)$.
After passing to a subsequence, there are centers
$\eta_{i,n}\to P_i$, $i\in J$, such that
\begin{equation}\label{pc:global-convergence}
 \begin{aligned}
 W_{J,n}(x)&=\sum_{i\in J}W_{a(P_i)}
                 \left(\frac{x-\eta_{i,n}}{\varepsilon_n}\right),\\
 \|z_n-W_{J,n}\|_\infty&\longrightarrow0.
 \end{aligned}
\end{equation}
Moreover,
\begin{equation}\label{pc:mass-quantization}
 \mathcal Q_{\varepsilon_n}(z_n)
       \longrightarrow\sum_{i\in J}m(a(P_i)).
\end{equation}
\end{proposition}

\begin{proof}
The definition \eqref{eq:positive-class} supplies a uniform mass bound,
so Lemma~\ref{pc:uniform-bound} applies.

For each $i\in J$, choose $x_{i,n}\in\overline{B_\rho(P_i)}$ where
$u_n+v_n\ge\delta$.  After passing to a subsequence,
$x_{i,n}\to Q_i\in\overline{B_\rho(P_i)}$.  These limits are distinct.
By Lemma~\ref{pc:autonomous-limits}, a common further subsequence has,
for every $i\in J$, a limit $W_{\widehat b_i}(\cdot-\zeta_i)$, where
$\widehat b_i=a(Q_i)$.  Set
$\eta_{i,n}=x_{i,n}+\varepsilon_n\zeta_i$.  Then
\[
 \eta_{i,n}\to Q_i,\qquad
 z_n(\eta_{i,n}+\varepsilon_n\,\cdot)
       \longrightarrow W_{\widehat b_i}
       \quad\hbox{in }C^2_{\rm loc}(\R^N)^2.
\]
At this stage $Q_i=P_i$ has not been proved or assumed.

We claim that
\begin{equation}\label{pc:outside-vanishing}
 \lim_{R\to\infty}\limsup_{n\to\infty}
 \sup_{x\notin\bigcup_{i\in J}B_{R\varepsilon_n}(\eta_{i,n})}
       (u_n(x)+v_n(x))=0.
\end{equation}
If this were false, there would be a further subsequence and points
$y_n$ such that
\[
 \min_{i\in J}\frac{|y_n-\eta_{i,n}|}{\varepsilon_n}
       \longrightarrow\infty,
 \qquad u_n(y_n)+v_n(y_n)\ge\delta_2>0.
\]
Take a subsequence with $a(y_n)\to\widehat b$.  Lemma~\ref{pc:autonomous-limits}
produces an additional positive profile of mass $m(\widehat b)\ge m_*$,
even if $|y_n|\to\infty$.  For each fixed scaled radius $R$, the balls
about $y_n$ and the $|J|$ retained centers are disjoint for large $n$.
Nonnegativity of the nonlinear density, local convergence, and then
$R\to\infty$ give
\[
 \liminf_n\mathcal Q_{\varepsilon_n}(z_n)
       \ge\sum_{i\in J}m(\widehat b_i)+m(\widehat b).
\]
The radius choice \eqref{pc:choice-radius} yields
\[
 \sum_{i\in J}m(\widehat b_i)+m(\widehat b)
 \ge\sum_{i\in J}m(a(P_i))+m_*
             -\frac{|J|\kappa}{2m}
 \ge\sum_{i\in J}m(a(P_i))+m_*-\frac\kappa2.
\]
This strictly exceeds the upper bound
$\sum_{i\in J}m(a(P_i))+m_*-\kappa$ in
\eqref{eq:positive-class}, a contradiction.  This proves
\eqref{pc:outside-vanishing}.

We next obtain an integrable uniform bound.  Choose $h>0$ as in
Lemma~\ref{pc:finite-energy}.  By \eqref{pc:outside-vanishing}, a fixed
large $R$ makes both components at most $h$ outside the balls
$B_{R\varepsilon_n}(\eta_{i,n})$ for all large $n$.  On that exterior,
\[
 (-\varepsilon_n^2\Delta+a_*/2)(u_n+v_n)\le0.
\]
If $C_1$ is the uniform amplitude bound and
$0<\sigma<\sqrt{a_*/2}$, the function
\[
 \Phi_n(x)=C_1e^{\sigma R}
       \sum_{i\in J}e^{-\sigma|x-\eta_{i,n}|/\varepsilon_n}
\]
is a supersolution there and dominates the solution on every core
sphere.  The exterior comparison of Lemma~\ref{lem:unique-exterior-comparison}
therefore applies.  Inside the core balls the same bound follows by
increasing the constant.  Rescaled interior estimates, exactly as in
Lemma~\ref{lem:unique-tails}, bound the gradients as well.  Consequently,
for some $C,\sigma>0$ independent of large $n$,
\begin{equation}\label{pc:uniform-tails}
 u_n+v_n+\varepsilon_n(|\nabla u_n|+|\nabla v_n|)
 \le C\sum_{i\in J}e^{-\sigma|x-\eta_{i,n}|/\varepsilon_n}
       \quad\hbox{on }\R^N.
\end{equation}

We now apply the local Pohozaev identity on $B_{2\rho}(P_i)$,
as in Lemma~\ref{ct:fixed-tube-compactness}. Its boundary stays a
fixed distance from all centers, so \eqref{pc:uniform-tails} gives
\[
 \varepsilon_n^{-N}\int_{B_{2\rho}(P_i)}
            \partial_l a(x)u_nv_n\,dx\longrightarrow0.
\]
Rescaling around $\eta_{i,n}$ gives the limit
$\partial_l a(Q_i)\int U_{\widehat b_i}V_{\widehat b_i}$.
The tail error after truncation at scaled radius $R$ is bounded by
\[
 C\int_{|y|>R}e^{-\sigma|y|}\,dy
       +C\varepsilon_n^{-N}e^{-c/\varepsilon_n};
\]
the second term comes from the other centers. Thus the product's
positive integral implies $\nabla a(Q_i)=0$. By
\eqref{pc:choice-radius}, $Q_i=P_i$ and $\widehat b_i=a(P_i)$.

Local convergence on each fixed scaled core, together with
\eqref{pc:outside-vanishing} and decay of the autonomous profiles, gives
$\|z_n-W_{J,n}\|_\infty\to0$. This proves the
admissibility condition \eqref{eq:unique-class}. Split each
nonlinear integral into the fixed scaled cores and their complement.
The core integrals converge to those of $W_{a(P_i)}$; their complements
are uniformly small by \eqref{pc:uniform-tails}.  This proves
\eqref{pc:mass-quantization} along the subsequence.

\end{proof}

\subsection{The exact count under the mass bound}

\begin{proof}[Proof of Theorem~\ref{thm:main-positive}]
Choose $\rho$ as in \eqref{pc:choice-radius}, and let
$\delta_0,\rho_0>0$ be the profile-neighborhood radii in
Theorem~\ref{thm:main-count}. For all sufficiently small
$\varepsilon$, simultaneously for every nonempty $J$,
\begin{equation}\label{pc:tube-inclusion}
 \mathcal A_{\varepsilon,J}(\rho,\delta,\kappa)
       \subset\mathcal S_\varepsilon\cap
                  \mathcal U_{\varepsilon,J}(\delta_0,\rho_0).
\end{equation}
Indeed, otherwise there would be a fixed nonempty $J$, a sequence
$\varepsilon_n\downarrow0$ and solutions in the mass class outside
the stated profile neighborhood. Proposition~\ref{pc:class-compactness}
gives, along a subsequence, centers tending to $P_i$ and uniform errors
tending to zero. For large $n$ the centers belong to
$B_{\rho_0}(P_i)$ and the errors are less than $\delta_0$, a contradiction.
Theorem~\ref{thm:main-count} now gives at most one solution in each
class, namely $z_{J,\varepsilon}$.

To show that this solution belongs to the class, recall that its centers
satisfy $\xi_{i,\varepsilon}=P_i+O(\varepsilon^2)$ and its correction
tends uniformly to zero. Thus
\[
 (u_{J,\varepsilon}+v_{J,\varepsilon})(\xi_{i,\varepsilon})
       =U_{a(P_i)}(0)+V_{a(P_i)}(0)+o(1)>\delta,
       \qquad i\in J,
\]
and these points lie in $B_\rho(P_i)$ for small $\varepsilon$.
Corollary~\ref{cor:potential-concentration-measures} gives
\[
 \mathcal Q_\varepsilon(z_{J,\varepsilon})
       =\sum_{i\in J}m(a(P_i))+o(1).
\]
Since $m_*-\kappa>0$, this also gives the mass bound in
\eqref{eq:positive-class}. Hence
$\mathcal A_{\varepsilon,J}=\{z_{J,\varepsilon}\}$.
There are finitely many nonempty subsets, so the parameter threshold
can be chosen uniformly. All remaining conclusions follow from
Theorems~\ref{thm:main-count} and \ref{thm:main-spectrum}, and
Section~\ref{sec:quantitative-concentration}.
\end{proof}
\par\medskip
\noindent\begin{minipage}{\textwidth}
\noindent\textbf{Code availability.} The Lean 4 formalization accompanying
this paper is available at \url{https://github.com/kfkeyan2026/LE-multipeaks-Lean}.
The repository contains the formalized statements and proofs, specifies the
versions of Lean and the required libraries, and provides instructions for
reproducing the verification.

\medskip
\noindent\textbf{Acknowledgments.}
The author acknowledges the assistance of ChatGPT (OpenAI) with the
preparation of Lean 4 code for formal verification 
and language polishing. The author takes full responsibility for the
mathematical content and the exposition
of this paper.
\end{minipage}

\end{document}